\documentclass[12pt,reqno]{amsart}
\usepackage{amssymb}
\usepackage{enumerate}
\usepackage[T1]{fontenc}
\usepackage[all]{xy}
\usepackage[usenames]{color}

\definecolor{darkgreen}{RGB}{65,165,118}
\definecolor{darkred}{RGB}{180,0,0}

\newcommand{\mynote}[1]{\noindent\textbf{[#1]}}

\newcommand{\sstar}{\textup{({\large$*$})}}
\newcommand{\sstars}{\textup{({\large$**$})}}
\newcommand{\pplus}{\textup{(\dag)}}
\newcommand{\ppluss}{\textup{(\ddag)}}

\newcommand{\8}{^\bullet}
\newcommand{\9}[1]{{}^{#1}\!}
\let\wh=\widehat

\newcommand{\dec}[1]{\widebar{#1}}
\newcommand{\ch}[1]{\widebar{#1}}
\newcommand{\xSFL}{(\dec{S},\dec\calf,\dec\call)}

\let\xto=\xrightarrow
\let\from=\leftarrow
\let\xfrom=\xleftarrow

\newcommand{\higherlim}[2]{\displaystyle\setbox1=\hbox{\rm lim}
	\setbox2=\hbox to \wd1{\leftarrowfill} \ht2=0pt \dp2=-1pt
	\setbox3=\hbox{$\scriptstyle{#1}$}
	\def\test{#1}\ifx\test\empty
	\mathop{\mathop{\vtop{\baselineskip=5pt\box1\box2}}}\nolimits^{#2}
	\else
	\ifdim\wd1<\wd3
	\mathop{\hphantom{^{#2}}\vtop{\baselineskip=5pt\box1\box2}^{#2}}_{#1}
	\else
	\mathop{\mathop{\vtop{\baselineskip=5pt\box1\box2}}_{#1}}%
	\nolimits^{#2}
	\fi\fi}

\let\oldcirc=\circ
\renewcommand{\circ}{\mathchoice
    {\mathbin{\scriptstyle\oldcirc}}{\mathbin{\scriptstyle\oldcirc}}
    {\mathbin{\scriptscriptstyle\oldcirc}}
    {\mathbin{\scriptscriptstyle\oldcirc}}}

\SelectTips{cm}{10} \UseTips   

\numberwithin{equation}{section}

\mathcode`\:="003A
\mathchardef\cdot="0201
\renewcommand{\*}{\mathop{\textup{\Large{$*$}}}\limits}

\renewenvironment{enumerate}[1][]
{\begin{enumerat}[#1]\setlength{\itemsep}{6pt}}{\end{enumerat}}

\newenvironment{enuma}{\begin{enumerate}[{\rm(a) }]}{\end{enumerate}}
\newenvironment{enumi}{\begin{enumerate}[{\rm(i) }]}{\end{enumerate}}

\renewenvironment{itemize}
{\begin{itemiz}\setlength{\itemsep}{6pt}\setlength{\itemindent}{-20pt}}
{\end{itemiz}}

\def\beq#1\eeq{\begin{equation*}#1\end{equation*}}
\def\beqq#1\eeqq{\begin{equation}#1\end{equation}}
\newcommand{\boldd}[1]{{\mathversion{bold}\textbf{#1}}}

\let\emptyset=\varnothing

\DeclareMathAlphabet\EuR{U}{eur}{m}{n}
\SetMathAlphabet\EuR{bold}{U}{eur}{b}{n}
\newcommand{\curs}{\EuR}

\renewcommand{\mod}{\textup{-}\curs{mod}}

\newcommand{\widebar}[1]{\overset{\mskip3mu\hrulefill\mskip3mu}{#1}
        \vphantom{#1}}

\renewcommand{\:}{\colon}

\newlength{\upto}\newlength{\dnto}
\newcommand{\I}[2]{\addtolength{\upto}{#1pt}\addtolength{\dnto}{#2pt}%
{\vrule height\upto depth\dnto width 0pt}}

\newcounter{let} 
\loop\stepcounter{let}
\expandafter\edef\csname cal\alph{let}\endcsname%
{\noexpand\mathcal{\Alph{let}}}
\ifnum\thelet<26\repeat

\newcommand{\tdef}[2][]{\expandafter\newcommand\csname#2\endcsname%
{#1\textup{#2}}}
\tdef{Iso}   \tdef{Aut}    \tdef{Out}    \tdef{Inn}    \tdef{Hom}
\tdef{End}   \tdef{Inj}    \tdef{map}    \tdef{Ker}    \tdef{Ob}
\tdef{Mor}   \tdef{Res}    \tdef{Spin}   \tdef{Id}     \tdef{Fr}
\tdef{rk}    \tdef{conj}   \tdef{incl}   \tdef{proj}   \tdef{diag} 
\tdef{trf}   \tdef{Sol}    \tdef{He}     \tdef{Sz}     \tdef{cj}
\tdef{free}  \tdef{ev}     \tdef{Map}    \tdef{Rep}    \tdef{pr}
\tdef{supp}  \tdef{res}    \tdef{hofibre}  
\tdef[_]{fc}  \tdef[_]{fn}  \tdef[_]{typ} \tdef[^]{op}   \tdef[^]{ab}

\newcommand{\fdef}[1]{\expandafter\newcommand\csname#1\endcsname%
{\mathfrak{#1}}}
\fdef{X}  \fdef{red}  \fdef{foc}  \fdef{hyp}  \fdef{Sub}

\newcommand{\bbdef}[1]{\expandafter\newcommand%
\csname#1\endcsname{\mathbb{#1}}}
\bbdef{C} \bbdef{F} \bbdef{R} \bbdef{Z} \bbdef{Q}

\newcommand{\gee}{\varepsilon}
\newcommand{\SFL}[1][]{(S#1,\calf#1,\call#1)}
\newcommand{\SFT}[1][]{(S#1,\calf#1,\calt#1)}

\newcommand{\gen}[1]{\langle{#1}\rangle}

\newcommand{\normal}{\trianglelefteq}
\let\nsg=\normal

\newcommand{\syl}[2]{\textup{Syl}_{#1}(#2)}
\newcommand{\sylp}[1]{\syl{p}{#1}}
\newcommand{\autf}{\Aut_{\calf}}

\newcommand{\outf}{\Out_{\calf}}

\newcommand{\homf}{\Hom_{\calf}}
\newcommand{\repf}{\Rep_{\calf}}
\newcommand{\isof}{\Iso_{\calf}}
\newcommand{\sminus}{\smallsetminus}
\newcommand{\defeq}{\overset{\textup{def}}{=}}
\newcommand{\zploc}{\Z_{(p)}}

\renewcommand{\Im}{\textup{Im}}

\renewcommand{\2}[2]{\underset{#1}{#2}}
\newcommand{\pcom}{{}^\wedge_p}

\newcommand{\4}[1]{\widehat{#1}}
\newcommand{\AUT}{\cala ut\typ}
\newcommand{\7}[1]{\boldsymbol{#1}}

\newcommand{\lra}{\longrightarrow}

\let\til=\widetilde

\newtheorem{Thm}{Theorem}[section]
\newtheorem{Prop}[Thm]{Proposition}
\newtheorem{Cor}[Thm]{Corollary}
\newtheorem{Lem}[Thm]{Lemma}

\newtheorem{Defi}[Thm]{Definition}

\newtheorem{Th}{Theorem}

\newcommand{\longleft}[1]{\;{\leftarrow%
\count255=0 \loop \mathrel{\mkern-6mu}%
    \relbar\advance\count255 by1\ifnum\count255<#1\repeat}\;}
\newcommand{\longright}[1]{\;{\count255=0 \loop \relbar\mathrel{\mkern-6mu}%
    \advance\count255 by1\ifnum\count255<#1\repeat\rightarrow}\;}
\newcommand{\Right}[2]{\overset{#2}{\longright{#1}}}
\newcommand{\RIGHT}[3]{\mathrel{\mathop{\kern0pt\longright{#1}}
    \limits^{#2}_{#3}}}
\newcommand{\Left}[2]{{\buildrel #2 \over {\longleft{#1}}}}
\newcommand{\LEFT}[3]{\mathrel{\mathop{\kern0pt\longleft{#1}}\limits^{#2}_{#3}}
}
\newcommand{\longleftright}[1]{\;{\leftarrow\mathrel{\mkern-6mu}%
    \count255=0\loop\relbar\mathrel{\mkern-6mu}%
    \advance\count255 by1\ifnum\count255<#1\repeat\rightarrow}\;}

\newcommand{\onto}[1]{\;{\count255=0 \loop \relbar\joinrel
    \advance\count255 by1
    \ifnum\count255<#1 \repeat \twoheadrightarrow}\;}
\newcommand{\Onto}[2]{\overset{#2}{\onto{#1}}}
\newcommand{\RLEFT}[3]{\mathrel{%
   \mathop{\vcenter{\baselineskip=0pt\hbox{$\kern0pt\longright{#1}$}%
   \hbox{$\kern0pt\longleft{#1}$}}}\limits^{#2}_{#3}}}

\theoremstyle{definition}
\newtheorem{Rmk}[Thm]{Remark}

\newcommand\5[1]{[\![#1]\!]}

\def\Top{\curs{Top}}
\def\hoTop{\curs{hoTop}}
\def\ho{\curs{ho}}
\def\hocolim{\mathop{\rm hocolim}}
\newcommand{\calhbullet}{\calf^\bullet}

\newcommand{\N}[2]{N_{#1}(#2)}

\title{An algebraic model for finite loop spaces}

\author{Carles Broto}
\address{Departament de Matemàtiques, Universitat Autònoma de
Barcelona, E--08193 Bellaterra, Spain}
\email{broto@mat.uab.es}
\thanks{C. Broto is partially supported by FEDER-MICINN grant MTM
2010-20692.}

\author{Ran Levi}
\address{Institute of Mathematics, University of Aberdeen,
Fraser Noble Building 138, Aberdeen AB24 3UE, U.K.}
\email{r.levi@abdn.ac.uk}
\thanks{R. Levi is partially supported by EPSRC grant EP/J014524/1.}

\author{Bob Oliver}
\address{Université Paris 13, Sorbonne Paris Cité, LAGA, UMR 7539 du CNRS, 
99, Av. J.-B. Clément, 93430 Villetaneuse, France.}
\email{bobol@math.univ-paris13.fr}
\thanks{B. Oliver is partially supported by UMR 7539 of the CNRS, and by 
project ANR BLAN08-2\_338236, HGRT}

\date{10 December 2012}



\begin{document}

\maketitle

\tableofcontents



The theory of $p$-local compact groups, developed in \cite{BLO3}, is 
designed to give a unified framework in which to study the $p$-local 
homotopy theory of classifying spaces of compact Lie groups and $p$-compact 
groups, as well as some other families of a similar nature. It also 
includes, and in many aspects generalizes, the earlier theory of $p$-local 
finite groups. 

Compact Lie groups, and the $p$-compact groups of Dwyer and Wilkerson 
\cite{DW}, are examples of loop spaces which are finite CW complexes, at 
least homologically. In this paper, we show that in fact, every 
homologically finite loop space $X$ gives rise to a $p$-local compact 
group, and thus enjoys all of the important properties (Sylow subgroups, 
Weyl groups, fusion relations) shared by these spaces.  We refer to 
\cite{BLO3} (where this was mentioned as an open question) for details. 

If $X$ is a finite loop space, then its group 
of components is finite, and its identity component is a connected finite 
loop space. As such, after $p$-completion in the sense of Bousfield and Kan 
\cite{BK}, the identity component is a $p$-compact group. This leads us to 
investigate a much more general question:  whether a space $B$ with a 
finite regular covering $E\to B$, such that $E$ is the classifying space of 
a $p$-local compact group, is itself the classifying space of a $p$-local 
compact group.

Before stating our main theorem, we briefly describe the objects of study. 
Let $\Z/p^\infty$ denote the union of all $\Z/p^n$ under the obvious 
inclusions. A \emph{discrete $p$-toral group} is a group $S$ containing a 
normal subgroup of the form $(\Z/p^\infty)^r$ ($r\ge0$) with $p$-power 
index. A (saturated) fusion system over $S$ is a category whose objects are 
the subgroups of $S$, and whose morphisms are monomorphisms of groups which 
are required to satisfy certain axioms. We refer the reader to Section 1 
(especially Definition \ref{sat.Frob.} and Corollary \ref{c:RS}) for 
details. 

Given a saturated fusion system $\calf$ over $S$, a centric linking system 
associated to $\calf$ is a category $\call$ whose objects are those 
subgroups of $S$ which are $\calf$-centric (Definition \ref{F-centric}), 
and whose morphism sets  are, in an appropriate sense, extensions of the 
corresponding morphism sets between the same objects in $\calf$. This extra 
structure allows us to associate a ``classifying 
space'' to the fusion system in question. Thus a $p$-local compact group is 
a triple $\SFL$, where $S$ is a discrete $p$-toral group, $\calf$ is a 
saturated fusion system over $S$, and $\call$ is a centric linking system 
associated to $\calf$. The classifying space of a $p$-local compact group 
is the $p$-completion of the geometric realisation of its linking system.

Our main theorem is the following.

\begin{Th} \label{ThA}
Assume that $f\:X\Right2{}Y$ is a finite regular covering space with covering 
group $G$, where $X$ is the classifying space of a $p$-local compact group 
$(\dec{S},\dec\calf,\dec\call)$.  Then $Y\pcom$ is the classifying space of a 
$p$-local compact group $\SFL$. 
\end{Th}

The proof of Theorem \ref{ThA} (restated as Theorem \ref{ThmA}) goes 
roughly as follows.  Algebraic extensions of $p$-local finite groups were 
described and classified in \cite{link-ext}, and this construction is 
generalized here (in Section 5) to the $p$-local compact case.  We then 
prove an equivalence between the data needed to define such an extension of 
$\dec\call$ by a finite group $G$, and that needed to define fibre bundles 
over $BG$ with fibre $|\dec\call|$.  Finally, in Section 7, this 
equivalence is used to prove that each $|\dec\call|\simeq X\Right3{}Y$ as 
in the hypotheses of the theorem can be realized by an extension of 
$p$-local compact groups as constructed in Section 5.

Special cases of Theorem \ref{ThA} (those where $G$ is a $p$-group or has 
order prime to $p$) are proven in \cite{BCGLO2}. But the proofs there 
cannot be applied to extensions by arbitrary finite groups, not even in the 
$p$-local finite case. 

Let $p$ be a prime. A space $Y$ is $p$-good if its $p$-completion $Y\pcom$ 
is $p$-complete \cite{BK}. A space $Y$ is $\F_p$-finite if $H_*(Y,\F_p)$ is 
a finite dimensional $\F_p$-vector space. In \cite{DW}, Dwyer and Wilkerson 
define a $p$-compact group to be an $\F_p$-finite loop space $X$ whose 
classifying space $BX$ is $p$-complete. In \cite{BLO3}, we show that every 
$p$-compact group $X$ gives rise to a $p$-local compact group $\SFL$ whose 
classifying space has the homotopy type of $BX$. 

\begin{Th} \label{ThB}
Let $X$ be any path connected space. Then for each prime $p$ such that 
$\Omega X$ is $\F_p$-finite, the space  $X\pcom$ has the homotopy type of 
the classifying space of a $p$-local compact group.  In particular, this 
holds for each prime $p$ if $\Omega X$ has the homotopy type of a finite 
complex.
\end{Th}

\begin{proof}
Fix a prime $p$ such that $\Omega X$ is $\F_p$-finite.  In particular, 
$H_0(\Omega X,\F_p)$ is finite, so $\pi_1(X)$ is a finite group, and hence 
$X$ is $p$-good by \cite[\S\,VII.5]{BK}. Set $\pi = \pi_1(X)$, let 
$\widetilde{X}$ be the universal cover of $X$, and 
consider the following fibration:
	\[\widetilde{X}\to X\to B\pi.\]
Then $\Omega\widetilde{X}$ is a connected component of $\Omega X$, and 
hence is also $\F_p$-finite. 
Applying fibrewise $p$-completion, we obtain a fibration 
	\[\widetilde{X}\pcom\to \widebar{X}\to B\pi,\]
where $\widetilde{X}\pcom$ is $p$-complete and $\Omega\widetilde{X}\pcom$ 
is $\F_p$-finite. Thus $\widetilde{X}\pcom$ is the classifying space of a 
$p$-compact group. By \cite[Theorem 10.7]{BLO3}, $\widetilde{X}\pcom$ is the 
classifying space of a $p$-local compact group.  So by 
Theorem \ref{ThA}, $\widebar{X}\pcom \simeq X\pcom$ is the classifying 
space of a $p$-local compact group.

In particular, if $X$ is a path connected space such that  $\Omega X$ has 
the homotopy type of a finite CW complex, then this holds for $X$ at all 
primes. 
\end{proof}

Since compact Lie groups are finite loop spaces, the above result applies 
to compact Lie groups, and thus provides an alternative argument to that in 
\cite[\S9]{BLO3}. 

The proof of Theorem \ref{ThA} required several results already known for 
$p$-local finite groups to be generalised to the $p$-local compact case. 
For example, Section 4 contains a proof that the family of $\calf$-centric 
$\calf$-radical subgroups determines saturation (compare \cite[Theorem 
A]{BCGLO1}), and Section 5 describes how to construct extensions of 
$p$-local compact groups by finite groups.  These results should be useful 
in other future developments of the theory of $p$-local compact groups.

The paper is organized as follows. We start in Section 1 with the basic 
definitions and general background on $p$-local compact groups. Sections 
2--5 consist mostly of generalizations to the $p$-local compact case of 
certain constructions and results which are already known for $p$-local 
finite groups. The equivalence between extensions of $\dec\call$ by $G$ and 
bundles over $BG$ with fibre $|\dec\call|$ is shown in Section 6.  All of 
this is then combined in Section 7 to prove Theorem \ref{ThA}. We end the 
paper with an appendix where we collect necessary results on transporter 
systems over discrete $p$-toral groups.

\textbf{Notation: }  When $G$ and $H$ are groups, $H<G$ always means that 
$H$ is a proper subgroup of $G$.  Whenever $F\:\calc\Right2{}\cald$ is a 
functor and $c,c'\in\Ob(\calc)$, we write $F_{c,c'}$ for the restriction of 
$F$ to $\Mor_{\calc}(c,c')$, and write $F_c=F_{c,c}$.  

The authors wish to acknowledge support from the Universitat Aut\`onoma de 
Barcelona, the Institute of Mathematics in Aberdeen and the University of 
Copenhagen for allowing the authors to meet and providing support at 
various times while this project was developing.

\bigskip


\section{Background on fusion and linking systems over discrete $p$-toral 
groups}

In this section, we collect  the definitions and  some basic facts on 
$p$-local compact groups which will be useful throughout the paper. We 
refer to \cite{BLO3} for more details on many of the results described 
here.

\begin{Defi}  
A \emph{discrete $p$-torus} is a group which is isomorphic to 
$(\Z/p^\infty)^n$ for some finite $n$.  
A \emph{discrete $p$-toral group} is a group $P$ with a 
normal subgroup $P_0\nsg{}P$ such that $P_0$ is a discrete $p$-torus and 
$P/P_0$ is a finite $p$-group. The 
subgroup $P_0$ will be called the \emph{identity component} of $P$, and $P$ 
will be called \emph{connected} if $P=P_0$. Set $\pi_0(P)=P/P_0$, 
and set $\rk(P)=n$ (the \emph{rank} of $P$) if $P_0\cong(\Z/p^\infty)^n$. 
The \emph{order} of $P$ is the 
pair $|P|\defeq (\rk(P), |\pi_0(P)|)$, regarded as an element of 
$\mathbb{N}^2$ ordered lexicographically. 
\end{Defi}

Let $\Sub(S)$ be the set of all subgroups of a (discrete $p$-toral) group 
$S$. For any group $G$ and any $H,K\le G$, 
$\Hom_G(H,K)\subseteq\Hom(H,K)$ denotes the set of homomorphisms induced by 
conjugation in $G$.

\begin{Defi} \label{Frob.def.}
A \emph{fusion system} $\calf$ over a discrete $p$-toral group $S$ is a 
category with $\Ob(\calf)=\Sub(S)$, whose morphism sets 
$\homf(P,Q)$ satisfy the following conditions:
\begin{enuma}
\item $\Hom_S(P,Q)\subseteq\homf(P,Q)\subseteq\Inj(P,Q)$ for all $P,Q\in\Sub(S)$.

\item Every morphism in $\calf$ factors as an isomorphism in $\calf$
followed by an inclusion.
\end{enuma}
\end{Defi}

With motivation from group theory, we make the following definition.

\begin{Defi}\label{F-conj}
Let $\calf$ be a fusion system over a discrete $p$-toral group $S$. Two subgroups $P ,Q\in \Sub(S)$ are \emph{$\calf$-conjugate} if they are
isomorphic as objects of the category $\calf$. Let $P^{\calf}$ denote the
set of all subgroups of $S$ which are $\calf$-conjugate to $P$. 
\end{Defi}

We are now ready to recall the definition of saturation of a fusion system. 

\begin{Defi} \label{sat.Frob.}
Let $\calf$ be a fusion system over a discrete $p$-toral group $S$.
\begin{itemize}
\item A subgroup $P\le{}S$ is \emph{fully centralized in $\calf$} if
$|C_S(P)|\ge|C_S(Q)|$ for all $Q\in P^\calf$.
\item A subgroup $P\le{}S$ is \emph{fully normalized in $\calf$} if
$|N_S(P)|\ge|N_S(Q)|$ for all $Q\in P^\calf$.
\item $\calf$ is a \emph{saturated fusion system} if the following
three conditions hold:
\begin{enumerate}[\rm(I) ]
\item For each $P\le{}S$ which is fully normalized in $\calf$, $P$ is fully
centralized in $\calf$, $\outf(P)$ is finite, and 
$\Out_S(P)\in\sylp{\outf(P)}$.  

\item If $P\le{}S$ and $\varphi\in\homf(P,S)$ are such that $\varphi(P)$ is
fully centralized, and if we set
    $$ N_\varphi = \{ g\in{}N_S(P) \,|\, \9\varphi c_g \in
    \Aut_S(\varphi(P)) \}, $$
then there is $\widebar{\varphi}\in\homf(N_\varphi,S)$ such that
$\widebar{\varphi}|_P=\varphi$.
\item If $P_1\le{}P_2\le{}P_3\le\cdots$ is an increasing sequence of 
subgroups of $S$, with $P_\infty=\bigcup_{n=1}^\infty{}P_n$, and if 
$\varphi\in\Hom(P_\infty,S)$ is any homomorphism such that 
$\varphi|_{P_n}\in\homf(P_n,S)$ for all $n$, then 
$\varphi\in\homf(P_\infty,S)$.
\end{enumerate}
\end{itemize}
\end{Defi}

We next define certain classes of subgroups which play an important 
role in generating the morphisms in a fusion system.

\begin{Defi}\label{F-centric}
Let $\calf$ be a fusion system over a discrete $p$-toral group $S$. A 
subgroup $P\in\Sub(S)$ is \emph{$\calf$-centric} if $C_S(Q)=Z(Q)$ for all 
$Q\in P^\calf$. A subgroup $P\in \Sub(S)$ is \emph{$\calf$-radical} if 
$\outf(P) \defeq \autf(P)/\Inn(P)$ is $p$-reduced; i.e., contains no 
nontrivial normal $p$-subgroup. Let $\calf^c$ denote the full 
subcategory of $\calf$ whose objects are the $\calf$-centric subgroups 
of $S$. 
\end{Defi}

In this paper it will be convenient to also use a different but 
equivalent definition of saturation, based on that due to Roberts 
and Shpectorov \cite{RS} in the finite case.  We recall 
their definitions.

\begin{Defi} \label{d:sat}
 Let $\calf$ be a fusion system over a $p$-group $S$.
 \begin{itemize}
\item A subgroup $P \leq S$ is \emph{fully automized} in $\calf$ if
the index of $\Aut_S(P)$ in $\autf(P)$ is finite and prime to $p$.

\item  A subgroup $P \leq S$ is \emph{receptive} in $\calf$ if it has the 
following property: for each $Q \leq S$ and each $\varphi\in \Iso_{\calf} 
(Q, P)$, if we set $N_\varphi = \{g \in N_S (Q) \;|\; \9{\varphi}c_g\in 
\Aut_S (P)\}$, then there is $\bar\varphi\in Hom_{\calf} (N_\varphi , S)$ 
such that $\bar\varphi|_P = \varphi$. 

\end{itemize}
\end{Defi}

In this terminology, axioms (I) and (II) in Definition \ref{sat.Frob.} 
say that each fully normalized subgroup is fully centralized and fully 
automized, and each fully centralized subgroup is receptive.

By definition, if $P$ is fully automized in $\calf$, then $\autf(P)$ 
is an extension of a discrete $p$-torus by a finite group.  Hence it does 
contain maximal discrete $p$-toral subgroups, unique up to conjugation, 
which we regard as its Sylow $p$-subgroups.  As usual, we let 
$\sylp{\autf(P)}$ denote the set of its Sylow $p$-subgroups.

The next lemma describes the relation between these concepts and those 
already defined.

\begin{Lem}\label{l:RS}
The following hold for any fusion system $\calf$ over a discrete $p$-toral 
group $S$.
\begin{enuma}
\item \label{RS-1} Every receptive subgroup of $S$ is fully centralized.

\item \label{RS-2} If $P\le S$ is fully automized and
receptive in $\calf$, then it is fully normalized.  

\item \label{RS-3} If $P\le{}S$ is fully automized and receptive in 
$\calf$, and $Q\in{}P^\calf$, then there is a morphism 
$\varphi\in\Hom_{\calf}(N_S(Q),N_S(P))$ such that 
$\varphi(Q)=P$. 

\end{enuma}
\end{Lem}

\begin{proof} The proofs are identical to those given in \cite{RS} and 
\cite[Lemma I.2.6]{AKO} in the finite case.  A generalization of these 
statements will be proven in Lemma \ref{l:RS-K} below.
\end{proof}

The following is an immediate consequence of Lemma \ref{l:RS}. 

\begin{Cor} \label{c:RS}
A fusion system $\calf$ over a discrete $p$-toral group $S$ is saturated 
if and only if 
\begin{itemize} 
\item each subgroup of $S$ is $\calf$-conjugate to one which is fully 
automized and receptive in $\calf$; and
\item axiom (III) holds for $\calf$:  if $P_1<P_2<\cdots$ are subgroups of 
$S$, $P=\bigcup_{i=1}^\infty{}P_i$, and $\varphi\in\Hom(P,S)$ is such that 
$\varphi|_{P_i}\in\homf(P_i,S)$ for each $i$, then $\varphi\in\homf(P,S)$.
\end{itemize}
\end{Cor}

When $G$ is a finite group and $\calh$ is a set of subgroups of $G$, 
$\calt_\calh(G)$ denotes the \emph{$\calh$-transporter category} for 
$G$:  $\Ob(\calt_\calh(G))=\calh$, and for $H,K\in\calh$,
	\[ \Mor_{\calt_\calh(G)}(H,K) = 
	\bigl\{ g\in{}G \,\big|\, \9g H\le K \bigr\}\,. \]

\begin{Defi}\label{L-cat}
Let $\calf$ be a fusion system over a discrete p-toral group $S$. Let 
$\calh$ be a set of subgroups of $S$ which is closed under 
$\calf$-conjugacy and overgroups,  and includes all subgroups which are 
$\calf$-centric and $\calf$-radical.  An $\calh$-linking system associated 
to $\calf$ is a category $\call$ with object set $\calh$, 
together with a pair of functors 
	\[ \calt_{\calh}(S) \Right5{\delta} \call
	\Right5{\pi} \calf\,, \]
such that each object  is
isomorphic (in $\call$) to one which is fully centralized in $\calf$, 
and such that  the following conditions are satisfied:
\begin{enumerate}[\rm(A) ]
\item The functor
$\delta$ is the identity on objects, and $\pi$ is the inclusion on 
objects.  For each $P,Q\in\calh$ such that $P$ is fully centralized 
in $\calf$, $C_S(P)$ acts freely on $\Mor_{\call}(P,Q)$ via $\delta_{P}$ 
and right composition, and $\pi_{P,Q}$ induces a bijection 
	\[ \Mor_\call(P,Q)/C_S(P) \Right5{\cong} \homf(P,Q) ~. \]

\item For each $P,Q\in\calh$ and each $g\in{}N_S(P,Q)$, $\pi_{P,Q}$ 
sends $\delta_{P,Q}(g)\in\Mor_\call(P,Q)$ to $c_g\in\homf(P,Q)$. 

\item For all $\psi\in\Mor_{\call}(P,Q)$ and all $g\in P$, 
$\psi\circ\delta_P(g)=\delta_Q(\pi(\psi)(g))\circ\psi$.
\end{enumerate}
A \emph{centric linking system} is an $\Ob(\calf^c)$-linking system; i.e., 
a linking system whose objects are the $\calf$-centric subgroups of $S$.
\end{Defi}

When $P\le Q$, we set $\iota_{P,Q}=\delta_{P,Q}(1)$. 
The morphisms $\iota_{P,Q}$ are regarded as the inclusions in $\call$.

\begin{Defi}\label{triple}
A \emph{linking triple} is a triple
of the form $\SFL$, where $S$ is a discrete $p$-toral 
group, $\calf$ is a saturated fusion system over $S$, and $\call$
is an $\calh$-linking system for some family $\calh$. 
A \emph{$p$-local compact group} is a linking triple 
where $\call$ is a centric linking system.
\end{Defi}

\begin{Defi} \label{d:H-sat}
Let $\calf$ be a fusion system over a discrete $p$-toral group $S$, and let  
$\calh\subseteq\Sub(S)$  be a family of subgroups. Then $\calh$ is 
\emph{closed} in $\Sub(S)$ if for each increasing sequence 
$P_1<P_2<P_3<\cdots$ of subgroups in $\calh$, $\bigcup_{i=1}^\infty{}P_i$ 
is also in $\calh$.

Let $\calh\subseteq\Sub(S)$ be a closed family.  Then 

\begin{enumerate}[\rm (a)]
\item $\calf$ is \emph{$\calh$-closed} if for each sequence $P_1<P_2<\cdots$  
in $\calh$ with $P=\bigcup_{i=1}^\infty{}P_i$, and each homomorphism 
$\varphi\in\Hom(P,S)$ such that 
$\varphi|_{P_i}\in\homf(P_i,S)$ for each $i$, $\varphi\in\homf(P,S)$.

\item $\calf$ is \emph{$\calh$-generated} if every morphism in $\calf$ is a 
composite of restrictions of morphisms in $\calf$ between subgroups in 
$\calh$.

\item $\calf$ is  \emph{$\calh$-saturated} if it is $\calh$-closed and 
if every subgroup of $\calh$ is 
$\calf$-conjugate to a subgroup which is fully automized and receptive.

\end{enumerate}
\end{Defi}

The following two results, both generalizations to discrete $p$-toral 
groups of well known properties of $p$-groups, will be useful.

\begin{Lem}[{\cite[Lemma 1.8]{BLO3}}] \label{N_Q(P)>P}
If $P<Q$ are discrete $p$-toral groups, then $P<N_Q(P)$.
\end{Lem}

\begin{Lem} \label{Aut(PrelQ)}
Let $Q\nsg P$ be discrete $p$-toral groups, where $|P/Q|<\infty$.  Then the 
group
	\[ \bigl\{ \alpha\in\Aut(P) \,\big|\, \alpha|_Q=\Id,~ [\alpha,P]\le 
	Q \bigr\} \]
is discrete $p$-toral.
\end{Lem}

\begin{proof} Let $A$ denote this group of automorphisms.  For each 
$\alpha\in{}A$ and each $g\in{}P$, $\alpha(g)=g\chi(g)$ for some 
$\chi(g)\in{}Q$ since $[\alpha,P]\le Q$, $c_g=c_{\alpha(g)}\in\Aut(Q)$ 
since $\alpha|_Q=\Id$, and hence $\chi(g)\in{}Z(Q)$.  Also, for $a\in{}Q$, 
$\alpha(ga)=g\chi(g)a=(ga)\chi(g)$ since $\alpha(a)=a$ and $[a,\chi(g)]=1$, 
so $\chi(ga)=\chi(g)$. Thus $\alpha$ is determined by the map 
$\chi\:P/Q\Right2{}Z(Q)$.  The resulting injection 
$A\Right2{}\Map(P/Q,Z(Q))$ is a homomorphism, and so $A$ is discrete 
$p$-toral since $\Map(P/Q,Z(Q))\cong Z(Q)^{|P/Q|}$ is.  
\end{proof}

\bigskip


\section{Normalizer fusion subsystems}

Let $\calf$ be a  fusion system over a discrete $p$-toral group $S$. For 
each $Q\leq S$ and each $K\leq\Aut(Q)$, define
	\begin{align*} 
	\autf^K(Q)&=K\cap\autf(Q)\\
	\Aut_S^K(Q)&=K\cap\Aut_S(Q) \\
	N_S^K(Q)&=\{\,x\in N_S(Q)\,|\, c_x\in K\,\} \quad \textup{(the 
	\emph{$K$-normalizer} of $Q$ in $S$).}
	\end{align*}

\begin{Defi} \label{K-norm.}
Let $\calf$ be a  fusion system over a discrete $p$-toral group $S$.  
Fix a subgroup $Q\leq S$ and a group of automorphisms $K\leq\Aut(Q)$.  
\begin{itemize} 
\item $Q$ is \emph{fully $K$-automized in $\calf$} if 
$\Aut_S^K(Q)\in\sylp{\autf^K(Q)}$.

\item $Q$ is \emph{fully $K$-normalized in $\calf$} if for each 
$\varphi\in\Hom_{\calf}(Q,S)$, 
	$$|N_S^K(Q)|\geq|N_S^{\9\varphi K}(\varphi(Q))|\,,$$ 
where ${\9\varphi K}=\bigl\{\varphi\alpha\varphi^{-1} 
\,\big|\, \alpha\in K\,\bigr\} \leq \Aut(\varphi(Q))$.

\item $N_\calf^K(Q)\subseteq\calf$ is the 
fusion system over $N_S^K(Q)$ where for $P,R\le{}N_S^K(Q)$, 
	\begin{multline*}  
	\Hom_{N_\calf^K(Q)}(P,R) = \bigl\{ \varphi\in\homf(P,R) \,\big|\,
	\textup{$\exists$ $\widebar\varphi\in\homf(PQ,RQ)$} \\
	\textup{with $\widebar\varphi|_P=\varphi$, 
	$\widebar{\varphi}(Q)=Q$, and $\widebar\varphi|_Q\in K$} \bigr\}. 
	\end{multline*}

\end{itemize}
\end{Defi}

As special cases of the above definition, set 
$N_\calf(Q)=N_\calf^{\Aut(Q)}(Q)$ and $C_\calf(Q)=N_\calf^{\{1\}}(Q)$:  the 
normalizer and centralizer fusion systems, respectively, of $Q$.  

The next lemma is a generalization of results in \cite{RS} to fusion 
systems over discrete $p$-toral groups.

\begin{Lem}\label{l:RS-K}
Let $\calf$ be a fusion system over a discrete $p$-toral group $S$, and let 
$\calp$ be an $\calf$-conjugacy class of subgroups of $S$.  Assume either that 
$\calf$ is saturated, or (more generally) that $\calp$ contains a subgroup 
which is fully automized and receptive in $\calf$. Then the following hold for 
each $P\in\calp$ and each $K\le\Aut(P)$.
\begin{enuma}

\item     \label{RS-lem1-3a} The subgroup $P$ is fully centralized if and 
only if it is receptive. 

\item \label{fullyKnorm} The subgroup $P$ is fully $K$-normalized if and 
only if it is fully $K$-automized and receptive. In this case, for each 
$Q\in\calp$ and each $\varphi\in\isof(Q,P)$, there are $\chi\in\autf^K(P)$ 
and $\widebar\varphi\in\homf(N^{K^\varphi}_S(Q)\cdot Q,N_S^K(P)\cdot P)$ 
such that $\widebar{\varphi}|_Q=\chi\circ\varphi$. 

\end{enuma}
\end{Lem}

\begin{proof} \textbf{(a) } For fusion systems over finite $p$-groups, this 
is shown in \cite[Propositions 3.7 \& 4.6]{RS} and in \cite[Lemma 
I.2.6(c)]{AKO}.  Those proofs carry over unchanged to the discrete $p$-toral 
case.

\smallskip

\noindent\textbf{(b) }  Assume $P$ is fully $K$-automized and receptive in 
$\calf$.  Fix $Q\in\calp$ and $\varphi\in\isof(Q,P)$.  Then 
$\9\varphi\Aut_S^{K^\varphi}(Q)$ is a discrete $p$-toral subgroup of 
$\autf^K(P)$, and since $P$ is fully $K$-automized, there is 
$\chi\in\autf^K(P)$ such that 
$\9{\chi\varphi}\Aut_S^{K^\varphi}(Q)\le\Aut_S^K(P)$.  Note that 
$K^\varphi=K^{\chi\varphi}$ since $\chi\in{}K$. 

Set $\psi={\chi\circ\varphi}$.  Then 
$\9\psi\Aut_S^{K^\psi}(Q)\le\Aut_S^K(P)$, and since $P$ is receptive, there 
is a morphism $\widebar{\psi}\in\homf(N_S^{K^\psi}(Q)\cdot Q,N_S^K(P)\cdot 
P)$ which extends $\psi$. In particular, 
$|N_S^{K^\varphi}(Q)|=|N_S^{K^\psi}(Q)|\le |N_S^K(P)|$.  Since $Q$ and 
$\varphi$ were arbitrary, this proves that $P$ is fully 
$K$-normalized in $\calf$.  

Now assume $P$ is fully $K$-normalized.  Fix $R\in\calp$ and   
$\psi\in\isof(P,R)$ such that $R$ is fully automized and receptive in 
$\calf$.
Fix $T\in\sylp{\autf^{\9\psi{}K}(R)}$ such that 
$T\ge\9\psi\Aut_S^K(P)$.  Since $R$ is fully automized, there is 
$\alpha\in\autf(R)$ such that $\9\alpha{}T\le\Aut_S(R)$.  
Set $\varphi={\alpha\circ\psi}$. Then
$\9\varphi\Aut_S^K(P)\le \9\alpha{}T\le\Aut_S(R)$, and since $R$ is 
receptive, $\varphi$ extends to 
$\widebar{\varphi}\in\homf(N_S^K(P)\cdot P,N_S(R))$, where 
$\Im(\widebar{\varphi})\le N_S^{\9\varphi K}(R)$.

Consider the following commutative diagram 
	\beq 
	\vcenter{\xymatrix{ 
	1 \ar[r] & C_S(P) \ar[r]\ar[d] & N_S^K(P) 
	\ar[r]\ar[d]_{\widebar\varphi}  & \Aut_S^K(P) \ar[d]^{c_\varphi}
	\ar[r] & 1 \\ 
	1 \ar[r] & C_S(R) \ar[r]    & N_S^{{}^\varphi\!K}(R) \ar[r]   & 
	\Aut_S^{{}^\varphi\!K}(R) \ar[r] & 1
	}} \eeq
where the rows are exact and all vertical maps are monomorphisms. 
Since $P$ is fully $K$-normalized, 
$\widebar\varphi(N_S^K(P))=N_S^{{}^\varphi\!K}(R)$, and hence 
$\widebar\varphi(C_S(P))=C_S(R)$ and $\9\varphi\Aut_S^K(P) =  
\Aut_S^{\9\varphi K}(R)$.  Thus $P$ is receptive by (a) (and since $R$ is 
receptive).  

Now, $T\in\sylp{\autf^{\9\psi K}(R)}$ by assumption, so 
$\9\alpha{}T\in\sylp{\autf^{\9\varphi K}(R)}$.  Also, 
$\9\alpha{}T\le\Aut_S(R)$, and hence 
$\9\alpha{}T=\Aut_S^{\9\varphi{}K}(R)$. Thus 
$\Aut_S^K(P)\in\sylp{\autf^K(P)}$, so $P$ is fully $K$-automized, and this 
finishes the proof of (b).
\end{proof}

The main result in this section is that normalizer fusion subsystems over 
discrete $p$-toral groups are saturated. The proof given here is modelled 
on that of \cite[Theorem I.5.5]{AKO}.

\newcommand{\xx}{_{+}}

\begin{Thm} \label{N_F(Q)-sat}
Fix a saturated fusion system $\calf$ over a discrete $p$-toral group $S$.  
Assume $Q\le{}S$ and $K\le\Aut(Q)$ are such that $Q$ is fully 
$K$-normalized in $\calf$.  Then $N_\calf^K(Q)$ is a saturated fusion 
system over $N_S^K(Q)$.  
\end{Thm}

\begin{proof}  Set $S\xx=N_S^K(Q)$ and $\calf\xx=N_\calf^K(Q)$ for short.  
For each $P\le S\xx$, set 
	\[ K_P = \bigl\{\alpha\in\Aut(PQ) \,\big|\, \alpha(P)=P,\ 
	\alpha(Q)=Q,\ \alpha|_Q\in{}K \bigr\}. \]
We need to show the following statements.
\begin{enuma}  
\item Each subgroup of $S\xx$ is $\calf\xx$-conjugate to a subgroup $P$ such 
that $PQ$ is fully $K_P$-normalized in $\calf$.
\item If $P\le{}S\xx$ and $PQ$ is fully $K_P$-normalized in $\calf$, then $P$ 
is fully automized in $\calf\xx$.
\item If $P\le{}S\xx$ and $PQ$ is fully $K_P$-normalized in $\calf$, then $P$ 
is receptive in $\calf\xx$.
\item Axiom (III) holds for $N_\calf^K(Q)$.  
\end{enuma}
The theorem will then follow immediately from Corollary \ref{c:RS}.

The proofs of (a), (b), and (c) are identical to the corresponding proofs 
in \cite[Theorem I.5.5]{AKO} (which is stated for fusion systems over 
finite $p$-groups).  It remains to prove (d).

Fix subgroups $P_1\le P_2\le P_3\le\cdots\le S\xx$, 
and set $P=\bigcup_{i=1}^\infty P_i$.  Assume $\varphi\in\Hom(P,S\xx)$ is 
such that $\varphi_i\defeq\varphi|_{P_i}\in\Hom_{\calf\xx}(P_i,S\xx)$ for 
each $i$.  

For each $i$, set 
	\[ X_i = \bigl\{ \psi\in\homf(P_iQ,S) \,\big|\, 
	\psi|_{P_i}=\varphi_i,~ \psi(Q)=Q,~ \psi|_Q\in K \bigr\} 
	\ne\emptyset \,, \]
and let $\widebar{X}_i$ be the image of $X_i$ in $\repf(P_iQ,S)$ 
(equivalently, the set of $C_S(P_i)$-conjugacy classes of morphisms in 
$X_i$).  Since $\repf(P_iQ,S)$ is finite by \cite[Lemma 2.5]{BLO3}, so is 
$\widebar{X}_i$.  There are natural restriction maps 
$\widebar{X}_i\Right2{}\widebar{X}_{i-1}$, and since the sets are finite 
and nonempty, the inverse limit is nonempty.  Fix an element 
$\bigl([\psi_i]\bigr)_{i=1}^\infty$ in the inverse limit.  Thus 
$\psi_i\in\homf(P_iQ,S)$, $\psi_i|_{P_i}=\varphi_i$, $\psi_i|_Q\in{}K$, and 
$c_{g_{i-1}}\circ\psi_i|_{P_{i-1}Q}=\psi_{i-1}$ for some 
$g_{i-1}\in{}C_S(P_{i-1})$.  

Since $S$ is artinian (cf. \cite[Proposition 1.2]{BLO3}), there is $N$ such 
that $C_S(P_i)=C_S(P)$ for each $i\ge N$.  For each $i\ge N$, set 
$\psi'_i=c_{g_N}\circ c_{g_{N+1}}\circ\cdots\circ c_{g_{i-1}}\circ\psi_i$.  Then 
$\psi'_i|_{P_{i-1}}=\psi'_{i-1}$, and $\psi'_i|_{P_i}=\varphi_i$. Set 
$\psi'=\bigcup_{i=1}^\infty\psi'_i$. Then $\psi'\in\homf(PQ,S)$ by axiom 
(III) for $\calf$, $\psi'|_P=\varphi$, $\psi'|_Q=\psi_N|_Q\in{}K$, and so 
$\varphi\in\Hom_{\calf\xx}(P,S)$.  
\end{proof}

The following is one easy application of Theorem \ref{N_F(Q)-sat}.

\begin{Lem} \label{l:centric}
Let $\calf$ be a saturated fusion system over a discrete $p$-toral group 
$S$.  Assume $Q\nsg{}P\le{}S$, where $Q$ is $\calf$-centric. Let 
$\varphi,\varphi'\in\homf(P,S)$ be such that $\varphi|_Q=\varphi'|_Q$.  
Then there is $x\in{}Z(Q)$ such that $\varphi'=\varphi\circ c_x$.
\end{Lem}

\begin{proof}  Since $\varphi\circ c_g=c_{\varphi(g)}\circ\varphi$ for each 
$g\in Q$, it suffices to show that $\varphi'=c_y\circ\varphi$ for some 
$y\in Z(\varphi(Q))$.  Upon replacing $P$ by $\varphi'(P)$, $Q$ by 
$\varphi(Q)=\varphi'(Q)$, and $\varphi$ by $\varphi\circ(\varphi')^{-1}$, 
we can assume that $\varphi'=\incl_P^S$ and $\varphi|_Q=\Id_Q$. We must 
show that $\varphi=c_x$ for some $x\in Z(Q)$.  

Set $K=\Aut_P(Q)$.  Since $Q$ is $\calf$-centric, it is fully centralized. 
Since $\autf^K(Q)=\Aut_S^K(Q)=K$, $Q$ is fully $K$-automized, and hence 
fully $K$-normalized by Lemma \ref{l:RS-K}(a,b).  Hence by Theorem 
\ref{N_F(Q)-sat}, the normalizer subsystem $N_\calf^K(Q)$ over 
$N_S^K(Q)=P{\cdot}C_S(Q)=P$ is saturated.  Also, since $\varphi|_Q=\Id$, 
$\Aut_{\varphi(P)}(Q)=\Aut_P(Q)=K$.  Thus $\varphi(P)\le N_S^K(Q)$, and 
$\varphi\in\Mor(N_\calf^K(Q))$.  
Set $\calf_0=N_\calf^K(Q)$ for short. 

It thus suffices to prove that $\calf_0=\calf_P(P)$.  
Assume otherwise:  then by \cite[Theorem 3.6]{BLO3}, 
there is $R\le P$ such that $O_p(\Aut_{\calf_0}(R))=1$, and in particular, 
such that $\Aut_{\calf_0}(R)$ is not a $p$-group.  By definition of 
$\calf_0$, we can assume that $R\ge Q$.  Hence there is 
$\Id\ne\alpha\in\Aut_{\calf_0}(R)$ which has finite order prime to $p$.  
Since $K$ is a discrete $p$-toral group and $\alpha|_Q\in{}K$, 
$\alpha|_Q=\Id_Q$.  Hence for $g\in{}R$, $g$ and $\alpha(g)$ have the same 
conjugation action on $Q$, and $g^{-1}\alpha(g)\in C_R(Q)\le{}Q$.  Thus 
$\alpha$ induces the identity on $R/Q$, so by \cite[Lemma 1.7(a)]{BLO3}, 
each $\alpha$-orbit in $R$ has $p$-power order.  This contradicts the 
assumption that $\alpha\ne\Id$ has order prime to $p$, so 
$\calf_0=\calf_P(P)$, and this finishes the proof.
\end{proof}

We will need the following application of Lemma \ref{l:centric}.

\begin{Lem} \label{Pfn:p'}
Let $\calf$ be a saturated fusion system over a discrete $p$-toral group 
$S$.  Then for each $P\le S$, the set of $S$-conjugacy classes of fully 
normalized subgroups $\calf$-conjugate to $P$ is finite of order prime to 
$p$.
\end{Lem}

\begin{proof} Since the conclusion depends only on the $\calf$-conjugacy 
class of $P$, we can assume that $P$ is fully normalized.  

Let $\repf(P,S)\fc\subseteq\repf(P,S)$ be the subset consisting of all 
classes of homomorphisms whose image is fully centralized.  Recall that 
$\repf(P,S)$ is finite by \cite[Lemma 2.5]{BLO3}.

\smallskip

\noindent\textbf{Step 1: } Assume first that $P$ is $\calf$-centric.  We 
first prove that $\repf(P,S)\fc=\repf(P,S)$ has order prime to $p$ in this 
case.  

Recall from \cite[\S\,3]{BLO3} that there is a functor $(P\mapsto P\8)$ 
from $\calf$ to itself, with the properties that $P\8\ge P$ for each $P\le 
S$, $(P\8)\8=P\8$, and the image $\calf\8$ of the functor contains finitely 
many $S$-conjugacy classes.  (We will be studying this functor in more 
detail in the next section.)

Let $P=P_0<P_1<P_2<\ldots<P_m=S$ be such that $P_i=N_S(P_{i-1})\8$ for each 
$i$.  The sequence is finite since there are only finitely many conjugacy 
classes in $\calf\8$ (and since $N_S(P_i)>P_i$ whenever $P_i<S$ by Lemma 
\ref{N_Q(P)>P}).  
Fix $0\le i<m$, and consider the restriction map 
	\[ \Res\: \repf(N_S(P_i),S) \Right5{} \repf(P_i,S). \]
This is injective by Lemma \ref{l:centric}, and $\Out_S(P_i)\cong 
N_S(P_i)/P_i\ne1$ is a finite $p$-group by \cite[Proposition 1.5(c)]{BLO3}.  
Let $Q/P_i$ be the stabilizer subgroup of $[\varphi]\in\repf(P_i,S)$ under 
the action of $N_S(P_i)/P_i$; then $\varphi$ extends to an 
$\calf$-morphism on $Q$ by axiom (II) in Definition \ref{sat.Frob.}.  Thus 
$\Im(\Res)$ is the fixed subset for the $N_S(P_i)/P_i$-action, and hence 
$|\repf(P_i,S)|\equiv|\repf(N_S(P_i),S)|$ (mod $p$).

By \cite[Proposition 3.3]{BLO3}, each $\psi\in\homf(N_S(P_i),S)$ extends to 
a unique morphism $\psi\8\in\homf(N_S(P_i)\8,S)$.  Thus restriction defines 
a bijection from $\homf(N_S(P_i)\8,S)$ to $\homf(N_S(P_i),S)$, and hence 
(after dividing out by the conjugation action of $S$) from 
$\repf(N_S(P_i)\8,S)$ to $\repf(N_S(P_i),S)$.  In other words, 
$|\repf(P_i,S)|\equiv|\repf(P_{i+1},S)|$ (mod $p$).  Thus
	\[ |\repf(P,S)| = |\repf(P_0,S)| \equiv |\repf(P_m,S)| = 
	|\outf(S)|, \]
where $|\outf(S)|$ is prime to $p$ by axiom (I) in Definition \ref{sat.Frob.}.

\smallskip

\noindent\textbf{Step 2: } Now assume $P<S$ is not $\calf$-centric; we 
claim that $|\repf(P,S)\fc|$ is prime to $p$.  Set 
	\[ \Gamma = \bigl\{ \gamma\in\autf(PC_S(P)) \,\big|\, 
	\gamma|_P=\Id_P \bigr\}
	\qquad\textup{and}\qquad
	T=\Aut_{C_S(P)}(PC_S(P))\nsg\Gamma \,. \] 
Then $\Gamma\cong\Aut_{C_\calf(P)}(C_S(P))$ and $T\cong\Inn(C_S(P))$ via 
restriction to $C_S(P)$, so $T\in\sylp{\Gamma}$.  

Now $\Gamma/T\le\outf(PC_S(P))$ acts on $\repf(PC_S(P),S)$ by right 
composition.  If $\varphi\in\homf(PC_S(P),S)$ and $\gamma\in\Gamma$ are 
such that $[\varphi]=[\varphi\circ\gamma]$ in $\repf(PC_S(P),S)$, then 
there is $g\in{}S$ such that $\varphi=c_g\circ\varphi\circ\gamma$, so 
$\varphi=\varphi\circ\gamma^{p^n}$ for some $n$, hence $|\gamma|$ is a 
power of $p$, and $\gamma\in{}T$.  Thus the $\Gamma/T$-action on 
$\repf(PC_S(P),S)$ is free.  The restriction map from $\repf(PC_S(P),S)$ to 
$\repf(P,S)\fc$ is surjective by axiom (II), and the inverse image 
of any $[\psi]\in\repf(P,S)\fc$ is one of the $\Gamma/T$-orbits.  Thus 
$|\repf(P,S)\fc|=|\repf(PC_S(P),S)|\big/|\Gamma/T|$.  Since 
$|\repf(PC_S(P),S)|$ is prime to $p$ by Step 1, $|\repf(P,S)\fc|$ is also 
prime to $p$.

\smallskip

\noindent\textbf{Step 3: } We are now ready to prove the lemma.  Let 
$\calp\fc$ be the set of $S$-conjugacy classes of subgroups fully 
centralized in $\calf$ and $\calf$-conjugate to $P$, and let 
$\calp\fn\subseteq\calp\fc$ be the subset of classes of fully normalized 
subgroups.  Let $\rho\:\repf(P,S)\fc\Right2{}\calp\fc$ be the map which 
sends the class of $\varphi$ to the class of $\varphi(P)$.

Let $N\big||\outf(P)|$ be the largest divisor prime to $p$.  For each $Q\in 
P^\calf$ such that $[Q]\in\calp\fc$, 
$|\rho^{-1}([Q])|=|\autf(Q)|\big/|\Aut_S(Q)|$.  Thus $p\big||\rho^{-1}(Q)|$ 
if $Q$ is not fully normalized, and $|\rho^{-1}(Q)|=N$ if $Q$ is fully 
normalized.  So $|\rho^{-1}(\calp\fn)|\equiv|\repf(P,S)\fc|$ (mod $p$), 
hence is prime to $p$ by Steps 1 and 2; and thus 
$|\calp\fn|=|\rho^{-1}(\calp\fn)|/N$ is also prime to $p$. 
\end{proof}

\bigskip


\section{The ``bullet" construction}

The bullet construction plays an important role in the theory of $p$-local 
compact groups, and is studied at length in \cite[Sec. 3]{BLO3}.  However, 
in that paper, we always assumed that the ambient fusion system is 
saturated. In this one, we will need some of these properties to prove 
saturation in certain situations, and hence must know them in a more 
general setting.   We start by recalling the definition.

\begin{Defi} \label{d:bullet} 
Fix a discrete $p$-toral group $S$, set $T=S_0$, and let $W\le\Aut(T)$ be a
subgroup which contains $\Aut_S(T)$. Let $m\ge0$ be such that
$\exp(S/T)=p^m$. 
\begin{itemize} 
\item For each $A\le{}T$, set $I(A)=C_T(C_W(A))\defeq \{t\in T\,|\, w(t) = 
t,\,\forall w\in C_W(A) \} \ge A$.  

\item For each $P\le{}S$, set $P^{[m]}=\gen{g^{p^m}\,|\,g\in{}P}$, and set
$P^\bullet=P{\cdot}I(P^{[m]})_0$.  

\item Let $P,Q\le{}S$ and $\varphi\in\Inj(P,Q)$ be such that 
$\varphi|_{P^{[m]}}=w|_{P^{[m]}}$ for some $w\in{}W$.  Then 
$\varphi^\bullet\:P^\bullet\Right2{}Q^\bullet$ denotes the unique map of 
sets, if it exists, such that $\varphi^\bullet(gh)=\varphi(g)w(h)$ for 
each $g\in{}P$ and $h\in{}I(P^{[m]})_0$. 

\end{itemize} 
\end{Defi}

Note that the definition of $P\8$ for $P\le{}S$ depends only on $S$ and on 
the choice of a subgroup $W\le\Aut(S_0)$.  For any fusion system 
$\calf$ over $S$, we associate to $\calf$ the construction ($P\mapsto P\8$) 
with $W= \Aut_\calf(S_0)$, and set $\calf^\bullet= \{\,P^\bullet\,|\, P\leq 
S\,\}$.

We will show later that under certain additional conditions, the map 
$\varphi\8$ is defined and a homomorphism for each morphism $\varphi$.

\begin{Lem} \label{dot-props}
Let $\calf$ be a fusion system over a discrete $p$-toral group $S$. 
Then 
\begin{enuma} 
\item $P\le{}Q\le{}S$ implies $P\8\le{}Q\8$, and 
\item $P\le{}S$ implies $(P\8)\8=P\8$ and $N_S(P)\le N_S(P\8)$.  
\end{enuma}
The following also hold if $\autf(S_0)$ is finite.
\begin{enuma} \setcounter{enumi}{2}
\item The set $\calf^\bullet$ contains finitely many $S$-conjugacy classes.
\item Assume $P_1\le P_2\le \cdots$ are subgroups of $S$, and 
$P=\bigcup_{i=1}^\infty{}P_i$.  Then there is $N>0$ such that 
$P_i{}^\bullet=P^\bullet$ for each $i>N$.
\end{enuma}
\end{Lem}

\begin{proof} By definition, the function $P\mapsto{}P^\bullet$ depends
only on $\autf(S_0)$ (and on $S$).  Hence (a), (b), and (c) hold by exactly
the same arguments as those used to prove them in \cite[Lemma 3.2]{BLO3}. 
The only exception is the last statement in (b) (not shown in 
\cite{BLO3}), which follows since for $g\in{}N_S(P)$, $g\in{}N_S(P^{[m]})$, 
and hence $g\in{}N_S(I(P^{[m]}))$.  
Note that the assumption $|\autf(S_0)|<\infty$ is needed to prove (c).

In the situation of (d), we have $P_1{}^\bullet\le P_2{}^\bullet\le\cdots$ 
by (a), and this sequence contains only finitely many subgroups by (c).  
Thus there are $N>0$ and $Q\le{}S$ such that $P_i{}^\bullet=Q$ for each 
$i>N$.  In particular, $P=\bigcup_{i=1}^\infty{}P_i\le{}Q$, so 
$P^\bullet\le{}Q^\bullet=Q$  by (a) and (b).  Since 
$P_i{}^\bullet\le{}P^\bullet$ for each $i$ by (a) again, this shows that 
$P^\bullet=Q=P_i{}^\bullet$ for $i>N$.
\end{proof}

While we do not assume here that fusion systems are saturated, we will, in 
most cases, assume the following condition on a fusion system $\calf$:
	\beq \textup{for each $P\le{}S_0$ and each 
	$\varphi\in\homf(P,S_0)$, $\varphi=w|_P$ for some 
	$w\in\autf(S_0)$.} \tag*{\sstar} \eeq
By \cite[Lemma 2.4(b)]{BLO3}, \sstar\ always holds if $\calf$ is 
saturated. 

\begin{Lem} \label{Fdot-invar}
Let $\calf$ be a fusion system over a discrete $p$-toral group $S$ such 
that $\calf$ satisfies {\sstar}.  Then for each $P\le{}S$ and each 
$\varphi\in\homf(P\8,S)$, $\varphi(P\8)=\varphi(P)\8\in\calf\8$.  In 
particular, $\calf\8$ is invariant under $\calf$-conjugacy.
\end{Lem}

\begin{proof} Let $m$ be such that $p^m=\exp(S/S_0)$, and set $W=\autf(S_0)$.  
Set $Q=\varphi(P)$.  We must show that $\varphi(P\8)=Q\8$.

Set $R=P\8\cap S_0\cap\varphi^{-1}(S_0)$.  Thus $R$ is the largest 
subgroup of $P\8$ such that $R$ and $\varphi(R)$ are both contained in 
$S_0$. By \sstar, there is $w\in\autf(S_0)$ such that 
$\varphi|_R=w|_R$.  In particular, $w(P^{[m]})=Q^{[m]}$, so 
$wC_W(P^{[m]})w^{-1}=C_W(Q^{[m]})$, and $w$ sends 
$I(P^{[m]})=C_T(C_W(P^{[m]}))$ isomorphically to $I(Q^{[m]})$.  Also, 
$\varphi(I(P^{[m]})_0)\le S_0$ since it is connected, so $I(P^{[m]})_0\le 
R$ and $\varphi(I(P^{[m]})_0)=w(I(P^{[m]})_0)=I(Q^{[m]})_0$.  So 
$\varphi(P\8)=Q\8$. 
\end{proof}

Throughout the rest of the section, we need to consider the following 
condition on a fusion system $\calf$ over a discrete $p$-toral group $S$:
	\beq \begin{split} 
	&\textup{for each $P\le S$ and each
	$\varphi\in\homf(P,S_0)$,}\\
	&\textup{there exists 
	$\widebar{\varphi}\in\homf(P\cdot C_S(P)_0,S_0)$ 
	such that $\widebar{\varphi}|_P=\varphi$.}
	\end{split} \tag*{\sstars} \eeq
As will be seen in the next lemma, this is a stronger version of \sstar.

\begin{Lem} \label{l:dots}
Fix a fusion system $\calf$ over a discrete $p$-toral group $S$ which 
satisfies \sstars. Then the following hold.
\begin{enuma} 
\item The fusion system $\calf$ satisfies condition \sstar.

\item For each $P,R\le{}S_0$ such that $P\le{}R\le{}P{\cdot}I(P)_0$, each 
$w\in{}W$, and each $\psi\in\homf(R,S)$ such that $\psi|_P=w|_P$, 
$\psi(R)\le{}S_0$ and $\psi=w|_R$.

\item For each $P,Q\le S$ and each $\varphi\in\homf(P,Q)$, $\varphi\8$ 
is defined and $\varphi\8\in\Hom(P\8,Q\8)$.

\item Assume $P\le{}Q\le{}P^\bullet$, $\varphi\in\homf(P,S)$, 
and $\widebar\varphi\in\homf(Q,S)$ are such that 
$\varphi=\widebar\varphi|_P$.  Then $\widebar\varphi=\varphi^\bullet|_Q$.
\end{enuma}
\end{Lem}

\begin{proof} Set $W=\autf(S_0)$.  Let $m$ be such that $p^m=\exp(S/S_0)$.

\smallskip

\noindent\textbf{(a) } For each $P\le S_0$ and 
$\varphi\in\homf(P,S_0)$, $\varphi$ extends to $w\in\autf(S_0)$ by 
\sstars.

\smallskip

\noindent\textbf{(b) }  Let $P\le{}S_0$, $P\le{}R\le{}P{\cdot}I(P)_0$, 
$w\in{}W$, and $\psi\in\homf(R,S)$ be as above.  
Set $P_*=\psi(P)=w(P)$ and $R_*=\psi(R)$. Then 
$P_*\le R_*$ are abelian, so each element of $\Aut_{R_*}(S_0)\le{}W$ acts 
via the identity on $P_*$, and hence via the identity on $I(P_*)$.  In 
other words, $[I(P_*),R_*]=1$.  By \sstars, $\psi^{-1}\in\isof(R_*,R)$ 
extends to a morphism $\chi\in\homf(R_*{\cdot}I(P_*)_0,S_0)$. 

Now, $\chi(P_*)=P$, and $\chi|_{P_*}=\psi^{-1}|_{P_*}=w^{-1}|_{P_*}$. 
Also, $\chi|_{P_*{\cdot}I(P_*)_0}=u|_{P_*{\cdot}I(P_*)_0}$ for some 
$u\in{}W$ by (a), $u|_{P_*}=w^{-1}|_{P_*}$, and hence 
$u|_{I(P_*)}=w^{-1}|_{I(P_*)}$ by definition of $I(-)$.  Thus 
$\chi|_{P_*{\cdot}I(P_*)_0}=w^{-1}|_{P_*{\cdot}I(P_*)_0}$, and so 
$\chi(P_*{\cdot}I(P_*)_0)=P{\cdot}I(P)_0\ge{}R$.  Since $\chi$ is 
injective and $\chi(R_*)=R$, this implies that 
$R_*\le{}P_*{\cdot}I(P_*)_0\le{}S_0$ and 
$\psi=(\chi|_{R_*})^{-1}=w|_R$.

\smallskip

\textbf{(c) }  Fix a morphism $\varphi\in\homf(P,Q)$, and let $w\in{}W$ be 
such that $\varphi|_{P^{[m]}}=w|_{P^{[m]}}$.  
Set $R=P^{[m]}\cdot (P\cap{}I(P^{[m]})_0)$. Then $P^{[m]} \leq R \leq 
P^{[m]}\cdot{}I(P^{[m]})_0$, and (b) implies $\varphi|_R=w|_R$.  
Hence $\varphi\8\:P\8\Right2{}Q\8$ is well defined as a map of sets by the 
formula $\varphi\8(gh)=\varphi(g)w(h)$ for all $g\in{}P$ and 
$h\in{}I(P^{[m]})_0$.  

To prove that $\varphi\8$ is a homomorphism, it remains to show that 
$w(ghg^{-1})=\varphi(g)w(h)\varphi(g)^{-1}$ for all $g\in{}P$ and 
$h\in{}I(P^{[m]})_0$.  Set $u=c_{\varphi(g)}^{-1}\circ w\circ c_g\in{}W$; 
we must show that $u|_{I(P^{[m]})_0}=w|_{I(P^{[m]})_0}$.  But 
$u|_{P^{[m]}}=w|_{P^{[m]}}$ since $\varphi|_{P^{[m]}}=w|_{P^{[m]}}$, and so 
$u|_{I(P^{[m]})}=w|_{I(P^{[m]})}$ by definition of $I(-)$.

\smallskip

\textbf{(d) }  Fix $P,Q\le{}S$ such that $P\le{}Q\le{}P^\bullet$.  Assume 
$\varphi\in\homf(P,S)$ and $\widebar\varphi\in\homf(Q,S)$ are such that 
$\varphi=\widebar\varphi|_P$.  Set $R=P^{[m]}$ and 
$Q'=Q\cap{}R{\cdot}I(R)_0$.  By definition, 
$\varphi\8|_{R{\cdot}I(R)_0}=w|_{R{\cdot}I(R)_0}$ for any $w\in{}W$ such 
that  $\varphi|_R=w|_R$.  By (b), 
$\widebar\varphi|_{Q'}=w|_{Q'}=\varphi^\bullet|_{Q'}$.   Since 
	\[ Q = Q\cap P\8 = Q\cap (P\cdot I(R)_0) 
	= P\cdot (Q\cap R\cdot I(R)_0) = PQ' \] 
and $\widebar\varphi|_P=\varphi$, this proves that 
$\widebar\varphi=\varphi^\bullet|_Q$. 
\end{proof}

Lemma \ref{l:dots}(c) allows us to extend each $\varphi\in\homf(P,Q)$ to a 
homomorphism  $\varphi\8\in\Hom(P\8,Q\8)$.  In the next lemma, we add 
some more hypotheses, enough to ensure that $\varphi\8\in\Mor(\calf)$.

\begin{Lem} \label{reduce2dots}
Fix a fusion system $\calf$ over a discrete $p$-toral group $S$.  Let 
$\calh$ be a family of subgroups of $S$ invariant under $\calf$-conjugacy.  
Assume the following:
\begin{enumerate}[\rm(i) ]
\item $\autf(S_0)$ is finite.

\item $\calf$ satisfies \sstars.

\item For all $P\in\calh$ and $P\le Q\le P^\bullet$, $Q\in\calh$. 

\item $\calf$ is $\calh$-generated and $\calh$-saturated.
\end{enumerate}
Then the following hold.
\begin{enuma} 
\item For each $\varphi\in\Mor(\calf)$, $\varphi^\bullet\in\Mor(\calf)$.
\label{4.6.a}

\item \label{4.6.b}
If, in addition, $\widehat\calh\supseteq\calh$ is also invariant under
$\calf$-conjugacy, and is such that $P\in\widehat\calh$ implies
$P^\bullet\in\calh$, then $\calf$ is $\widehat\calh$-saturated.
\end{enuma}
\end{Lem}

\begin{proof} \textbf{(a) }  Fix $\varphi\in\homf(P,Q)$; we must show 
$\varphi\8\in\homf(P\8,Q\8)$. It suffices to do this when 
$\varphi$ is an isomorphism.  Since $\calf$ is $\calh$-generated, $\varphi$ 
is a composite of isomorphisms which are restrictions of 
$\calf$-isomorphisms between subgroups in $\calh$.  So it suffices to 
consider the case where $\varphi$ extends to 
$\widebar\varphi\in\isof(\widebar{P},\widebar{Q})$ for some 
$\widebar{P},\widebar{Q}\in\calh$.  If 
$\widebar\varphi\8\in\Mor(\calf)$, then 
$\varphi\8\in\Mor(\calf)$ since it is a restriction of 
$\widebar\varphi\8$.  It thus suffices to prove (a) when $P$ and $Q$ 
are in $\calh$.  

By axiom (III), which holds by hypothesis for all subgroups in $\calh$, 
there is a subgroup $P'\le{}P\8$ containing $P$ which is maximal 
among subgroups  of $P\8$ satisfying 
$\varphi'\defeq\varphi\8|_{P'}\in\homf(P',Q\8)$.  Assume 
$P'<{}P\8$, and set $Q'=\varphi'(P')$. 

Since $P,Q\in\calh$, and since $P\le{}P'\le{}P\8$ and $Q\le{}Q'\le{}Q\8$,
$P',Q'\in\calh$ by (iii).  Since $\calf$ is $\calh$-saturated, there are
$R'\le{}S$ and $\psi'\in\isof(Q',R')$ such that $R'$ is receptive in
$\calf$.  Set $R=\psi'(Q)$ and $\psi=\psi'|_Q\in\isof(Q,R)$.  
By Lemma \ref{l:dots}(d), $\psi'=\psi\8|_{Q'}$, and
$\psi'\varphi'=(\psi\varphi)\8|_{P'}$.  We thus have isomorphisms of 
triples 
	\[ (P\8\ge P' \ge P) \RIGHT6{(\varphi\8,\varphi',\varphi)}{\cong} 
	(Q\8\ge Q'\ge Q) \RIGHT6{(\psi\8,\psi',\psi)}{\cong} 
	(R\8\ge R' \ge R) \,. \]

Set $P''=N_{P\8}(P')$ and $Q''=N_{Q\8}(Q')$. Since $P\8>P'$, $P''>P'$ and 
$Q''>Q'$ by Lemma \ref{N_Q(P)>P}. Also, $N_{\psi'}\ge{}Q''$ and 
$N_{\psi'\varphi'}\ge{}P''$ since there are (abstract) homomorphisms 
which extend $\psi'$ to $Q''$ and $\psi'\varphi'$ to $P''$. Hence $\psi'$ 
and $\psi'\varphi'$ extend to homomorphisms $\psi''\in\homf(Q'',S)$ and 
$\chi''\in\homf(P'',S)$, and $\psi''=\psi\8|_{Q''}$ and 
$\chi''=(\psi\varphi)\8|_{P''}$ by Lemma \ref{l:dots}(d) again. So 
$\psi''(Q'')=\chi''(P'')=N_{R\8}(R')$, and 
$\varphi\8|_{P''}=(\psi'')^{-1}\chi''$ is a morphism in $\calf$, which 
contradicts the maximality of $P'$.  We now conclude that $P'=P\8$, and 
hence that $\varphi\8\in\isof(P\8,Q\8)$.

\smallskip

\noindent\textbf{(b) } We are assuming $\calf$ is $\calh$-saturated, and 
want to show it is $\widehat\calh$-saturated.  Fix $P$ in 
$\widehat\calh{\sminus}\calh$.  Thus $P\8\in\calh$ by assumption.  Choose 
$R$ which is $\calf$-conjugate to $P\8$, and receptive and fully automized.  
Fix $\varphi\in\isof(P\8,R)$ and set $Q=\varphi(P)$.  By Lemma 
\ref{l:dots}(d), $\varphi=(\varphi|_P)\8$, and so $Q\8=R$.  

Set $H=\bigl\{\alpha\in\autf(Q\8)\,\big|\,\alpha(Q)=Q\bigr\}$. Since 
$Q\8=R$ is fully automized, $\Aut_S(Q\8)\in\sylp{\autf(Q\8)}$.  Thus $H$ 
acts by translation on $\autf(Q\8)/\Aut_S(Q\8)$, a finite set of order 
prime to $p$, and there is $\beta\in\autf(Q\8)$ such that the coset 
$\beta^{-1}\Aut_S(Q\8)$ is in an $H$-orbit of order prime to $p$.  In other 
words, the stabilizer subgroup $\beta^{-1}\Aut_S(Q\8)\beta\cap{}H$ has 
index prime to $p$ in $H$, and thus lies in $\sylp{H}$.  

By (a) (and the definition of $H$), restriction to $Q$ sends 
$\beta{}H\beta^{-1}$ surjectively onto $\autf(\beta(Q))$. Also,  
$\Aut_S(\beta(Q))$ is the image of $\Aut_S(Q\8)\cap\beta{}H\beta^{-1}$ 
under this surjection, and hence is a Sylow $p$-subgroup in 
$\autf(\beta(Q))$. Thus $\beta(Q)$ is fully automized.  (Note that 
$\beta(Q)\8=Q\8$.)

For any $R\in{}P^\calf=Q^\calf$ and any $\varphi\in\isof(R,Q)$, 
$N_\varphi\le N_{\varphi\8}$: if $\varphi 
c_g\varphi^{-1}=c_h\in\Aut_S(Q)$, then 
$\varphi\8c_g\varphi\8{}^{-1}=c_h\in\Aut_S(Q\8)$ by Lemma 
\ref{dot-props}(d).  Hence $Q$ and $\beta(Q)$ are both receptive in 
$\calf$ since $Q\8$ is.  Thus the conjugacy class $P^\calf$ contains a 
subgroup which is receptive and fully automized.  

Assume $P_\infty=\bigcup_{i=1}^\infty{}P_i$, where $\{P_i\}$ is an 
increasing sequence of subgroups in $\widehat\calh$.  Fix 
$\varphi\in\homf(P_\infty,S)$ such that $\varphi|_{P_i}\in\homf(P_i,S)$ for 
each $i$.  By Lemma \ref{dot-props}(d), for $i$ large enough, 
$P_i{}\8=P_\infty\8$.  By Lemma \ref{l:dots}(d), $\varphi$ is the 
restriction of $(\varphi|_{P_i})\8$, which is in $\calf$ by (a).  So 
$\varphi\in\homf(P_\infty,S)$.  
\end{proof}

\newcommand{\XXX}{\mathfrak{H}}

\bigskip


\section{Centric and radical subgroups determine saturation}

The main result in this section is Theorem \ref{H-sat-gen}, which gives 
sufficient conditions for a fusion system over a discrete $p$-toral group to 
be saturated.  

We will frequently refer to the conditions \sstar\ and \sstars\ of 
Section 3, which we recall here:
	\beqq \textup{$\forall$ $P\le{}S_0$, $\forall$
	$\varphi\in\homf(P,S_0)$, $\varphi=w|_P$ for some 
	$w\in\autf(S_0)$.} \tag*{\sstar} \eeqq
	\beqq 
	\textup{$\forall$ $P\le S$ and $\varphi\in\homf(P,S_0)$, $\exists$ 
	$\widebar{\varphi}\in\homf(P\cdot C_S(P)_0,S_0)$
	with $\widebar{\varphi}|_P=\varphi$.}
	\tag*{\sstars} \eeqq
By Lemma \ref{l:dots}(a), condition \sstars\ implies \sstar.

The following finiteness result will be needed.

\begin{Lem} \label{P^F/S-finite}
Let $\calf$ be a fusion system over a discrete $p$-toral group $S$ such 
that $\autf(S_0)$ is finite and \sstar\ holds. 
Then for each $P\le{}S$, there are only finitely many $S$-conjugacy classes
of subgroups of $S$ which are $\calf$-conjugate to $P$. 
\end{Lem}

\begin{proof} 
By \sstar\ and since 
$|\autf(S_0)|<\infty$, $P_0{}^\calf$ is finite.  By 
\cite[Lemma 1.4(a)]{BLO3}, for each $R\in{}P_0{}^\calf$, there are only 
finitely many $N_S(R)/R$-conjugacy classes of finite subgroups of 
$N_S(R)/R$ of any given order. (Compare \cite[Lemma 2.5]{BLO3}.)
\end{proof}

The main result in this section is the following theorem.  We refer to 
Definition \ref{d:H-sat} for the definitions of $\calh$-saturated and 
$\calh$-generated fusion systems.  

\begin{Thm} \label{H-sat-gen}
Fix a fusion system $\calf$ over a discrete $p$-toral group $S$
such that $\autf(S_0)$ is finite and condition \sstars\ holds.
Let $\calh$ be a family of subgroups of $S$ which satisfies 
the following:
\begin{enumi} 
\item $\calh$ is invariant under $\calf$-conjugacy. 

\item $\calh$ is closed in $\Sub(S)$, and $\calf$ is $\calh$-generated and 
$\calh$-saturated.

\item For all $P\in\calh$ and $P\le Q\le P^\bullet$, $Q\in\calh$.

\item If $P\in\calf^\bullet$ is $\calf$-centric and $P\notin\calh$, then there is 
$Q\in{}P^\calf$ such that 
	\[ O_p(\outf(Q)) \cap \Out_S(Q) \ne 1. \]
\end{enumi}
Then $\calf$ is saturated.
\end{Thm}

Theorem \ref{H-sat-gen} will be shown by following as closely as possible 
the proof of \cite[Theorem 2.2]{BCGLO1}.  The main difference is that since 
a discrete $p$-toral group can have infinitely many subgroups, the 
induction arguments used in \cite{BCGLO1} cannot be used here.  But the 
beginning steps are mostly unchanged: they are based on the concept of 
\emph{proper $\calp$-pairs} for an $\calf$-conjugacy class $\calp$ of 
subgroups of $S$.

\begin{Defi}
Let $\calf$ be a fusion system over a discrete $p$-toral group $S$, 
and let $\calp$ be an $\calf$-conjugacy class of subgroups of $S$.
\begin{itemize}
\item  A \emph{proper $\calp$-pair} is a pair of subgroups $(Q,P)$ 
such that $P<Q\le{}N_S(P)$ and $P\in\calp$.  

\item Two proper $\calp$-pairs $(Q,P)$ and $(Q',P')$ are
\emph{$\calf$-conjugate} if there is an isomorphism $\varphi\in\isof(Q,Q')$
such that $\varphi(P)=P'$. We let $(Q,P)^\calf$ denote the set of proper 
$\calp$-pairs which are $\calf$-conjugate to $(Q,P)$.

\item A proper $\calp$-pair $(Q,P)$ is \emph{fully normalized}
if $|N_{N_S(P)}(Q)|\ge|N_{N_S(P')}(Q')|$ for all 
$(Q',P')\in(Q,P)^\calf$.
\end{itemize}
\end{Defi}

Some basic properties of proper $\calp$-pairs are shown in the 
following lemma.

\begin{Lem}\label{l:pairs} 
Fix a fusion system $\calf$ over a discrete $p$-toral group $S$.  Assume
$\calh\subseteq\Sub(S)$ is closed, and invariant under $\calf$-conjugacy.  
Assume also that $\calf$ is $\calh$-generated and $\calh$-saturated.  
Let $\calp$ be an $\calf$-conjugacy class of subgroups of $S$, maximal
among those not in $\calh$.  
\begin{enuma}
\item \label{l:pairs1} If $(Q,P)$ is a fully normalized proper 
$\calp$-pair, then $Q$ is receptive in $\calf$ and
	\[ \Aut_{N_S(P)}(Q) \in \sylp{\Aut_{N_\calf(P)}(Q)}. \]

\item \label{l:pairs2} For each proper $\calp$-pair $(Q,P)$, and 
each $(Q',P')\in(Q,P)^\calf$ which is fully normalized,
there is a morphism 
	\[ \varphi\in\homf(N_{N_S(P)}(Q),N_S(P')) \] 
such that $\varphi(P)=P'$ and $\varphi(Q)=Q'$.

\item  \label{l:pairs3} Assume that $\Aut_\calf(S_0)$ is finite, condition 
\sstars\ holds, and $P\in\calh$, $P\le Q\le P\8$ imply $Q\in\calh$.  Then 
for each $P,\widehat{P}\in\calp$ such that $\widehat{P}$ is fully 
normalized in $\calf$, there is a morphism 
$\varphi\in\homf(N_S(P),N_S(\widehat{P}))$ such that 
$\varphi(P)=\widehat{P}$.  

\end{enuma}
\end{Lem}

\begin{proof}
\noindent\textbf{(a) }  For each proper $\calp$-pair $(Q,P)$, define 
        $$ K_P = \{\varphi\in \Aut(Q) \,|\, \varphi(P)=P\} \le \Aut(Q) \,.$$
Notice that 
	\begin{align} 
	\begin{split} 
	\Aut_S^{K_P}(Q) &= \Aut_S(Q)\cap K_P =\Aut_{N_S(P)}(Q) \\
	\autf^{K_P}(Q) &= \autf(Q)\cap K_P =\Aut_{N_{\calf}(P)}(Q)
	\label{e:pairs-a1} \end{split} \\
	N_S^{K_P}(Q) &= N_{N_S(P)}(Q)\,. \label{e:pairs-a2} 
	\end{align}


Assume the pair $(Q,P)$ is fully normalized.  Fix $Q'\in{}Q^\calf$ 
and 
$\alpha\in \Iso_{\calf}(Q,Q')$, and set $P'=\alpha(P)$. Then 
$(Q',P')\in(Q,P)^\calf$, and 
	\[ \9\alpha K_P = K_{P'} \defeq 
	\{\varphi\in \Aut(Q') \,|\, \varphi(P')=P'\}.\] 
Hence by \eqref{e:pairs-a2},
	$$|N_S^{K_P}(Q)| = |N_{N_S(P)}(Q)| \geq |N_{N_S(P')}(Q')| = 
	|N_S^{K_{P'}}(Q')| = |N_S^{{}^\alpha\!K_P}(Q')|, $$
and so $Q$ is fully $K_P$-normalized in $\calf$.
By Lemma~\ref{l:RS-K}(\ref{fullyKnorm}) (and since $Q\in\calh$ and $\calf$ 
is $\calh$-saturated), $Q$ is receptive and
fully $K_P$-automized in $\calf$.  So by \eqref{e:pairs-a1},
	$$ \Aut_{N_S(P)}(Q)  \in \sylp{\Aut_{N_\calf(P)}(Q)}. $$

\smallskip


\noindent\textbf{(b) }  Let $(Q,P)$, and $(Q',P')\in(Q,P)^\calf$, be proper 
$\calp$-pairs such that $(Q',P')$ is fully normalized.  Let $K_P\le\Aut(Q)$ 
and $K_{P'}\le\Aut(Q')$ be as in (a).  Since $(Q',P')$ is fully normalized, 
$Q'$ is receptive and fully $K_{P'}$-automized by (a).  

Choose some $\psi\in\isof(Q,Q')$ such that $\psi(P)=P'$.  By Lemma 
\ref{l:RS-K}(b), there are $\chi\in\autf^{K_{P'}}(Q')$ and 
$\varphi\in\homf(N_S^{K_P}(Q),N_S^{K_{P'}}(Q'))$ such that 
$\varphi|_Q=\chi\circ\psi$.  Then $\varphi(P)=\chi(P')=P'$, 
$\varphi(Q)=Q'$, and $N_S^{K_P}(Q)=N_{N_S(P)}(Q)$ and 
$N_S^{K_{P'}}(Q')=N_{N_S(P')}(Q')$ by \eqref{e:pairs-a2}.

\smallskip

\noindent\textbf{(c) }  Fix subgroups $P,\widehat{P}\in\calp$ such 
that $\widehat{P}$ is fully normalized in $\calf$.  We will construct a 
morphism in $\calf$ from $N_S(P)$ to $N_S(\widehat{P})$ which sends $P$ to 
$\widehat{P}$.  In particular, $\widehat{P}$ will be fully centralized in 
$\calf$, since its centralizer contains an injective image of the 
centralizer of any other subgroup in $\calp$.

Let $\widehat\calt$ be the set of all sequences 
	\beqq \xi=\bigr(P=P_0,Q_0,\varphi_0;P_1,Q_1,\varphi_1;\dots;
	P_{k-1},Q_{k-1},\varphi_{k-1};P_k=\widehat{P}\bigr) 
	\label{e:pairs5} \eeqq 
such that for each $0\le i\le k-1$, $(Q_i, P_i)$ is a proper $\calp$-pair, 
$\varphi_i\in\homf(Q_i,\N{S}{P_{i+1}})$, and $\varphi_i(P_i)=P_{i+1}$.  
Let $\calt\subseteq\widehat{\calt}$ be the subset of those sequences 
in which each $\varphi_i$ is maximal, in the sense that it cannot be 
extended in $\calf$ to a subgroup of $N_S(P_i)$ which properly contains 
$Q_i$. We give $\widehat{\calt}$ the partial ordering by inclusion (of 
sequences of the same length with the same $P_i$'s); then $\calt$ is the 
set of maximal elements in $\widehat{\calt}$.  Since $\calf$ is 
$\calh$-saturated, axiom (III) ensures that each element of 
$\widehat{\calt}$ is contained in an element of $\calt$.

We first check that $\calt\ne\emptyset$.  Choose any 
$\varphi\in\isof(P,\widehat{P})$.  Since $\calf$ is 
$\calh$-generated, there are subgroups $P_i<R_i\in\calh$ for $0\le i\le k$, and 
morphisms $\phi_i\in\homf(R_i,R_{i+1})$ for $i\le k-1$, such that 
$P_0=P$, $P_k=\widehat{P}$, $\phi_i(P_i)=P_{i+1}$, and $\varphi= 
\phi_{k-1}|_{P_{k-1}}\circ\cdots\circ\phi_1|_{P_1}\circ\phi_0|_{P_0}$. For 
each $i$, let $P_{i+1} = \phi_i(P_i)$. Then 
$N_{R_i}(P_i)>P_i$ by Lemma \ref{N_Q(P)>P} (and since $R_i>P_i$).  Set $Q_i = 
N_{R_i}(P_i)$ and $\varphi_i = \phi_i|Q_i$.  Then $(Q_i, P_i)$ is a proper 
$\calp$-pair, $\varphi_i\in\Hom_\calf(Q_i,N_S(P_{i+1}))$, and 
$(P_0,Q_0,\varphi_0;P_1,Q_1,\varphi_1;\ldots;P_k)\in\widehat{\calt}$.  
Thus $\widehat{\calt}\ne\emptyset$, and $\calt\ne\emptyset$ since each 
element is contained in a maximal element.

Lemma \ref{reduce2dots}(\ref{4.6.a}) applies in this situation by the extra 
hypotheses which were assumed.  Hence for any $\xi\in\calt$ as in 
\eqref{e:pairs5}, $\varphi_i$ extends to 
$\varphi_i^\bullet\in\homf(Q_i{}^\bullet,N_S(P_{i+1})^\bullet)$ for each 
$i$, and therefore to $Q_i{}^\bullet\cap \N{S}{P_i}$. The maximality of 
$\varphi_i$ implies that
	\begin{equation} Q_i=Q_i{}^\bullet\cap \N{S}{P_i}. 
	\label{e:pairs1} \end{equation}
The same argument applied to $\varphi_i^{-1}$ shows that 
\begin{equation} \varphi_i(Q_i)=\varphi_i(Q_i)^\bullet\cap \N{S}{P_{i+1}}.
	\label{e:pairs2} \end{equation}

Let $\calt_r\subseteq\calt$ be the subset of those $\xi$ for which there is 
no $1\le{}i\le{}k-1$ such that $Q_i=\N{S}{P_i}=\varphi_{i-1}(Q_{i-1})$. Let 
$\calt\Right1{R}\calt_r$ be the ``reduction'' map, which for every $1\le i 
\le k-1$ such that $Q_i=\N{S}{P_i}=\varphi_{i-1}(Q_{i-1})$, removes from 
the sequence the triple $(P_i, Q_i,\varphi_i)$, and replaces 
$\varphi_{i-1}$ by the composite $\varphi_i\circ\varphi_{i-1}$. Since 
$\calt\ne\emptyset$, the existence of the reduction map $R$ shows that 
$\calt_r\ne\emptyset$.  

For each $\xi\in\calt$, define 
	$$ I(\xi) = \{i=0,\ldots,k-1\,|\, Q_i<\N{S}{P_i} 
	\textup{ and } \varphi_i(Q_i)<\N{S}{P_{i+1}} \}. $$ 
If $\xi\in\calt$ and $I(\xi)\ne\emptyset$, define 
	$$ \lambda(\xi) =\min_{i\in I(\xi)}\bigl\{\,|Q_i{}^\bullet|\,\bigr\} 
	\ge (0,p)\,. $$
Since $\Aut_\calf(S_0)$ is assumed to be finite, $\calhbullet$ has finitely 
many $\calf$-conjugacy classes by Lemma \ref{dot-props}(c), so 
there are only a finite number of values that $\lambda(\xi)$ can take.

Assume there is a sequence $\xi\in\calt_r$ as in \eqref{e:pairs5} 
such that $I(\xi)=\emptyset$.  For each $0\le i<k$, if $Q_i<N_S(P_i)$, 
then $\varphi_i(Q_i)=N_S(P_{i+1})$ since $i\notin I(\xi)$, and hence 
$Q_{i+1}<N_S(P_{i+1})$ if $i+1\le k-1$ since $\xi$ is reduced.  Thus 
$|N_S(P_i)|>|Q_i|=|N_S(P_{i+1})|>\cdots>|N_S(P_k)|$ in this case, which 
contradicts the assumption that $P_k=\widehat{P}$ is fully normalized. 
Hence $Q_i=N_S(P_i)$ for each $i<k$, and 
$\varphi_{k-1}\circ\cdots\circ\varphi_0$ lies in 
$\homf(N_S(P),N_S(\widehat{P}))$ and sends $P$ to $\widehat{P}$.  

It remains to show that there exists $\xi\in\calt_r$ such that 
$I(\xi)=\emptyset$.  To see this, fix $\xi\in\calt_r$ as in 
\eqref{e:pairs5} 
such that $I(\xi)\ne\emptyset$.  We will construct 
$\widehat{\xi}\in\calt_r$ such that either $I(\widehat{\xi})=\emptyset$ or 
$\lambda(\widehat{\xi})>\lambda(\xi)$. Since $\lambda(\xi)$ can only take 
finitely many values, this will prove our claim.

Fix $i\in{}I(\xi)$, set $R_{i+1}=\varphi_i(Q_i)$, and choose a fully 
normalized proper $\calp$-pair $(U_i,V_i)$ which is $\calf$-conjugate to 
$(Q_i,P_i)$ (hence also to $(R_{i+1}, P_{i+1})$).  By (b), 
there are morphisms
	$$ \psi_i\in\homf\bigl(N_{N_S(P_i)}(Q_i{}),N_S(V_i)\bigr),
	\quad\textup{and}\quad
	\theta_i\in\homf\bigl(N_{N_S(P_{i+1})}(R_{i+1}),N_S(V_i)\bigr) $$
such that $\psi_i(P_i)=\theta_i(P_{i+1})=V_i$ and 
$\psi_i(Q_i)=\theta_i(R_{i+1})=U_i$.  
Upon replacing $(P_i,Q_i,\varphi_i)$ by 
	\[ \bigl(P_i,N_{N_S(P_i)}(Q_i),\psi_i; 
	V_i,\theta_i(N_{N_S(P_{i+1})}(R_{i+1})),
	\theta_i^{-1}\bigr)\,, \]
we get a new sequence $\xi_1\in\widehat{\calt}$. Upon applying axiom (III) 
again, we get a sequence $\xi_2\in\calt$ where $(P_i,Q_i,\varphi_i)$ 
has been replaced by 
	\[ \bigl(P_i,\widetilde{Q}_i,\widetilde\psi_i; 
	V_i,\widetilde{R}_{i+1},\widetilde\theta_i^{-1}\bigr) \]
for some maximal extensions $\widetilde{\psi}_i$ of $\psi_i$ and 
$\widetilde{\theta}_i$ of $\theta_i$. Also, $Q_i<N_S(P_i)$ and 
$R_{i+1}<N_S(P_{i+1})$ since $i\in{}I(\xi)$, so 
	\[ \widetilde{Q}_i\ge N_{N_S(P_i)}(Q_i) >Q_i 
	\qquad\textup{and}\qquad
	\widetilde{R}_{i+1}\ge N_{N_S(P_{i+1})}(R_{i+1})> R_{i+1}\,, \]
and $\widetilde{Q}_i{}\8>Q_i{}\8$ and $\widetilde{R}_{i+1}{}\8>R_{i+1}{}\8$ 
by \eqref{e:pairs1} and \eqref{e:pairs2}.  In particular, 
$|\widetilde{Q}_i{}\8|>\lambda(\xi)$ and 
$|\widetilde{R}_{i+1}{}\8|>\lambda(\xi)$.  

Upon repeating this procedure for  all $i\in I(\xi)$, we obtain a new 
element $\xi'\in\calt$ such that either $I(\xi')=\emptyset$ or 
$\lambda(\xi')>\lambda(\xi)$. Set  $\widehat{\xi}=R(\xi')\in\calt_r$.  Then 
either $I(\widehat{\xi})=\emptyset$ or 
$\lambda(\widehat{\xi})>\lambda(\xi)$.
Since the function $\lambda$ can only take a finite number of possible 
values, it follows by induction that there is $\xi \in\calt_r$ such that 
$I(\xi)=\emptyset$. 
\end{proof}

Very roughly, Lemma \ref{l:pairs} allows us to reduce the proof of Theorem 
\ref{H-sat-gen} to showing that the saturation properties hold for certain 
subgroups that are normal in the fusion system.  This case is handled by 
the following lemma.

\begin{Lem} \label{P*}
Let $\calf$ be a fusion system over a discrete $p$-toral group $S$ such 
that $\autf(S_0)$ is finite and condition \sstar\ holds. Fix a subgroup 
$Q\nsg{}S$, set $\calh=\{P\le S\,|\,P>Q\}$, and assume 
\begin{enumi}
\item $Q\nsg\calf$;
\item $\calf$ is $\calh$-generated and $\calh$-saturated; and 
\item either $Q$ is not $\calf$-centric, or 
$\Out_S(Q)\cap{}O_p(\outf(Q))\ne1$.
\end{enumi}
Then $Q$ is fully automized and receptive in $\calf$.
\end{Lem}

\begin{proof} When $\Gamma$ is a group containing a normal discrete 
$p$-torus $P$ of finite index, we let $O_p(\Gamma)$ be the inverse 
image in $\Gamma$ of the maximal normal $p$-subgroup $O_p(\Gamma/P)$ under 
the obvious projection. Equivalently, this is the largest normal 
discrete $p$-toral subgroup of $\Gamma$. Define
        $$ \widehat{Q} = \{x\in{}S\,|\, c_x\in{}O_p(\autf(Q)) \} \,. $$
Then $\widehat{Q} \nsg S$ by definition.  We claim
that $\widehat{Q}$ is strongly closed in $\calf$. Assume that $x\in \widehat{Q}$
is $\calf$-conjugate to $y\in S$. Since $Q$ is normal in $\calf$,
there exists $\psi \in \homf(\gen{x,Q},\gen{y,Q})$ which satisfies 
$\psi(Q)=Q$ and $\psi(x)=y$.  In particular, 
$\psi\circ c_x \circ \psi^{-1} = c_y$. It follows that $y\in \widehat{Q}$, since
$c_x\in{}O_p(\autf(Q))$.

Note also that $Q{\cdot}C_S(Q)\nsg \widehat{Q}$, and $\widehat{Q}/Q\cdot 
C_S(Q)\cong\Out_S(Q)\cap{}O_p(\outf(Q))$. If $Q$ is $\calf$-centric, 
then this last group is nontrivial by (iii), and if not, then 
$Q{\cdot}C_S(Q)>Q$ by definition. Thus $\widehat{Q}>Q$ in either case, and so 
$\widehat{Q}\in\calh$.

Consider the following statement:
	\beqq \textup{each $\varphi\in\autf(Q)$ extends to some 
	$\widebar{\varphi}\in\autf(\widehat{Q})$.} \label{ext-P*} \eeqq
We first prove that \eqref{ext-P*} implies the lemma, and then prove 
\eqref{ext-P*}.

\smallskip

\noindent\boldd{Point \eqref{ext-P*} implies the lemma: } Since $Q$ is 
normal in $\calf$ and $\widehat{Q}$ is strongly closed, each of them is the 
only subgroup in its $\calf$-conjugacy class.  So $Q$ and $\widehat{Q}$ are both 
fully centralized and fully normalized in $\calf$. Also, $\widehat{Q}$ is 
receptive and fully automized in $\calf$, since $\calf$ is 
$\calh$-saturated and $\widehat{Q}{}^\calf=\{\widehat{Q}\}\subseteq\calh$.  

By \eqref{ext-P*}, the restriction map from $\autf(\widehat{Q})$ to $\autf(Q)$ is 
surjective, and so $\Aut_S(Q)\in\sylp{\autf(Q)}$ since 
$\Aut_S(\widehat{Q})\in\sylp{\autf(\widehat{Q})}$.  Thus $Q$ is fully automized in $\calf$.

Next we prove that $Q$ is receptive in $\calf$.  Fix $\varphi\in\autf(Q)$. 
As usual, let $N_\varphi$ be the group of all $g\in{}N_S(Q)$ such that 
$\varphi c_g\varphi^{-1}\in\Aut_S(Q)$. By \eqref{ext-P*}, $\varphi$ extends 
to some $\psi\in\autf(\widehat{Q})$. Consider the groups of automorphisms
	\begin{align*} 
	K &= \bigl\{\chi\in\Aut_S(\widehat{Q})\,\big|\, \chi|_Q=c_x 
	\textup{ some }x\in{}N_\varphi \bigr\} \\ 
	K_0 &= \bigl\{\chi\in\autf(\widehat{Q})\,\big|\, 
	\chi|_Q=\Id_Q \bigr\} \nsg \autf(\widehat{Q}). 
	\end{align*}
By definition, for all $x\in{}N_\varphi$, we have
$(\psi{}c_x\psi^{-1})|_{Q}=\chi|_{Q}$ for some
$\chi\in\Aut_S(\widehat{Q})$.  In other words, as subgroups of
$\Aut(\widehat{Q})$,
	$$ \9\psi K
	\le \Aut_S(\widehat{Q}) \cdot K_0. $$

Now, $\Aut_S(\widehat{Q})\in\sylp{\Aut_S(\widehat{Q})K_0}$ since 
$\widehat{Q}$ is fully automized, so there are 
$\omega\in\Aut_S(\widehat{Q})$ and $\chi\in{}K_0$ such that 
$\9{\omega\chi}(\9\psi K)\le \Aut_S(\widehat{Q})$.  Hence 
$\9{\chi\psi}K\le\Aut_S(\widehat{Q})$.  Since $\widehat{Q}$ is receptive in 
$\calf$, $\chi\psi\in\autf(\widehat{Q})$ extends to a morphism 
$\widebar{\varphi}$ defined on $N_{\chi\psi}\geq 
N_S^K(\widehat{Q})\ge{}N_\varphi$, and 
$\widebar{\varphi}|_Q=\psi|_Q=\varphi$ since $\chi|_Q=\Id_Q$.

\smallskip

\noindent\boldd{Proof of \eqref{ext-P*}: } Since $\calf$ is 
$\calh$-generated, each $\varphi\in\autf(Q)$ is a composite of 
automorphisms of $Q$ which extend to strictly larger subgroups.  So it 
suffices to show \eqref{ext-P*} when $\varphi$ itself extends to 
some $P>Q$.  

Let $\calx$ be the set of all subgroups $P\in\calh$ such that 
$\varphi=\widebar\varphi|_Q$ for some $\widebar\varphi\in\homf(P,S)$.  We 
are assuming that $\calx\ne\emptyset$. We claim that
	\beqq \textup{$P\in\calx$\quad $\implies$\quad 
	$N_{\4QP}(P)\in\calx$.} \label{e:P*2} \eeqq
Assume this, and fix $P_1\in\calx$.  If 
$|\4QP_1/P_1|<\infty$, then by repeated application of \eqref{e:P*2} and 
Lemma \ref{N_Q(P)>P}, we get that $\4Q\in\calx$.  If not, then 
$|\4Q/Q|=\infty$, and $\4Q_0Q/Q$ is a nontrivial discrete $p$-torus.  Set 
$P_2=N_{\4Q_0P_1}(P_1)$ and $P_3=P_2\cap\4Q_0Q$.  Then $P_2\in\calx$ by 
\eqref{e:P*2}, $P_1\le P_2\le\4Q_0P_1$, and $P_2>P_1$ if $P_2\ngeq\4Q_0$ by 
Lemma \ref{N_Q(P)>P}.  So $P_3>Q$, and hence $P_3\in\calx$.  Now set 
$P_4=N_{\4Q}(P_3)$, so $P_4\in\calx$ by \eqref{e:P*2} again, 
$P_4\ge\4Q_0Q$ since $\4Q_0Q/Q$ is abelian, hence 
$[\4Q:P_4]<\infty$, and $\4Q\in\calx$ by earlier remarks.  Since
$\4Q\in\calx$, $\varphi$ extends to some $\4\varphi\in\homf(\4Q,S)$, and 
$\4\varphi(\4Q)=\4Q$ since $\4Q$ is strongly closed in $\calf$. This shows 
that \eqref{e:P*2} implies \eqref{ext-P*}.

It remains to prove \eqref{e:P*2}.  Fix $P\in\calx\subseteq\calh$, and 
let $\til\varphi\in\isof(P,P_2)$ be an extension of $\varphi$.  
Choose some $P_3\in P^\calf=P_2{}^\calf$ which is fully automized 
and receptive in $\calf$.  Notice that $P_3>Q$, since $Q\nsg\calf$ by 
(i). Let $\psi\in\isof(P_2,P_3)$ be any isomorphism.  Upon replacing 
$\widetilde{\varphi}$ by $\psi\circ\widetilde{\varphi}$, we can arrange 
that $P_2$ be fully automized and receptive.  

Consider the groups of automorphisms
	\begin{align*}  
	L&=\bigl\{\chi\in\autf(P_2)\,\big|\, \chi|_{Q}\in O_p(\autf(Q)) 
	\bigr\} \\
	L_0&=\bigl\{\chi\in\autf(P_2)\,\big|\, \chi|_{Q}=\Id_Q \bigr\} \,.
	\end{align*}
Both $L$ and $L_0$ are normal subgroups of $\autf(P_2)$. Also, $L/L_0$ is 
a discrete $p$-toral group, since there is a monomorphism $L/L_0\rightarrow 
O_p(\autf(P_2))$.  Since $P_2$ is fully automized, 
$\Aut_S(P_2)\in \sylp{\autf(P_2)}$, and hence 
$\Aut_S^L(P_2)\in\sylp{L}$ and $L=\Aut_S^L(P_2)L_0$.

Thus $P_2$ is fully $L$-automized in $\calf$. By Lemma 
\ref{l:RS-K}(b), there are $g\in N_S^L(P_2)=N_{\widehat{Q}}(P_2)$, $\chi\in L_0$, 
and $\psi\in\homf(N_S^{L^\varphi}(P)P,N_S^L(P_2)P_2)$, such that 
$\psi|_{P}=(c_g\circ\chi)\circ\varphi$. Also, 
$N_S^{L^\varphi}(P)P=N_{\widehat{Q}P}(P)$.  Upon replacing $\psi$ by 
$c_g^{-1}\circ\psi$, we can assume that $g=1$.  Then $\psi|_Q=\varphi|_Q$, 
and thus $N_{\4QP}(P)\in\calx$.  
\end{proof}

The following lemma combines Lemmas \ref{l:pairs} and \ref{P*}.

\begin{Lem} \label{P**}
Let $\calf$ be a fusion system over a discrete $p$-toral group $S$  such 
that $\autf(S_0)$ is finite and condition \sstars\ holds.  Fix a family 
$\calh\subsetneqq\Sub(S)$, and a subgroup $P\le S$ which is maximal in 
$\Sub(S){\sminus}\calh$.  Assume 
\begin{enumi} 
\item $\calh$ is invariant under $\calf$-conjugacy; 

\item $\calh$ is closed in $\Sub(S)$, and $\calf$ is $\calh$-generated and 
$\calh$-saturated;

\item $P\in\calh$ and $P\le Q\le P\8$ imply $Q\in\calh$; and 
\item either $P$ is not $\calf$-centric, or 
$\Out_S(P)\cap{}O_p(\outf(P))\ne1$.
\end{enumi}
Then $\calf$ is $(\calh\cup P^\calf)$-saturated.
\end{Lem}

\begin{proof} 
By assumption, all overgroups of subgroups in $P^\calf$ are in 
$\calh$.  Since $P{}^\calf$ contains only finitely many 
$S$-conjugacy classes by Lemma \ref{P^F/S-finite}, there is a subgroup 
$\widehat{Q}\in P{}^\calf$ which is fully normalized in $\calf$.  By 
Lemma \ref{l:pairs}(\ref{l:pairs3}), for each $Q\in P^\calf$, there is 
$\chi_Q\in\homf(N_S(Q),N_S(\widehat{Q}))$ such that $\chi_Q(Q) = 
\widehat{Q}$. Let $\caln\subseteq\calh$ be the family of all subgroups 
of $N_S(\widehat{Q})$ which strictly contain $\widehat{Q}$.  We claim that 
the normalizer system $N_\calf(\widehat{Q})$ is $\caln$-saturated.  

For each $Q\in\caln$, the pair $(Q,\widehat{Q})$ is a proper 
$P^\calf$-pair. Let $(Q',P')$ be a proper $P^\calf$-pair which is 
$\calf$-conjugate to $(Q,\widehat{Q})$ and fully normalized in $\calf$. 
Then $\widehat{Q} = \chi_{P'}(P')< Q''\defeq\chi_{P'}(Q') \le 
N_S(\widehat{Q})$, so $Q''\in\caln$, and $(Q'',\widehat{Q})$ is a proper 
$P^\calf$-pair. Furthermore, since $(Q',P')$ is fully normalized, and since 
$\chi_{P'}$ is a monomorphism, $(Q'',\widehat{Q})$ is also fully normalized 
in $\calf$. If $\alpha\in\isof(Q,Q')$ is such that 
$\alpha(\widehat{Q})=P'$, then $(\chi_{P'}|_{Q'})\circ \alpha$ is a 
morphism in $N_\calf(\widehat{Q})$ which sends $(Q,\widehat{Q})$ onto 
$(Q'',\widehat{Q})$. Hence by Lemma \ref{l:pairs}(\ref{l:pairs1}), applied 
with with $\calp = P^\calf$, the subgroup $Q''$ is 
receptive in $\calf$ (hence in $N_\calf(\widehat{Q})$), and is fully 
automized in $N_\calf(\widehat{Q})$. This shows that every $Q\in\caln$ is 
$N_\calf(\widehat{Q})$-conjugate to some $Q''\in\caln$ which is fully 
automized and receptive in the normalizer fusion system. Axiom (III) holds 
for $N_\calf(\widehat{Q})$ with respect to the family $\caln$, since 
$\caln$ is closed under overgroups, $\caln\subseteq\calh$, and $\calf$ is 
$\calh$-saturated. Thus $N_\calf(\widehat{Q})$ is $\caln$-saturated.

Conjugation by $(\chi_P)|_P\in\isof(P,\widehat{Q})$ sends $\autf(P)$ 
isomorphically to $\autf(\widehat{Q})$ and $\Aut_S(P)$ into 
$\Aut_S(\widehat{Q})$.  So if $P$ and $\widehat{Q}$ are not 
$\calf$-centric, then (iv) implies that $\Out_S(\widehat{Q})\cap 
O_p(\outf(\widehat{Q}))\ne1$.  Thus by Lemma \ref{P*}, $\widehat{Q}$ is 
receptive and fully automized in $N_\calf(\widehat{Q})$.  Hence 
$\widehat{Q}$ is fully automized in $\calf$.  It is receptive in $\calf$ 
since for each $Q\in P^\calf$, there is an $\calf$-morphism which sends $Q$ 
onto $\widehat{Q}$ and $N_S(Q)$ into $N_S(\widehat{Q})$.  

Axiom (III) holds for $\calh\cup{}P^\calf$ since it holds for $\calh$ and 
since no subgroup in $P^\calf$ contains any subgroups in $\calh$.  
Thus $\calf$ is $(\calh\cup{}P^\calf)$-saturated. 
\end{proof}

We are now ready to prove Theorem \ref{H-sat-gen}. The inductive Lemma 
\ref{P**} would suffice to prove the theorem for fusion systems over finite 
$p$-groups.  But for fusion systems over discrete $p$-toral groups, because 
our groups have infinite chains of subgroups, the results of Section 3 are 
needed to allow an induction proof.

\begin{proof}[\textbf{Proof of Theorem \ref{H-sat-gen}}]  Let 
$\mathfrak{K}$ be the set of all closed families of subgroups 
$\calk\subseteq\Sub(S)$ such that 
\begin{itemize} 
\item $\calk\supseteq\calh$, 
\item $\calk$ is invariant under $\calf$-conjugacy, 
\item $P\in\calk$ and $P\le Q\le P\8$ implies $Q\in\calk$, and 
\item $\calf$ is $\calk$-saturated.  
\end{itemize}
We must show that $\Sub(S)\in\mathfrak{K}$.

Assume otherwise, and choose $\calk_0\in\mathfrak{K}$ for which 
$\calk_0\cap\calf\8$ contains the largest possible number of 
$\calf$-conjugacy classes.  Set $\calk_1=\{P\le{}S\,|\,P\8\in\calk_0\}$.  Then 
$\calk_1\in\mathfrak{K}$ by Lemma \ref{reduce2dots}(b).  If 
$\calk_1\subsetneqq\Sub(S)$, then $\calf\8\not\subseteq\calk_0$; let $P$ be 
maximal among subgroups in $\in\calf\8{\sminus}\calk_0$.  Then $P$ is 
maximal among subgroups in $\Sub(S){\sminus}\calk_1$, $\calk_1\cup 
P^\calf\in\mathfrak{K}$ by Lemma \ref{P**}, and this contradicts the 
maximality assumption on $\calk_0$.
\end{proof}

\bigskip


\section{Extensions of $p$-local compact groups}
\label{s:ext}

The main result in this section is a version of \cite[Theorem 9]{link-ext} 
which describes how to extend a $p$-local compact group by a finite group.  
But before proving this, we need to show some of the basic properties of 
linking systems over discrete $p$-toral groups.  All of the results in this 
section are generalizations of results in \cite{link-ext} or earlier papers 
about linking systems over finite $p$-groups.

We first look at automorphisms of fusion and linking systems.  The 
definitions are the same as in the finite case (e.g., \cite[Definition 
1.13]{AOV1}).

\begin{Defi} \label{d:AutI}
\begin{enuma} 
\item For any fusion system over a discrete $p$-toral group $S$, an 
automorphism $\alpha$ of $S$ is \emph{fusion preserving} if there is an 
automorphism $\widehat{\alpha}$ of $\calf$ which sends an object $P$ to 
$\alpha(P)$ and sends a morphism $\varphi$ to $\alpha\varphi\alpha^{-1}$ 
(after restricting $\alpha$ in the obvious way).  Let $\Aut(S,\calf)$ be 
the group of fusion preserving automorphisms of $S$, and set 
$\Out(S,\calf)=\Aut(S,\calf)/\autf(S)$.

\item For any linking system $\call$ over a discrete $p$-toral group $S$, 
an automorphism of categories $\alpha\:\call\Right2{\cong}\call$ is 
\emph{isotypical} if for each $P\in\Ob(\call)$, 
$\alpha(\delta_P(P))=\delta_{\alpha(P)}(\alpha(P))$.  Let 
$\Aut\typ^I(\call)$ be the group of isotypical automorphisms of $\call$ 
which send inclusions to inclusions.  

\item Each $\gamma\in\Aut_{\call}(S)$ acts on the set $\Mor(\call)$ by 
composing on the left or right with $\gamma$ and its restrictions.  More 
precisely, for any $\varphi\in\Mor_{\call}(P,Q)$, set 
	\beqq \gamma\varphi=\gamma|_{Q,\pi(\gamma)(Q)}\circ\varphi 
	\in\Mor_{\call}(P,\pi(\gamma)(Q)) \label{e:lact} \eeqq
and
	\beqq \varphi\gamma=\varphi\circ\gamma|_{\pi(\gamma)^{-1}(P),P} 
	\in \Mor_{\call}(\pi(\gamma)^{-1}(P),Q). \label{e:ract} \eeqq

\item For each $\gamma\in \Aut_\call(S)$, let 
$c_\gamma\in\Aut\typ^I(\call)$ be the automorphism which sends an object 
$P$ to $\pi(\gamma)(P)$, and a morphism $\varphi$ to 
$\gamma\varphi\gamma^{-1}$. Set 
	\[ \Out\typ(\call) = 
	\Aut\typ^I(\call)\big/\{c_\gamma\,|\,\gamma\in\Aut_\call(S)\}\,. \]

\end{enuma}
\end{Defi}

The argument that $\Aut\typ^I(\call)$ is a group in this situation is 
exactly the same as that used in \cite[Lemma 1.14]{AOV1} when $S$ is 
finite.  

\begin{Prop} \label{AutI}
Let $\SFL$ be a linking triple, with structure functors 
	\[ \calt_{\Ob(\call)}(S)\Right2{\delta} \call \Right2{\pi}\calf\,. 
	\]  
Fix $\alpha\in\Aut\typ^I(\call)$.  Let $\beta\in\Aut(S)$ be the restriction of 
$\alpha$ to $S$ under the identification $S\cong\delta_S(S)$; thus 
$\alpha(\delta_S(g))=\delta_S(\beta(g))$ for all $g\in{}S$.  Then 
$\beta\in\Aut(S,\calf)$.  Furthermore, 
$\pi\circ\alpha=\widehat{\beta}\circ\pi$, where 
$\widehat{\beta}\in\Aut(\calf)$ is the automorphism which sends $P$ to 
$\beta(P)$ and $\varphi$ to $\beta\varphi\beta^{-1}$.  
\end{Prop}

\begin{proof}  Clearly, $\alpha(S)=S$, and hence $\alpha$ sends 
$\delta_S(S)$ to itself.  Thus $\beta$ is well defined.  

For each $P\in\Ob(\call)$ and $g\in{}P$, since $\alpha$ sends inclusions to 
inclusions, it sends 
	\[ 
	\vcenter{\xymatrix@C=50pt{
	\I02P \ar[r]^{\iota_{P,S}} \ar[d]_{\delta_P(g)} & \I02S 
	\ar[d]^{\delta_S(g)} \\
	\I02P \ar[r]^{\iota_{P,S}} & \I02S }} 
	\qquad\textup{to}\qquad
	\vcenter{\xymatrix@C=50pt{
	\alpha(P) \ar[r]^{\iota_{\alpha(P),S}} 
	\ar[d]_{\alpha(\delta_P(g))} 
	& S \ar[d]^{\delta_S(\beta(g))} \\
	\alpha(P) \ar[r]^{\iota_{\alpha(P),S}} & S \rlap{\,.} }}  \]
The first square commutes by axiom (C) in Definition \ref{L-cat}, so the 
second also commutes.  Since restrictions in $\call$ are uniquely defined 
by Proposition \ref{L-prop}(\ref{x:uniqrest}), this shows that 
$\alpha(\delta_P(g))=\delta_S(\beta(g))|_{\alpha(P)}= 
\delta_{\alpha(P)}(\beta(g))$.  Hence
	\[ \delta_{\alpha(P)}(\beta(P)) = \alpha(\delta_P(P))
	= \delta_{\alpha(P)}(\alpha(P))\,, \]
where the second equality holds since $\alpha$ is isotypical.  Thus 
$\alpha(P)=\beta(P)$ since $\delta_{\alpha(P)}$ is injective.

Fix $P,Q\in\Ob(\call)$ and $\psi\in\Mor_\call(P,Q)$, and set 
$\varphi=\pi(\psi)\in\homf(P,Q)$.  For each $g\in{}P$, $\alpha$ sends 
	\[ 
	\vcenter{\xymatrix@C=30pt{
	\I02P \ar[r]^{\psi} \ar[d]_{\delta_P(g)} & \I02Q 
	\ar[d]^{\delta_Q(\varphi(g))} \\
	\I02P \ar[r]^{\psi} & \I02Q }} 
	\qquad\textup{to}\qquad
	\vcenter{\xymatrix@C=30pt{
	\beta(P) \ar[r]^{\alpha(\psi)} \ar[d]_{\delta_{\beta(P)}(\beta(g))} 
	& \beta(Q) \ar[d]^{\delta_{\beta(Q)}(\beta(\varphi(g)))} \\
	\beta(P) \ar[r]^{\alpha(\psi)} & \beta(Q) \rlap{\,.} }}  \]
The first square commutes by axiom (C), so the second also commutes.  By 
(C) again,  $\delta_{\beta(Q)}(\beta(\varphi(g)))$ can be replaced by 
$\delta_{\beta(Q)}(\pi(\alpha(\psi))(\beta(g)))$, leaving the second square 
commutative. Since $\delta_{\beta(Q)}$ is a monomorphism, and since 
morphisms in $\call$ are epimorphisms by Proposition 
\ref{L-prop}(\ref{x:monoepi}), it 
follows that $\beta(\varphi(g))=\pi(\alpha(\psi))(\beta(g))$.  Thus 
	\beqq \pi(\alpha(\psi)) = \beta\varphi\beta^{-1} 
	= \beta\pi(\psi)\beta^{-1}. \label{e:7.3} \eeqq

In particular, $\beta\varphi\beta^{-1}\in\homf(\beta(P),\beta(Q))$ for each 
$P,Q\in\Ob(\call)$ and $\varphi\in\homf(P,Q)$.  Since $\Ob(\call)$ includes 
all subgroups which are $\calf$-centric and $\calf$-radical, all morphisms 
in $\calf$ are composites of restrictions of morphisms between objects of 
$\call$ by Alperin's fusion theorem in the version of \cite[Theorem 
3.6]{BLO3}.  So $\beta\varphi\beta^{-1}\in\Mor(\calf)$ for all 
$\varphi\in\Mor(\calf)$, and $\beta\in\Aut(S,\calf)$.  

Thus there is a well defined functor $\widehat{\beta}$ from $\calf$ to 
itself which sends each $P\le{}S$ to $\beta(P)$ and sends each 
$\varphi\in\homf(P,Q)$ to $\beta\varphi\beta^{-1}$.  This is an 
automorphism of the category $\calf$ by the same argument applied to 
$\alpha^{-1}$, and $\pi\circ\alpha=\widehat{\beta}\circ\pi$ by 
\eqref{e:7.3}.  
\end{proof}

We are now ready to define the structures which will be needed to construct 
extensions of linking systems.

\begin{Defi} \label{d:U(L,G)}
Fix a linking triple $\SFL$ and a finite group $G$.  
\begin{enuma} 

\item An \emph{extension pair for $\call$ and $G$} is a pair 
$(\widehat{\Gamma},\tau)$, where $\widehat{\Gamma}$ is an extension of 
$\Gamma\defeq\Aut_{\call}(S)\nsg\widehat{\Gamma}$ by $G$, and where 
$\tau\:\widehat{\Gamma}\Right4{}\Aut\typ^I(\call)$ is a homomorphism which 
makes both triangles in the following diagram commute:
    \[ \vcenter{\xymatrix@C=60pt@R=40pt{
    \I12\Gamma \ar[r]^-{\conj}
    \ar[d]_{\incl}^(.4){\quad \textup{\normalsize\pplus}} &
    \I12\Aut\typ^I(\call) \ar[d]^{(\alpha\mapsto\alpha_{S})}
    _(.6){\textup{\normalsize\ppluss \quad }} \\
    \I12\widehat{\Gamma} \ar[r]^(.45){\conj} \ar[ur]^{\tau} & \I12\Aut(\Gamma)
    \rlap{~,}
    }} \]
\end{enuma}
Fix an extension pair $U=(\widehat{\Gamma},\tau)$ for $\call$ and $G$.  
Let $\rho\:\widehat{\Gamma}\Right2{}G$ be the surjection with kernel $\Gamma$. 
\begin{enuma} \stepcounter{enumi}
\item Let $\call_U=\call_{(\widehat{\Gamma},\tau)}$ be the category with 
$\Ob(\call_U)=\Ob(\call)$, and with 
	$$ \Mor(\call_U) = \Mor(\call) \times_{\Gamma} \widehat{\Gamma}  
	= \bigl(\Mor(\call) \times \widehat{\Gamma}\bigr)\big/{\sim}, $$
where $(\varphi,\gamma)\sim(\varphi',\gamma')$ if and only if there is 
$\lambda\in\Gamma$ such that $\varphi'=\varphi\lambda$ and 
$\gamma'=\lambda^{-1}\gamma$.  Here, $\varphi\lambda$ 
is as defined in \eqref{e:ract}.  When 
$\varphi\in\Mor_\call(\gamma(P),Q)$, the equivalence class of 
$(\varphi,\gamma)$ is denoted $\5{\varphi,\gamma}\in\Mor_{\call_U}(P,Q)$. 
Composition in $\call_U$ is defined by
	\[ \5{\psi,\eta}\circ\5{\varphi,\gamma}
	=\5{\psi\circ\tau(\eta)(\varphi),\eta\gamma}. \]
Here, $\tau(\eta)(\varphi)\in\Mor_{\call}(\eta\gamma(P),\eta(Q))$ when 
$\varphi\in\Mor_{\call}(\gamma(P),Q)$ (where we write 
$\eta(P)=\tau(\eta)(P)$, etc.).  

\end{enuma}
\end{Defi}

Thus when $U=(\widehat{\Gamma},\tau)$ is an extension pair for $\call$ and $G$, 
$\5{\varphi\lambda,\gamma}=\5{\varphi,\lambda\gamma}$ in $\Mor(\call_U)$ 
for all $\varphi\in\Mor(\call)$, $\lambda\in\Gamma$, and 
$\gamma\in\widehat{\Gamma}$. 

To show composition in $\call_U$ is well defined, we note 
that for all $\psi,\varphi\in\Mor(\call)$, $\mu,\lambda\in\Gamma$, 
and $\eta,\gamma\in\widehat{\Gamma}$ with appropriate domain and range,
	\begin{multline*}  
	\5{\psi\mu,\eta}\circ\5{\varphi\lambda,\gamma} 
	=\5{\psi\mu\circ\tau(\eta)(\varphi\lambda),\eta\gamma} 
	=\5{\psi\mu\circ\tau(\eta)(\varphi),
	(\eta\lambda\eta^{-1})\eta\gamma} \\
	=\5{\psi\mu\circ\tau(\eta)(\varphi)\mu^{-1},
	\mu\eta\lambda\gamma}
	=\5{\psi\circ\tau(\mu\eta)(\varphi),
	\mu\eta\lambda\gamma}
	=\5{\psi,\mu\eta}\circ\5{\varphi,\lambda\gamma} 
	\end{multline*}
The second equality follows from the commutativity of triangle \ppluss\ in 
Definition \ref{d:U(L,G)}, and the fourth from that of \pplus.

We are now ready to state and prove the main result of this section.  The 
following theorem  is a generalization to $p$-local compact groups of 
\cite[Theorem 9]{link-ext} (which in turn was a generalization of 
\cite[Theorem 4.6]{BCGLO2}).

\begin{Thm}[{\cite[Theorem 9]{link-ext}}] \label{T:new-ext}
Fix a linking triple $\xSFL$.  Set  
$\dec\calh=\Ob(\dec\call)$, and assume it is closed under overgroups.  
Let
	\[ \calt_{\dec\calh}(\dec{S}) \Right5{\dec\delta} \dec\call 
	\Right5{\dec\pi} \dec\calf \]
be the structure functors for $\dec\call$.  
Set $\dec\Gamma=\Aut_{\dec\call}(\dec{S})$, and regard $\dec{S}$ as a subgroup of 
$\dec\Gamma$ via the inclusion 
$\dec\delta\:\calt_{\dec\calh}(\dec{S})\lra\dec\call$. 
Fix a finite group $G$ and an extension pair $U=(\Gamma,\tau)$ for 
$\dec\call$ and $G$, and choose 
$S\in\sylp{\Gamma}$.  Then there is a saturated fusion system 
$\calf$ over $S$ containing $\dec\calf$, and a transporter system 
$\calt$ associated to $\calf$ and containing $\dec\call$, such that the 
following hold: 
\begin{enuma} 
\item 
$\Ob(\calt)=\calh\defeq\{P\le{}S\,|\,P\cap{}\dec{S}\in\dec\calh\}$, 
and this set contains all subgroups of $S$ which are 
$\calf$-centric and $\calf$-radical.  Also, $\calt$ contains 
$\call_U$ as a full subcategory. 

\item The group $\Gamma$ can be identified with $\Aut_\calt(\dec{S})$ 
in a way so that the inclusion of $\dec\call$ in $\calt$ induces the inclusion 
of $\dec\Gamma=\Aut_{\dec\call}(\dec{S})$ in $\Gamma$.

\item For each $\gamma\in\Gamma$, 
$c_\gamma=\tau(\gamma)\in\Aut\typ^I(\dec\call)$.

\item $\dec\call$ is normal in $\calt$ (cf. Definition \ref{d:T0<|T}).

\item The space $|\call_U|$ is a deformation retract of $|\calt|$. The 
inclusion of geometric realizations $|\dec\call|\subseteq|\call_U|$ 
($\simeq|\calt|$) is homotopy equivalent to a regular covering space 
$X\Right2{}|\calt|$ with covering group $G\cong\Gamma/\dec\Gamma$.
\end{enuma}
\end{Thm}

\begin{proof}  The categories $\calt$ and $\calf$ will be constructed in 
Step 2, after preliminary constructions in Step 1. We show that $\calt$ is 
a transporter system in Steps 3 and 4, and prove that $\calf$ is saturated 
in Step 5.  Finally we prove (d) and (e) in Step 6.

Note that $\dec{S}=O_p(\dec\Gamma)\in\sylp{\dec\Gamma}$, since 
$\dec\Gamma/\dec{S}=\Aut_{\dec\call}(\dec{S})/\dec\delta_{\dec{S}}(\dec{S})$ has order 
prime to $p$ by Proposition \ref{L-prop}(\ref{x:fnorm}). For each 
$\gamma\in\Gamma$, let $c_\gamma\in\Aut(\dec{S})$ denote conjugation by 
$\gamma$ on $\dec{S}=O_p(\dec\Gamma)\nsg\Gamma$.  By the commutativity of 
triangle \ppluss\ in Definition \ref{d:U(L,G)}(a), this is the restriction 
to $\dec{S}$ of $\tau(\gamma)_{\dec{S}}\in\Aut(\dec\Gamma)$.  Hence by 
Proposition \ref{AutI}, $c_\gamma$ is fusion preserving (induces an 
automorphism of the category $\dec\calf$), and 
$\tau(\gamma)(P)=c_\gamma(P)$ for all $P\in\dec\calh$.  To simplify 
notation below, we write $\gamma(P)=\tau(\gamma)(P)$ to denote this action 
of $\gamma$ on $\dec\calh$.

\smallskip

\noindent\textbf{Step 1: }  Set $\call_1=\call_U$.  Thus by Definition 
\ref{d:U(L,G)}, $\Ob(\call_1)=\Ob(\dec\call)=\dec\calh$, and 
	\[ \Mor(\call_1) = \Mor(\dec\call)\times_{\dec\Gamma}\Gamma = 
	\bigl\{\5{\varphi,\gamma}\,\big|\, \varphi\in\Mor(\dec\call),~ 
	\gamma\in\Gamma\bigr\}. \]

We claim that 
	\beqq \textup{all morphisms in $\call_1$ are monomorphisms and 
	epimorphisms.} \label{T:e5} \eeqq
For any $\5{\varphi,\gamma}$, $\5{\varphi',\gamma'}$, and $\5{\psi,\eta}$ 
with appropriate domain and range,
	\begin{align*}  
	\5{\psi,\eta}\circ&\5{\varphi,\gamma} =
	\5{\psi,\eta}\circ\5{\varphi',\gamma'} \ \Longrightarrow\ 
	\5{\psi\circ\tau(\eta)(\varphi),\eta\gamma} = 
	\5{\psi\circ\tau(\eta)(\varphi'),\eta\gamma'} \\ 
	&\Longrightarrow\ \exists\,\lambda\in\dec\Gamma,\ 
	\eta\gamma=\lambda^{-1}\eta\gamma' \textup{ and } 
	\psi\circ\tau(\eta)(\varphi)=
	\psi\circ\tau(\eta)(\varphi')\circ\lambda \\
	&\Longrightarrow\  \gamma=(\eta^{-1}\lambda\eta)^{-1}\gamma',
	\textup{ and } 
	\varphi=\varphi'\circ\tau(\eta^{-1})(\lambda),
	\end{align*}
where the second equality in the last line holds since morphisms in 
$\dec\call$ are monomorphisms (Proposition \ref{L-prop}(\ref{x:monoepi})).  Also, 
$\tau(\eta^{-1})(\lambda)=\eta^{-1}\lambda\eta$ by the commutativity of 
\ppluss, so $\5{\varphi,\gamma}=\5{\varphi',\gamma'}$, and hence 
$\5{\psi,\eta}$ is a monomorphism.  The proof that morphisms are 
epimorphisms is similar.

Set $\Aut_{\Gamma}(\dec{S})=\{c_\gamma\in\Aut(\dec{S})\,|\,\gamma\in\Gamma\}$.  
Let $\calf_1$ be the smallest fusion system over $\dec{S}$ which contains 
$\dec\calf$ and $\Aut_{\Gamma}(\dec{S})$.  Define 
	$$ \pi_1 \: \call_1 \Right4{} \calf_1 $$
to be the identity on objects, while setting 
$\pi_1(\5{\varphi,\gamma}) =\dec\pi(\varphi)\circ{}c_\gamma$.  The proof that 
this is a functor (i.e., that it preserves composition) reduces to showing 
that the following square commutes
	\beq \vcenter{\xymatrix@C=50pt{
	\I12P \ar[r]^-{\dec\pi(\varphi)} \ar[d]_{c_\gamma} & \I12Q 
	\ar[d]^{c_\gamma} \\
	\I12\gamma(P) \ar[r]^-{\dec\pi(\tau(\gamma)(\varphi))} & \I12\gamma(Q) 
	}} \eeq
for each $\varphi\in\Mor_{\dec\call}(P,Q)$ and each $\gamma\in\Gamma$. By 
the above remarks, $c_\gamma\in\Aut(\dec{S},\dec\calf)$ is the restriction 
to $\dec{S}$ of $\tau(\gamma)_{\dec{S}}$.  By the last statement in 
Proposition \ref{AutI}, applied with $\alpha=\tau(\gamma)$ and 
$\beta=c_\gamma$, $\dec\pi\circ\tau(\gamma)= 
\widehat{c_\gamma}\circ\dec\pi$, where 
$\widehat{c_\gamma}\in\Aut(\dec\calf)$ is such that 
$\widehat{c_\gamma}(\dec\pi(\varphi))=c_\gamma\dec\pi(\varphi)c_\gamma^{-1}$.  
Hence the square commutes.  Since $\pi_1(\call_1)$ contains 
$\dec\calf|_{\dec\calh}$ and $\Aut_{\dec\Gamma}(\dec{S})$, and is closed 
under restrictions of morphisms to subgroups in $\dec\calh$ (Proposition 
\ref{L-prop}(\ref{x:uniqrest})), $\pi_1$ maps onto $\calf_1|_{\dec\calh}$.  

We regard $\dec\call$ as a subcategory of $\call_1$ by identifying each 
morphism $\varphi\in\Mor_{\dec\call}(P,Q)$ with 
$\5{\varphi,1}\in\Mor_{\call_1}(P,Q)$.  By construction, 
$\dec\pi=\pi_1|_{\dec\call}$.  For $P\le{}Q$ in $\dec\calh$, the inclusion morphism 
$\iota_{P,Q}=\dec\delta_{P,Q}(1)$ for $\dec\call$ is also considered as an 
inclusion morphism in $\call_1$.  The existence of restricted 
morphisms in $\dec\call$ (Proposition \ref{L-prop}(\ref{x:uniqrest})) carries 
over easily to the existence of restricted morphisms in $\call_1$, 
and they are unique by \eqref{T:e5}.

For all $P,Q\in\dec\calh$, define
	$$ (\delta_1)_{P,Q} \: N_{S}(P,Q) \Right5{} 
	\Mor_{\call_1}(P,Q) $$
by setting $(\delta_1)_{P,Q}(s)=\5{\iota_{{sPs^{-1}},Q},s}$. When 
$s\in{}N_{\dec{S}}(P,Q)$, 
$\5{\iota_{{sPs^{-1}},Q},s}=\5{\dec\delta_{P,Q}(s),1}$; and thus 
$(\delta_1)_{P,Q}$ extends the monomorphism $\dec\delta_{P,Q}$ from 
$N_{\dec{S}}(P,Q)$ to $\Mor_{\dec\call}(P,Q)$, under the identification of 
$\dec\call$ as a subcategory of $\call_1$.  To simplify the notation, we write 
$\delta_1(x)=(\delta_1)_{P,Q}(x)$ when $P$ and $Q$ are understood.

We claim that for all $P,Q\in\dec\calh$, $\psi\in\Mor_{\call_1}(P,Q)$, and 
$x\in{}P$, 
	\[ (\delta_1)_Q(\pi_1(\psi)(x))\circ\psi
	=\psi\circ(\delta_1)_P(x). \]
Set $\psi=\5{\varphi,\gamma}$, where $\gamma\in\Gamma$ and 
$\varphi\in\Mor_{\dec\call}(\gamma(P),Q)$.  Then 
	\begin{align*}  
	\psi\circ\delta_1(x) &= 
	\5{\varphi,\gamma}\circ\5{\Id_P,x} = \5{\varphi,\gamma x} 
	= \5{\varphi,c_\gamma(x)\gamma} 
	= \5{\varphi\circ\dec\delta(c_\gamma(x)),\gamma} \\
	&= \5{\dec\delta(\dec\pi(\varphi)(c_\gamma(x)))\circ\varphi,\gamma} 
	= \5{\dec\delta(\pi_1(\psi)(x))\circ\varphi,\gamma} \\
	&= \5{\dec\delta(\pi_1(\psi)(x)),1}\circ\5{\varphi,\gamma} 
	= \delta_1(\pi_1(\psi)(x))\circ\psi ~,
	\end{align*}
where the fifth equality holds by axiom (C) for the linking system 
$\dec\call$.  

We next show that morphisms in $\call_1$ have the following extension 
property: 
        \beqq
	\begin{split}  
	&\textup{$\forall\,P,Q\in\dec\calh$, 
	$\psi\in\Iso_{\call_1}(P,Q)$, and 
	$P',Q'\le{}\dec{S}$ 
	for which $P\nsg P'$, $Q\nsg Q'$, and}\\
	&\textup{$\psi\delta_1(P')\psi^{-1}\le
	\delta_1(Q')$, \ there is a unique 
	$\psi'\in\Mor_{\call_1}(P',Q')$ 
	where $\psi'|_{P,Q}=\psi$.}
	\end{split}
	\label{T:e4}
	\eeqq 
Set $\psi=\5{\varphi,\gamma}$, where 
$\varphi\in\Mor_{\dec\call}(\gamma(P),Q)$.  For all $x\in P'$, 
	\[ \5{\varphi,\gamma}\circ\5{\dec\delta(x),1}\circ 
	\5{\varphi,\gamma}^{-1} = 
	\5{\varphi\circ\tau(\gamma)(\dec\delta(x))\circ\varphi^{-1},1} =
	\5{\varphi\circ\dec\delta(c_\gamma(x))\circ\varphi^{-1},1}
	\in \delta_1(Q')~, \]
where $\tau(\gamma)(\dec\delta(x))=\dec\delta(c_\gamma(x))$ by the 
commutativity of \ppluss.  Thus 
$\varphi\dec\delta(\gamma(P'))\varphi^{-1}\le 
\dec\delta(Q')$, so $\varphi$ extends to 
$\varphi'\in\Mor_{\dec\call}(\gamma(P'),Q')$ by 
Proposition \ref{L-prop}(\ref{x:uniqext}).  Set 
$\psi'=\5{\varphi',\gamma}$.  Then 
$\psi'|_{P,Q}=\psi$ since $\tau(\gamma)(\iota_{{P},{P'}})= 
\iota_{{\gamma(P)},{\gamma(P')}}$ (i.e., $\tau(\gamma)$ sends inclusions 
to inclusions), and this proves \eqref{T:e4}.

\smallskip

\noindent\textbf{Step 2: }  Throughout the rest of the proof, for 
each $P\le{}S$, we set $\ch{P}=P\cap \dec{S}$. 
We next construct categories $\calt$ and 
$\calf_2$, both of which have object sets $\calh$, and which contain 
$\call_1$ and the restriction of $\calf_1$ to $\dec\calh$, respectively.  
Afterwards, we let $\calf$ be the fusion system over $S$ generated by 
$\calf_2$ and restrictions of morphisms.  

Let $\calt$ be the category with $\Ob(\calt)=\calh$, and where 
for all $P,Q\in\calh$,
	\beqq \Mor_{\calt}(P,Q) = 
	\bigl\{ \psi\in\Mor_{\call_1}(\ch{P},\ch{Q}) \,\big|\,
	\forall\, x\in{}P,\ \exists\, y\in{}Q \textup{ such that }
	\psi\circ\delta_1(x)=\delta_1(y)\circ\psi \bigr\}. 
	\label{T:e4a} \eeqq
If $\psi\in\Mor(\calt)$, then we denote the corresponding morphism in 
$\call_1$ by $\ch\psi$. Let 
	$$ \delta_{P,Q} \: 
	\underset{\subseteq{}N_{S}(\ch{P},\ch{Q})}{N_{S}(P,Q)} 
	\Right5{} 
	\underset{\subseteq\Mor_{\call_1}(\ch{P},\ch{Q})}{\Mor_{\calt}(P,Q)} $$
be the restriction of $(\delta_1)_{\ch{P},\ch{Q}}$.  Let $\calf_2$ be the 
category with $\Ob(\calf_2)=\calh$, and where
\begin{small}  
	$$ \Mor_{\calf_2}(P,Q) = 
	\bigl\{ \varphi\in\Hom(P,Q) \,\big|\, 
	\exists\,\psi\in\Mor_{\call_1}(\ch{P},\ch{Q}), \textup{ where }
	\psi\circ\delta_1(x)=\delta_1(\varphi(x))\circ\psi,\ \forall\, x\in{}P
	\bigr\}. $$
\end{small}
Define $\pi\:\calt\Right2{}\calf_2$ to be the identity on objects, and 
to send $\psi\in\Mor_{\calt}(P,Q)$ to the homomorphism $\pi(\psi)(x)=y$ 
whenever $\psi\circ\delta_1(x)=\delta_1(y)\circ\psi$ (uniquely defined by 
\eqref{T:e4a} and \eqref{T:e5}).  This is clearly a functor:  it is 
seen to preserve composition by juxtaposing the commutative squares which 
define $\pi$ on morphisms.  

Let $\calf$ be the fusion system over $S$ generated by $\calf_2$ and 
restriction of homomorphisms.  Since $\calh=\Ob(\calf_2)$ is closed under 
overgroups, $\calf_2$ is a full subcategory of $\calf$.  Since $\call_1$ 
is a full subcategory of $\calt$, 
$\Hom_{\calf_1}(P,Q)=\Hom_{\calf_2}(P,Q)$ for all $P,Q\in\dec\calh$.  If 
$P,Q\le{}\dec{S}$ are any subgroups and $\varphi\in\Hom_{\calf}(P,Q)$, then $\varphi$ 
is a composite of restrictions of morphisms in $\calf_2$, and hence (since 
$P\in\Ob(\calf_2)=\calh$ implies $\ch{P}\in\dec\calh$) a composite of 
restrictions of morphisms in $\calf_2$ (equivalently $\calf_1$) between 
subgroups in $\dec\calh$.  Thus $\varphi\in\Hom_{\calf_1}(P,Q)$; and we 
conclude that $\calf_1$ is also a full subcategory of $\calf$.  

\smallskip

\noindent\textbf{Step 3: } 
We next prove that 
	\beqq \textup{each $P\in\calh$ is $\calf$-conjugate to some 
	$P'\in\calh$ such that 
	$\delta_{\ch{P'}}(N_{S}(\ch{P'}))\in 
	\sylp{\Aut_\calt(\ch{P'})}$.} \label{T:e6} \eeqq

Fix $P\in\calh$. Let $\calp$ be the set of all $\dec{S}$-conjugacy 
classes $[\ch{Q}]$ of subgroups $\ch{Q}\in\ch{P}{}^{\dec\calf}$ (recall 
$\ch{P}=P\cap\dec{S}$) which are fully normalized in $\dec\calf$.  (If 
$\ch{Q}$ is fully normalized in $\dec\calf$, then so is every subgroup in 
$[\ch{Q}]$.) By Lemma \ref{Pfn:p'}, $|\calp|$ is finite and prime to $p$.

We claim that in general, for each $\gamma\in\Gamma$ and each 
$\ch{Q},\ch{R}\in\dec\calh$, 
	\begin{align} 
	\textup{$\ch{Q}$ and $\ch{R}$ $\dec\calf$-conjugate} 
	\quad&\iff\quad
	\textup{$\gamma(\ch{Q})$ and $\gamma(\ch{R})$ $\dec\calf$-conjugate} 
	\label{T:e3a} \\
	\textup{$\ch{Q}$ and $\ch{R}$ $\dec{S}$-conjugate} 
	\quad&\iff\quad
	\textup{$\gamma(\ch{Q})$ and $\gamma(\ch{R})$ $\dec{S}$-conjugate.} 
	\label{T:e3b} 
	\end{align}
The first holds since $\gamma$ acts on $\dec\call$ and hence on $\dec\calf$ as 
a group of automorphisms (Proposition \ref{AutI}), and the second since 
$\dec{S}\nsg\Gamma$.

Let $\Gamma'\subseteq\Gamma$ be the subset of those 
$\gamma\in\Gamma$ such that $\gamma(\ch{P})\in\ch{P}^{\dec\calf}$. 
Then for $\gamma_1,\gamma_2\in\Gamma'$, 
$\gamma_1\gamma_2(\ch{P})\in\gamma_1(\ch{P})^{\dec\calf}$ by \eqref{T:e3a} 
and since $\gamma_2(\ch{P})\in\ch{P}^{\dec\calf}$, and hence 
$\gamma_1\gamma_2\in\Gamma'$.  By \eqref{T:e3b} and 
since each $\gamma\in\Gamma$ acts on $\dec{S}$ via the fusion preserving 
automorphism $c_\gamma\in\Aut(\dec{S})$ as shown above, $\gamma$ permutes 
the $S$-conjugacy classes of subgroups which are fully normalized in 
$\dec\calf$.  Thus each element of $\Gamma'$ permutes the set $\calp$.  

Fix $S'\in\sylp{\Gamma'}$.  Let $\eta\in\Gamma$ be such 
that $S''\defeq\eta{}S'\eta^{-1}\le{}S$.   Since $\calp$ has 
order prime to $p$ by Lemma \ref{l:centric}, there is some 
$[\ch{Q}]\in\calp$ fixed by $S'$.  In other words, for each 
$\gamma\in{}S'$, $\gamma(\ch{Q})$ is $\dec{S}$-conjugate to 
$\ch{Q}$.  So by \eqref{T:e3b}, for each $s=\eta\gamma\eta^{-1}\in{}S''$ 
(where $\gamma\in{}S'$), $s(\eta(\ch{Q}))$ is $\dec{S}$-conjugate to 
$\eta(\ch{Q})$.  Set $\ch{R}=\eta(\ch{Q})$. Then each coset in 
$S''/\dec{S}$ contains some element $s$ which normalizes $\ch{R}$, i.e., the 
obvious homomorphism $ N_{S''}(\ch{R})\to S''/\dec{S}$ is onto with kernel 
$N_{\dec{S}}(\ch{R})$. Since $S$ and $\dec{S}$ have the same 
identity component,
	\beqq |\pi_0(N_{S}(\ch{R}))| \ge |\pi_0(N_{S''}(\ch{R}))| 
	=|\pi_0(N_{\dec{S}}(\ch{R}))|\cdot|S''/\dec{S}| = 
	|\pi_0(N_{\dec{S}}(\ch{Q}))|\cdot|S'/\dec{S}|~. \label{T:e3c} \eeqq

Since $\Gamma'$ is the subgroup of elements of $\Gamma$ which send 
$\ch{P}$ to a subgroup in its $\dec\calf$-conjugacy class, 
	\beqq |\pi_0(\Aut_\calt(\ch{R}))| = |\pi_0(\Aut_\calt(\ch{P}))| = 
	|\pi_0(\Aut_{\dec\call}(\ch{P}))|\cdot|\Gamma'/\dec\Gamma| = 
	|\pi_0(\Aut_{\dec\call}(\ch{Q}))|\cdot|\Gamma'/\dec\Gamma|~. 
	\label{T:e3d} \eeqq
Since $\ch{Q}$ is fully normalized in $\dec\calf$, $S'\in\sylp{\Gamma'}$, 
and $\dec{S}\in\sylp{\dec\Gamma}$, \eqref{T:e3c} and \eqref{T:e3d} imply 
that $\delta_{\ch{R}}(N_{S}(\ch{R}))$ is a Sylow $p$-subgroup 
of $\Aut_\calt(\ch{R})$.  

Choose any $\psi\in\Iso_\calt(\ch{P},\ch{R})$.  Then 
$\psi\delta_{\ch{P}}(N_{S}(\ch{P}))\psi^{-1}$ is a $p$-subgroup 
of $\Aut_\calt(\ch{R})$.  Choose $\chi\in\Aut_\calt(\ch{R})$ such that 
$(\chi\psi)\delta_{\ch{P}}(N_{S}(\ch{P}))(\chi\psi)^{-1}\le 
\delta_{\ch{R}}(N_{S}(\ch{R}))$.  By definition of the category 
$\calt$, $\chi\psi$ extends to a morphism 
$\widebar{\psi}\in\Mor_\calt(P,N_{S}(\ch{R}))$.  Set 
$P'=\pi(\widebar{\psi})(P)$.  Then $\ch{P'}=\ch{R}$, $P'$ is 
$\calf$-conjugate to $P$, and $P'\in\calh$ since 
$\ch{P'}\in\dec\calh$ ($\dec\calh$ is invariant under 
$\dec\calf$-conjugacy).  This finishes the proof of \eqref{T:e6}.

\smallskip

\noindent\textbf{Step 4: } We are now ready to show that $\calt$ is a 
transporter system.  
For each $P\in\dec\calh$, set
	\[ E(P) = \Ker\bigl[\Aut_\calt(P) \Right4{\pi_P} \autf(P) \bigr] 
	~. \]
For each $P,Q\in\Ob(\calt)$, $E(P)$
acts on $\Mor_\calt(P,Q)$ by right composition and $E(Q)$ by left
composition. Both actions are free since all morphisms in $\call_1$ (hence
in $\calt$) are monomorphisms and epimorphisms by \eqref{T:e5}.  We claim
that $\pi_{P,Q}$ is the orbit map of the action of $E(P)$ on
$\Mor_\calt(P,Q)$.  Since every morphism in $\calt$ (and also by definition in $\calf$)
factors uniquely as the composite of an isomorphism followed by an
inclusion, it suffices to prove this when $P$ and $Q$ are
$\calf$-conjugate.  It thus suffices to prove it when $P=Q$, and this holds
by definition of $E(P)$.  

This proves axiom (A2). Axioms (A1) and (B) hold by construction, and (C) 
holds by definition of the functor $\pi\:\calt\to{}\calf$.  It 
remains to prove axioms (I) -- (III).

Fix $P\in\calh$ such that $\delta_{\ch{P}}(N_{S}(\ch{P}))\in\sylp{\Aut_\calt(\ch{P})}$.  By 
\eqref{T:e6}, every subgroup in $\calh$ is $\calf$-conjugate to some such 
$P$.  Write $G=\Aut_\calt(\ch{P})$, $T=\delta_{\ch{P}}(N_{S}(\ch{P}))$, and 
$P'=\delta_{\ch{P}}(P)$ for short, where $\delta_{\ch{P}}$ is 
injective by construction.  Thus $P'\le{}T\in\sylp{G}$.  Fix 
$R\in\sylp{N_G(P')}$, and choose $\alpha\in{}G$ such that 
$\alpha{}R\alpha^{-1}\le{}T$.  Then 
	$\alpha{}R\alpha^{-1}\in
	\sylp{N_G(\alpha{}P'\alpha^{-1})}$, and so 
\[\alpha P'\alpha^{-1}\le\alpha{}R\alpha^{-1} = N_T(\alpha{}P'\alpha^{-1}).\] 
Set 
$Q=\delta_{\ch{P}}^{-1}(\alpha P'\alpha^{-1})\le\delta_{\ch{P}}^{-1}(T) = N_{S}(\ch{P})$.  Then $Q$ is $\calf$-conjugate to $P$, $N_{S}(Q)\le N_{S}(\ch{P})$, and so
	\begin{align*}  
	N_T(\alpha P'\alpha^{-1}) &= 
	N_{\delta_{\ch{P}}(N_{S}(\ch{P}))}(\delta_{\ch{P}}(Q)) = \delta_{\ch{P}}(N_{S}(Q)),~ \textup{and} \\
	N_G(\alpha P'\alpha^{-1}) &= 
	N_{\Aut_\calt(\ch{P})}(\delta_{\ch{P}}(Q)) = \Aut_\calt(Q)~. 
	\end{align*}
Thus $\delta_Q(N_{S}(Q))\in\sylp{\Aut_\calt(Q)}$, and this proves axiom (I).

We next claim that 
	\beqq \begin{split} 
	P\nsg P'\le{}S,\quad Q\nsg Q'\le{}S,\quad 
	&\psi\in\Iso_{\calt}(P,Q), \quad
	\psi(\delta_P(P'))\psi^{-1}\le\delta_Q(Q')
	\\ &\implies\quad \exists~
	\psi'\in\Mor_\calt(P', Q') 
	\textup{ with } \psi'|_{P,Q}=\psi. 
	\end{split} \label{T:e7a} \eeqq
Set $\ch\psi=\psi|_{\ch{P},\ch{Q}}=\5{\varphi,\gamma}$, where 
$\varphi\in\Iso_{\dec\call}(\gamma(\ch{P}),\ch{Q})$ and $\gamma\in\Gamma$.  
By Proposition \ref{L-prop}(\ref{x:uniqext}), $\varphi$ extends in a 
unique way to $\gamma(\ch{P'})$. Hence $\ch{\psi}$ extends to some unique 
$\ch{\psi'}\in\Mor_{\call_1}(\ch{P'},\ch{Q'})$.  By definition of morphisms 
in $\calt$ and the original hypothesis on $\psi$, $\ch{\psi'}$ extends to a 
morphism $\psi'\in\Mor_\calt(P',Q')$ which extends $\psi$.  This proves 
\eqref{T:e7a}, and thus proves axiom (II).

It remains to prove axiom (III).  Fix $P_1\le P_2\le P_3\le\cdots$ in 
$\Ob(\calt)$ and $\psi_i\in\Mor_\calt(P_i,S)$ such that for all 
$i\ge1$, $\psi_i=\psi_{i+1}\circ\delta_{P_i,P_{i+1}}(1)$.  Set 
$P=\bigcup_{i=1}^\infty P_i$.  Let $\varphi_i\in\Mor(\dec\call)$ and 
$\gamma_i\in\dec\Gamma$ be such that $\ch{\psi}_i=\5{\varphi_i,\gamma_i}$ 
for each $i$.  Since $|\Gamma/\dec\Gamma|<\infty$, we can assume, after 
passing to a subsequence, that the $\gamma_i$ all lie in the same coset 
$\dec\Gamma\gamma$.  If $\gamma_i=\nu_i\gamma$ for some 
$\nu_i\in\dec\Gamma$, then $\ch{\psi}_i=\5{\varphi_i\nu_i,\gamma}$; we can 
thus assume $\gamma_i=\gamma$ for each $i$.  But then 
$\varphi_i=\varphi_{i+1}|_{P_i,P_{i+1}}$ for each $i$, so by Proposition 
\ref{L-prop}(\ref{x:axiomIII}) applied to $\dec\call$, there is 
$\varphi\in\Mor_{\dec\call}(P,\dec{S})$ such that 
$\varphi|_{P_i,\dec{S}}=\varphi_i$ for each $i$.  Set 
$\ch{\psi}=\5{\varphi,\gamma}\in\Mor_{\calt}(\ch{P},\dec{S})$; then 
$\ch{\psi}|_{\ch{P}_{i}}=\ch{\psi}_{i}\defeq\psi_i|_{\ch{P_i}}$ for each 
$i$.

Fix $g\in{}P$, and let $i$ be such that $g\in{}P_i$.  Set 
$h=\pi(\psi_i)(g)\in S$.  Then 
$\delta_{S}(h)\circ\psi_i=\psi_i\circ\delta_{P_i}(g)$ by 
axiom (C) for $\calt$, so 
$(\delta_1)_{\dec{S}}(h)\circ\ch\psi=\ch\psi\circ(\delta_1)_{\ch{P}_i}(g)$ 
by \eqref{T:e5} (the morphisms are epimorphisms in $\call_1$).  Hence by 
definition of $\calt$, there is a unique morphism 
$\psi\in\Mor_\calt(P,S)$ such that $\psi|_{\ch{P}}=\ch{\psi}$.  By 
the uniqueness of extensions again, $\psi|_{P_i}=\psi_i$ for each $i$.

\smallskip


\noindent\textbf{Step 5: }  We are now ready to show that $\calf$ is 
saturated.  By Theorem \ref{H-sat-gen}, it suffices to prove the following 
statements.
\begin{enumi} 
\item The group $\Aut_{\calf}(S_0)$ is finite.

\item For all $P\le S$, $Q\le S_0$, and 
$\varphi\in\Hom_{\calf}(P,Q)$, there exists 
$\widebar{\varphi}\in\Hom_{\calf}(P\cdot C_{S}(P)_0,S_0)$ 
such that $\widebar{\varphi}|_P=\varphi$.

\item $\calh$ is closed in $\Sub(S)$, and $\calf$ is 
$\calh$-generated and $\calh$-saturated.

\item For all $P\in\calh$ and $P\le Q\le P\8$, $Q\in\calh$.

\item\label{T:e8} 
If $P\in\calf^\bullet$ is $\calf$-centric and $P\notin\calh$, then there is 
$Q\in{}P^{\calf}$ such that 
	\[ O_p(\Out_{\calf}(Q)) \cap \Out_{S}(Q) \ne 1. \]
\end{enumi}
Point (i) holds by construction (and since $\dec\calf$ is saturated).  
By Step 4 and Proposition \ref{T=>Fsat}, $\calf$ is 
$\calh$-saturated; i.e., it satisfies the saturation axioms for subgroups 
in $\calh$.  It is also $\calh$-generated by definition:  each morphism 
in $\calf$ is a composite of restrictions of morphisms between subgroups in 
$\calh$.  Since $\dec\calh$ is closed under overgroups and
$\dec\calf$-conjugacy, $\calh$ is closed under overgroups and
$\calf$-conjugacy by definition, and this finishes the proofs of (iii) 
and (iv).

If $P\le S$, $Q\le S_0$, and $\varphi\in\Hom_{\calf}(P,Q)$, 
then $P\le\dec{S}$ since $\dec{S}$ is strongly closed in $\calf$ (and 
$\dec{S}_0=S_0$).  By definition of $\calf_1$, 
$\varphi=\chi\circ\psi$ for some $\chi\in\Aut_{\calf}(S_0)$ 
($\chi=\pi(\gamma)|_{S_0}$ for some 
$\gamma\in\Gamma=\Aut_\calt(\dec{S})$) and some 
$\psi\in\Mor(\dec\calf)$.  Set $P'=\psi(P)\le S_0$, and choose 
$R\in{}P^{\dec\calf}=P'{}^{\dec\calf}$ which is receptive in $\dec\calf$.  
By \cite[Lemma 2.4(a)]{BLO3}, $R\le S_0$, and there is 
$w\in\Aut_{\dec\calf}(S_0)$ such that $w(P')=R$.  Since $R$ is 
receptive, $w\circ\psi$ extends to 
$\widebar{\psi}\in\Hom_{\dec\calf}(P\cdot C_{S}(P)_0,S_0)$, so 
$\varphi$ extends to $\chi\circ w^{-1}\circ\widebar{\psi}$.  This proves 
(ii).

It remains to prove (v).  
Let $\calk$ be the set of all $P\in\calf\8$ such that the 
saturation axioms hold for subgroups $\calf$-conjugate to $P$ and all 
of their overgroups.  Since $\calh$ is closed under overgroups and 
$\calf$-conjugacy, $\calk\supseteq\calh\8$.  
Set $\calk'=\calf\8{\sminus}\calk$, and let $\ch\calk'$ be the 
set of all $\ch{P}=P\cap{}\dec{S}$ for $P\in\calk'$.  We will show that for 
all $P\le{}S$, 
	\beqq \textup{$P\in\calk$ or $\ch{P}$ maximal in $\ch\calk'$ 
	$\Longrightarrow$ \eqref{T:e8} holds for $P$} 
	\label{T:e9} \eeqq
Having done that, we will prove that $\calk'=\emptyset$, so 
$\calk\supseteq\calf\8$, and hence (using \eqref{T:e9}) that 
\eqref{T:e8} holds for all $P\in\calf\8$.


We first show that
	\beqq 
	\textup{$P\in\calk$ or $\ch{P}$ maximal in $\ch\calk'$} 
	~\implies~ \exists\, Q\in P^{\calf}
	\textup{ such that $\ch{Q}$ is fully normalized in $\dec\calf$\,.}
	\label{T:e11} \eeqq
If $\ch{P}$ is fully normalized, we are done, so assume otherwise.  Let 
$\ch{P'}$ be $\dec\calf$-conjugate to $\ch{P}$ and fully normalized in 
$\dec\calf$.  Since $\dec\calf$ is saturated, $\ch{P'}$ is fully 
automized and receptive, and so by \ref{l:RS}(c), there is 
$\rho\in\Hom_{\dec\calf}(N_{\dec{S}}(\ch{P}),N_{\dec{S}}(\ch{P'}))$ such 
that $\rho(\ch{P})=\ch{P'}$.  Clearly, $\ch{P}<{}\dec{S}$, so 
$N_{\dec{S}}(\ch{P})>{}\ch{P}$.  If $\ch{P}$ is maximal in $\ch\calk'$, 
then by Lemma \ref{reduce2dots}(b), the saturation axioms hold for all $Q$ 
such that $\ch{Q}>\ch{P}$.  So whether $\ch{P}\in\calk$ or $\ch{P}$ is 
maximal in $\ch\calk'$, the saturation axioms hold for 
$N_{\dec{S}}(\ch{P})$, $N_{\dec{S}}(\ch{P'})$, and all subgroups of 
$S$ which contain them.  

Set $R=N_{\dec{S}}(\ch{P})$, and let 
$K=\{\alpha\in\Aut(R)\,|\,\alpha(\ch{P})=\ch{P}\}$.  Set $R'=\rho(R)$ and 
$K'=\9{\rho}K$.  Choose $R''\le{}\dec{S}$ and 
$\tau\in\Iso_{\calf}(R,R'')$ such that $R''$ is fully 
$\9{\tau}K$-normalized in $\calf$, and set $\ch{P}''=\tau(\ch{P})$ and 
$K''=\9{\tau}K$.  Thus $K'$ and $K''$ are the groups of automorphisms of 
$R'$ and $R''$, respectively, which send $\ch{P}'$ and $\ch{P}''$ to themselves.  
By Lemma \ref{l:RS-K}(b), there are automorphisms 
$\chi,\chi'\in\Aut_{\calf}^{K''}(R'')$, and morphisms 
	\[ \widebar{\tau}\in\Hom_{\calf}\bigl(N_{S}^K(R),
	N_{S}^{K''}(R'')\bigl) \qquad\textup{and}\qquad
	\widebar{\rho}\in\Hom_{\calf}\bigl(N_{S}^{K'}(R'),
	N_{S}^{K''}(R'')\bigl) \]
such that $\widebar{\tau}|_{R}=\chi\tau$ and 
$\widebar{\rho}|_{R'}=\chi'\tau\rho^{-1}$.  

We claim that 
	\beqq |N_{\dec{S}}(\ch{P})| < |N_{\dec{S}}^{K'}(R')| 
	\le |N_{\dec{S}}^{K''}(R'')| \le |N_{\dec{S}}(\ch{P''})| 
	\,. \label{T:e10} \eeqq
Since $\ch{P}$ is not fully normalized in $\dec\calf$, 
$R'=\rho(N_{\dec{S}}(\ch{P}))<{}N_{\dec{S}}(\ch{P'})$, and hence 
	\[ R'< N_{N_{\dec{S}}(\ch{P'})}(R') = N_{\dec{S}}^{K'}(R')  \]
by Lemma \ref{N_Q(P)>P}.  
This proves the first inequality in \eqref{T:e10}.  The next one holds 
since $\bar\rho$ sends $N_{\dec{S}}^{K'}(R')$ into 
$N_{\dec{S}}^{K''}(R'')$, and the last since all elements of 
$N_{\dec{S}}^{K''}(R'')$ normalize $\ch{P''}$. Thus $P''$ is 
$\calf$-conjugate to $P$ and 
$|N_{\dec{S}}(\ch{P''})|>|N_{\dec{S}}(\ch{P})|$. If $\ch{P''}$ is not fully 
normalized in $\dec\calf$, then since $P^{\calf}$ contains finitely 
many $S$-conjugacy classes (Lemma \ref{P^F/S-finite}), we can 
repeat this procedure, until we find a subgroup $Q$ which satisfies 
\eqref{T:e11}.

We are now ready to prove \eqref{T:e9}.  Assume $P$ is $\calf$-centric 
and $P\notin\calh$ (otherwise the statement is empty), and let $Q$ be 
$\calf$-conjugate to $P$ such that $\ch{Q}$ is fully 
$\dec\calf$-normalized.  Thus $\ch{P}\notin\dec\calh=\Ob(\dec\call)$.  Since by 
definition, a linking system must contain all centric and radical 
subgroups with respect to the underlying fusion system, either $\ch{P}$ and 
$\ch{Q}$ are not $\dec\calf$-centric or they are not $\dec\calf$-radical.  
If $\ch{Q}$ is not $\dec\calf$-centric, then there is 
$g\in{}C_{\dec{S}}(\ch{Q}){\sminus}\ch{Q}$ (since $\ch{Q}$ is fully 
centralized). If $\ch{Q}$ is not $\dec\calf$-radical, then
$O_p(\Aut_{\dec\calf}(\ch{Q}))>\Inn(\ch{Q})$ and is
contained in the Sylow subgroup $\Aut_{\dec{S}}(\ch{Q})$
($\ch{Q}$ is fully normalized), and thus there is
$g\in{}N_{\dec{S}}(\ch{Q}){\sminus}\ch{Q}$ such that
$c_g\in{}O_p(\Aut_{\dec\calf}(\ch{Q}))$.  In either case, 
	\[ g\in{}Q' \defeq \bigl\{g\in{}N_{\dec{S}}(\ch{Q}) \,\big|\, 
	c_g\in{}O_p(\Aut_{\dec\calf}(\ch{Q})) \bigr\}
	\qquad\textup{and}\qquad g\notin \ch{Q} ~, \]
and hence $Q'>{}\ch{Q}$.  Also, $Q$ normalizes $Q'$ and $QQ'>{}Q$, so 
$N_{QQ'}(Q)>{}Q$, and there is $x\in{}Q'{\sminus}Q$ such that 
$x\in{}N_{S}(Q)$. For any such $x$, $c_x\notin\Inn(Q)$ since 
$C_{S}(Q)\le{}Q$ (recall $P$ is $\calf$-centric), and 
$c_x|_{\ch{Q}}$ is in $O_p(\Aut_{\dec\calf}(\ch{Q}))$.  Also, 
since $x\in{}\dec{S}\nsg S$, $[x,S]\le \dec{S}$, and so 
$c_x$ induces the 
identity on $Q/\ch{Q}$. Thus $c_x$ is in the subgroup 
	\[ U = \bigl\{\alpha\in\autf(Q)\,\big|\,
	\alpha|_{\ch{Q}}\in{}O_p(\Aut_{\dec\calf}(\ch{Q})),
	\textup{ $\alpha$ induces the identity on $Q/\ch{Q}$}
	\bigr\} \nsg \autf(Q). \]
Since the group of all $\alpha\in\Aut(Q)$ which induce the identity on 
$\ch{Q}$ and on $Q/\ch{Q}$ is discrete $p$-toral by Lemma \ref{Aut(PrelQ)}, 
and is contained in $U$ with $p$-power index, $U$ is a nontrivial normal 
discrete $p$-toral subgroup of $\autf(Q)$, and $U\le 
O_p(\autf(Q)$.  Since $x\in N_{\dec{S}}(Q)\sminus Q$, $c_x$ 
represents a nontrivial element of $\Out_{S}(Q)\cap 
O_p(\outf(Q))$, so \eqref{T:e8} holds for $P$, and the proof 
of \eqref{T:e9} is complete.

We want to show that $\calk'=\emptyset$.  Assume otherwise; then 
$\ch\calk'\ne\emptyset$ since $P\in\calk'$ implies $\ch{P}\in\ch{\calk}'$.  
Choose $Q$ to be maximal in $\ch{\calk'}$, and choose $P$ to be maximal 
among those $P\in\calk'$ such that $\ch{P}=Q$.  Then $P$ is also maximal in 
$\calk'$.  So by Lemma \ref{reduce2dots}(b)), $P$ is maximal among 
subgroups not satisfying the saturation axioms.  By Lemma \ref{P**}, this 
maximality of $P$ implies that \eqref{T:e8} does not hold for $P$.  Since 
this contradicts \eqref{T:e9}, we now conclude that $\calk'=\emptyset$, and 
hence (by \eqref{T:e9}) that \eqref{T:e8} holds for all $P\in\calf\8$.  

Thus $\calf$ is saturated.  Also, \eqref{T:e8} implies that $\calh$ 
contains all subgroups which are $\calf$-centric and $\calf$-radical.  

\smallskip

\noindent\textbf{Step 6: } By \cite[Corollary 3.5]{BLO3}, $\calf\8$ 
contains all subgroups which are $\calf$-centric and 
$\calf$-radical, so they are all contained in $\calh$ by point (v) 
in Step 5.  Point (a) holds by this together with the definition of 
$\calh=\Ob(\calt)$.  Point (b) holds by the definition of $\calt$ in 
Step 2, and (c) holds by the definition of composition (of morphisms 
between subgroups in $\dec\calh$) in Step 1.  

Condition (i) in Definition \ref{d:T0<|T} (for the inclusion 
$\dec\call\subseteq\calt$) holds by (a), and since $\dec{S}$ is strongly closed 
in $\calf$ by construction.  Condition (ii) (the Frattini condition) 
holds by the construction in Step 1, and (iii) (invariance of $\dec\call$ under 
$\Aut_\calt(\dec{S})$-conjugacy) holds by (c).  Thus $\dec\call\nsg\calt$.


Let $r\:\calt\Right2{}\calt$ be the retraction $r(P)=\ch{P}$ and 
$r(\varphi)=\ch{\varphi}$ with image $\call_1$.  There is a natural 
transformation of functors $r\Right2{}\Id_\calt$ which sends an object $P$ 
to the inclusion $\iota_{{\ch{P}},P}$.  Hence $|r|\simeq\Id_{|\calt|}$, and so 
$|\call_U|$ is a deformation retract of $|\calt|$.  Also, the inclusion of 
$|\dec\call|$ into $|\call_U|$ is homotopy equivalent to a regular covering 
with covering group $G\cong\Gamma/\dec\Gamma$ by \cite[Proposition 
A.4]{OV1}, and this finishes the proof of (e).
\end{proof}


\bigskip


\section{The category $\AUT(\call)$}

Fix a $p$-local compact group $\SFL$.  Let $\AUT(\call)$ be the groupoid 
with object set $\Aut\typ^I(\call)$ and with morphisms the natural 
transformations. Since a natural transformation $\eta\:\alpha\Right1{}\beta$ 
is determined by $\eta_S\in\Aut_\call(S)$ (since $\alpha(S)=\beta(S)=S$), 
morphisms can be described as
	\begin{align} 
	\Mor_{\AUT(\call)}(\alpha,\beta) 
	&= \{\,\chi\in \Aut_{\call }(S) \,|\, 
	\chi\alpha(\varphi)=\beta(\varphi)\chi, \textup{ all 
	$\varphi\in\Mor(\call)$} \,\} \tag*{} \\
	&= \{\,\chi\in \Aut_{\call }(S) \,|\, \beta=c_\chi\circ\alpha \,\}.
	\label{e:Mor(Aut)} \end{align}
Thus $\chi\in\Aut_\call(S)$ corresponds to a natural transformation $\eta$ 
with $\chi=\eta_{S}$. With this notation, composition of morphisms takes 
the form
	$$ \xymatrix{\alpha\ar@/_1pc/[rr]_{\omega\circ\chi} \ar[r]^{\chi} 
	& \beta \ar[r]^{\omega} & \gamma } $$
because $\gamma=c_\omega\circ\beta=c_\omega\circ c_\chi\circ\alpha = 
c_{\omega\circ\chi}\circ\alpha$.  Here, $\omega\circ\chi$ is 
the composite of $\omega$ and $\chi$ in $\Aut_{\call}(S)$.

Composition in $\Aut\typ^I(\call)$ gives $\AUT(\call)$ the 
structure of a discrete strict monoidal category, where 
$\alpha\cdot \beta = \alpha\circ\beta$, and where 
	\begin{equation}\label{composition}
	(\alpha\Right1{\chi}{}\beta)\cdot(\alpha'\Right1{\chi'}{}\beta')=
	\bigl(\alpha\circ\alpha'\RIGHT7{\beta(\chi')\circ\chi}
	{=\chi\circ\alpha(\chi')} \beta\circ\beta'\bigr) \,.
	\end{equation} 
This structure makes the nerve $\caln\AUT(\call)$ into a simplicial 
group, and its geometric realization $|\AUT(\call)|$ becomes a 
topological group. The projection $\AUT(\call) 
\Right2{}\Out\typ(\call)$ induces a map of simplicial groups 
$\caln\AUT(\call)\Right2{pr}{}\Out\typ(\call)$, where now 
$\Out\typ(\call)$, by abuse of language, denotes the discrete simplicial 
group with vertex set $\Out\typ(\call)$. This projection is explicitly 
given by 
	$$ \pr\bigl(\alpha_0\Right1{}{}\alpha_1\Right1{}{}\dots 
	\Right1{}{}\alpha_n\bigr) = [\alpha_0]
	\qquad( = [\alpha_1] = \dots = [\alpha_n])\,,$$
and it sends $\pi_0(|\AUT(\call)|)$ isomorphically onto $\Out\typ(\call )$. 

The evaluation functor $\AUT(\call)\times\call \Right2{}\call$ induces an 
action of the simplicial group $\caln\AUT(\call)$ on the simplicial set 
$\caln\call$ as follows. In dimension 0, 
$\caln_0\AUT(\call)=\Aut\typ^I(\call)$ acts on  $\caln_0\call=\Ob(\call)$ 
in the obvious way. For $\chi\in\caln_1\AUT(\call)$ and 
$\varphi\in\Mor_\call(P,Q)=\caln_1\call$, 
	\begin{equation}\label{e:action} 
	\bigl( \alpha\Right3{\chi}\beta\bigr) 
	\bigl(P\Right3{\varphi}Q\bigr) = \bigl(\alpha(P) 
	\RIGHT9{\beta(\varphi)\circ\chi_P}{=\chi_Q\circ\alpha(\varphi)} 
	\beta(Q) \bigr), \end{equation}
and this extends naturally to higher dimensional sequences. 
The simplicial action  induces an action of the topological group 
$|\AUT(\call)|$ on the space $|\call|$.

Our aim is to describe   maps  $BG\Right1{}B|\caln\AUT(\call)|$ via 
\emph{twisting functions} (see, e.g., \cite[Definition 3.14]{Curtis}). 
We need to show that extensions of the type constructed in Theorem 
\ref{T:new-ext} realize certain types of topological fibre bundles, and the 
relevant obstruction theory is encoded in the simplicial equalities that 
characterise the twisting functions.

\newcommand{\sbt}
{\,\begin{picture}(-1,1)(-1,-3)\circle*{3}\end{picture}\,}

For any (discrete) group $G$, let $\cale(G)$ denote the category whose 
objects are the elements of $G$, with a unique morphism between each 
pair of objects.  Let $\calb(G)$ be the category with a single 
object $\sbt$\,, and with $G$ as the automorphism group of that object. 
Then $G$ acts on $\cale(G)$ by translation, and the quotient category 
can be identified with $\calb(G)$. The geometric realizations of 
$\cale(G)$ and $\calb(G)$ are the universal contractible free $G$-space 
$EG$, and the classifying space $BG$, respectively. 

For a discrete group $G$, we will use the simplicial set 
$\caln\calb(G)^{\mathrm{op}}$ as a model for $BG$. This allows us to 
conveniently denote simplices in $\caln_n\calb(G)^{\mathrm{op}}$  by the 
usual bar notation 
	\[\7g = [g_1|g_2|\cdots|g_n]\defeq \bigl( \sbt \Left3{g_1} 
	\sbt \Left3{g_2} \sbt \cdots \sbt \Left3{g_n} \sbt 
	\bigr),\] 
where $g_i\in G$. We will generally omit the superscript ``op'' from the 
notation. More generally, for any small category $\calc$, we will consider 
the nerve of the opposite category $\calc\op$ as a model for the nerve of 
$\calc$.  Thus, for example,  we consider face maps as
	\begin{align*} 
	&d_0\bigl(c_0\Left2{f_1}c_1\Left2{f_2}c_2\cdots \Left2{f_n}c_n\bigl) = 
	\bigl(c_1\Left2{f_2}c_2\cdots \Left2{f_n}c_n\bigl) \,, \\
	&d_1\bigl(c_0\Left2{f_1}c_1\Left2{f_2}c_2\cdots \Left2{f_n}c_n\bigl) = 
	\bigl(c_0\Left2{f_1\circ f_2}c_2\cdots \Left2{f_n}c_n\bigl) \,, 
	\end{align*}
etc.

To each simplicial group $K$, one associates a simplicial set 
$\widebar{W}(K)$ (cf. \cite[Definition 3.20]{Curtis} or 
\cite[\S\,V.4]{GJ}), by setting $\widebar{W}_0(K)=*$,
	\[ \widebar{W}_n(K) = K_{n-1} \times K_{n-2} \times \cdots K_0 
	\qquad\textup{for $n>0$,} \]
and with face and degeneracy maps 
	\beqq\begin{split} 
	d_i(\kappa_{n-1},\ldots,\kappa_0) &= \begin{cases} 
	(\kappa_{n-2},\ldots,\kappa_0) & \textup{if $i=0$} \\
	(d_{i-1}\kappa_{n-1},\ldots,(d_0\kappa_{n-i})\kappa_{n-i-1},\ldots,\kappa_0) & 
	\textup{if $0<i<n$} \\
	( d_{n-1}\kappa_{n-1},\ldots,d_1\kappa_1) & \textup{if $i=n$}
	\end{cases} \\
	s_i(\kappa_{n-1},\ldots,\kappa_0) &= (s_{i-1}\kappa_{n-1},\ldots, 
	s_0\kappa_{n-i},1,\kappa_{n-i-1},\ldots,\kappa_0) \qquad \textup{all $0\le 
	i\le n$}
	\end{split} \label{e:Wbar} \eeqq
for $\kappa_i\in{}K_i$. Then $|\widebar{W}(K)|\simeq B|K|$, and so we can 
take $\overline{W}(\caln\AUT(\call))$ as a simplicial model for classifying 
space of the topological group $|\AUT(\call)|$.

A \emph{twisting function} $t$ from a simplicial set $X$ to a simplicial 
group $U$ is a collection of maps  $t_n\colon X_n\to U_{n-1}$ satisfying 
certain identities formulated in \cite[Definition 3.14]{Curtis}. 
When $X=\caln\calb(G)$ and $U=\caln\AUT(\call)$, a twisting function 
$\phi\:\caln\calb(G)\lra \caln\AUT(\call )$ is a collection of maps
	$$\phi_n \: \caln_n\calb(G) \Right2{} \caln_{n-1}\AUT(\call)
	\qquad\qquad \textup{(all $n\ge1$)} $$
satisfying the relations 
	\begin{equation}\label{twisting.function}
	\begin{split}
	\phi_{n-1} (d_i\7{g}) &=
	d_{i-1}\phi_n (\7{g}) 
	\qquad\qquad\qquad \text{for $2\le i\le n$,} \\
	\phi_{n-1} (d_1\7{g}) &=
	d_0\phi_n (\7{g}) \cdot
	\phi_{n-1} (d_0\7{g}) \,,\\
	\phi_{n+1} (s_{i}\7{g})&=s_{i-1}\phi_{n}(\7{g})  
	\qquad\qquad\qquad\text{for  $i\geq1$,} \\
	\phi_{n+1} (s_0\7{g})  &  =\mathbf{1}
\end{split}
\end{equation}
for all $n\geq1$ and all $\7g\in \caln_n\calb(G)$.  Here, 
$\mathbf{1}$ denotes the identity element in $\caln_{n}\AUT(\call)$. 

To a twisting function $\phi=\{\phi_n\}_{n\ge1}$ as above, 
one associates the simplicial map 
	\begin{equation}\label{assoc.map}\varphi\colon  \caln\calb(G) \Right2{}\overline{W}(\caln\AUT(\call ))\end{equation}
where for each $\7{g}=[g_1|\ldots|g_n]\in\caln_n\calb(G)$, 
	\begin{align*} 
	\varphi(\7g) &= \bigl(\phi_n(\7g), \phi_{n-1}(d_0\7g), \ldots, 
	\phi_1(d_0^{n-1}\7g)\bigr) \\
	&= \bigl(\phi_n([g_1|\ldots|g_n]), 
	\phi_{n-1}([g_2|\ldots|g_n]), \dots, 
	\phi_1([g_n])\bigr). 
	\end{align*}
It is not hard to see that the simplicial equalities 
(\ref{twisting.function}) are designed so that $\varphi$ is a 
simplicial map. For example, $\varphi$ commutes with the face maps by 
\eqref{twisting.function} together with the face relations 
	\[ d_0^jd_i = \begin{cases} d_{i-j}d_0^j & \textup{if $j\le i$} \\ 
	d_0^{j+1} & \textup{if $j\ge i$\,,} 
	\end{cases} \]
and formal manipulations.

The following lemma will be needed.

\begin{Lem} \label{l:chi}
Fix a linking triple $\SFL$, and an extension pair 
$U=(\wh{\Gamma},\tau)$ for $\call$ and $G$.  Let $\rho\:\wh\Gamma\Right2{}G$ be the 
surjection with kernel $\Gamma\defeq\Aut_\call(S)$, and let $t_U\colon G\Right2{}\wh\Gamma$ be a regular section of $\rho$ (i.e., $t_U$ is a right inverse for $\rho$ as a map of sets, and $t_U(1) = 1$).  
Define maps of sets
	\begin{align} 
	t &\: G\Right3{}\Aut\typ^I(\call) & &\textup{by setting} & 
	t&=\tau\circ{}t_U \tag*{} \\
	v &\: G\times G \Right3{} \Inn(\call) & &\textup{such that} &
	t(g)t(h)&=v(g,h)t(gh) ~\forall\,g,h\in{}G \label{e:tvchi} \\
	\chi &\: G\times G \Right2{} \Aut_\call(S) & &\textup{such that} & 
	t_U(g)t_U(h)&=\chi(g,h)t_U(gh) ~\forall\,g,h\in{}G \tag*{} \,.
	\end{align} 
Then the following hold.
\begin{enuma}
\item For each $g,h\in G$,  $v(g,h)=c_{\chi(g,h)}\in\Inn(\call)$.
\item For each $g\in G$, $\chi(1,g)=\chi(g,1)=1$.
\item For each $g,h,k\in G$, $\chi(g,h) \chi(gh,k) = t(g)(\chi(h,k))  
\chi(g,hk)$.
\end{enuma}
\end{Lem}

\begin{proof} Part (a) follows by the commutativity of triangle \pplus\ in 
Definition \ref{d:U(L,G)}(a), and (b) is immediate from the definition 
of $\chi$. It remains to prove (c).

By definition, 
	\[
	(t_U(g)t_U(h))t_U(k) = \chi(g,h)\chi(gh,k)\cdot t_U(ghk),\]
	and
	\[t_U(g)(t_U(h)t_U(k))= \9{t_U(g)}\chi(h,k) \cdot \chi(g,hk) 
	\cdot t_U(ghk)
	= t(g)(\chi(h,k)) \cdot \chi(g,hk) \cdot t_U(ghk),\]
where the last equality holds because of the commutativity of triangle 
\ppluss\ in Definition \ref{d:U(L,G)}(a). The claim now follows by 
the associativity of multiplication in $\Gamma$.
\end{proof}

We are now ready to show that for given $\call$ and $G$, there is a 
bijective correspondence between twisting functions from $\caln\calb(G)$ to 
$\caln\AUT(\call)$ and extension pairs for $\call$ and $G$.  

\begin{Prop}\label{constwistingfunc}
Fix a linking triple $\SFL$ and a finite group $G$. 
\begin{enuma} 

\item Let $U = (\widehat{\Gamma}, \tau)$ be  an extension pair for $\call$ 
and $G$, and let $t_U\colon G\to \widehat{\Gamma}$ be a regular section. 
Then there is a unique twisting function 
$\phi=\phi_U\: \caln\calb(G)\to\caln\AUT(\call)$ such that  
\smallskip

\begin{enumi}
\item for all $g\in G$, $\tau\circ t_U(g) = \phi_1([g])$; and
\item for all $g,h\in G$, $t_U(g)t_U(h) = \phi_2([g|h])t_U(gh)$.
\end{enumi} 

\item Let $\phi=\{\phi_n\}$ be any twisting function from $\caln\calb(G)$ 
to $\caln\AUT(\call)$.  There is an extension pair $U  = 
(\widehat{\Gamma}, \tau)$, and a regular section $t_U\colon 
G\to\widehat{\Gamma}$, such that $\phi_U=\phi$.

\end{enuma}
\end{Prop}

\newcommand{\g}[2]{g_{#1,#2}}

\begin{proof} 
\noindent\textbf{(a) } Let $U=(\widehat{\Gamma}, \tau)$ be an extension 
pair for $\call$ and $G$, and fix a regular section $t_U\colon G\to 
\widehat{\Gamma}$ of the natural projection. Since a simplex in the nerve 
of a category is determined by its 1-faces, there is at most one twisting 
function $\phi$ which satisfies (i) and (ii). We will prove that such a 
twisting function exists.

Define $t\colon G\to \Aut\typ^I(\call)$,  $v\colon G\times G\to\Inn(\call)$ 
and $\chi\colon G\times G\to\Aut_\call(S)$, as in \eqref{e:tvchi}.
For $\7g = [g_1|\cdots|g_n]$ and $i\le{}j$, write 
$\g{i}{j}=g_ig_{i+1}\cdots{}g_j$. Define $\phi$ by setting 
	\beqq\label{formulaphin} \begin{split}
	\phi_n([g_1|g_2|\dots & |g_n]) 
	= \Bigl( t(g_1) \Left5{\chi(g_1,g_2)} 
	\2{=v(g_1,g_2)^{-1}t(g_1)}{t(\g12)t(g_2)^{-1}} 
	\xfrom{\chi(g_1,g_2)^{-1} \chi(g_1,\g23)} \cdots \\ 
	& \Left4{} 
	\2{=v(g_1,\g2{n-1})^{-1}t(g_1)}{t(\g1{n-1})t(\g2{n-1})^{-1}}
	\xfrom{\chi(g_1,\g2{n-1})^{-1} \chi(g_1,\g2n)} 
	\2{=v(g_1,\g2n)^{-1}t(g_1)}{t(\g1n)t(\g2n)^{-1}} \Bigr) \,. 
	\end{split} \eeqq
The equalities between the 
objects hold by \eqref{e:tvchi} (the definition of $v$).

The first and third relations in (\ref{twisting.function}) clearly 
hold for $\phi$.  Also,  $\phi_n([g_1|\cdots|g_n])$ is the identity 
sequence if $g_1=1$, which implies the fourth identity in 
(\ref{twisting.function}). Thus we need only check the second 
relation, namely that
	\beqq d_0\phi_n([g_1|\cdots|g_n]) \,\cdot\, \phi_{n-1}([g_2|\cdots|g_n]) 
	= \phi_{n-1}([g_1g_2|g_3|\cdots|g_n]) \label{e:tp3} \eeqq
for each $g_1,\dots,g_n\in{}G$.  Each side of \eqref{e:tp3} lies in 
$\caln_{n-2}\AUT(\call)$; i.e., a sequence of $(n-1)$ objects and 
$(n-2)$ morphisms. The objects on each side of \eqref{e:tp3} are the same, 
since
	\[ t(\g1m)t(\g2m)^{-1} \circ t(\g2m)t(\g3m)^{-1} 
	= t(\g1m)t(\g3m)^{-1}  \]
for all $m\ge2$ (where we now set $\g32=1$ and $\g22=g_2$).  To show that 
the morphisms in the sequences coincide, we must show that 
	\[ \bigl(\chi(g_1,\g2m)^{-1}\chi(g_1,\g2{m+1})\bigl)\cdot
	\bigl(\chi(g_2,\g3m)^{-1} \chi(g_2,\g3{m+1})\bigr) 
	= \chi(\g12,\g3m)\chi(\g12,\g3{m+1}) \]
for all $m\ge 2$ (where again $\g32=1$). By \eqref{composition}, 
the left hand side of this equation takes the form
	\beqq\label{e:6.2a}\begin{split}
	(\chi(g_1,\g2m)^{-1} &\circ \chi(g_1,\g2{m+1})) \circ 
	(t(\g1{m+1})t(\g2{m+1}^{-1}))(\chi(g_2,\g3m)^{-1} \circ
	\chi(g_2,\g3{m+1})) \\
	=&\chi(g_1,\g2m)^{-1} \chi(g_1,\g2{m+1})\circ  
	c_{\chi(g_1,\g2{m+1})}^{-1}t(g_1)
	(\chi(g_2,\g3m)^{-1}\chi(g_2,\g3{m+1}))\\
	=& \chi(g_1,\g2m)^{-1} t(g_1)(\chi(g_2,\g3m))^{-1} \circ 
	t(g_1)(\chi(g_2,\g3{m+1})) \chi(g_1,\g2{m+1}) \\
	=& \bigl(\chi(g_1,g_2) \chi(g_1g_2,\g3m)\bigr)^{-1} \chi(g_1,g_2) 
	\chi(\g12,\g3{m+1}) \\
	=& \chi(\g12,\g3m)^{-1} \chi(\g12,\g3{m+1}) \,.
	\end{split}\eeqq
Here, the first equality in \eqref{e:6.2a} follows from 
\eqref{e:tvchi} and the definition of $v$, and  the third follows 
upon applying Lemma \ref{l:chi}(c) twice. This finishes the proof that 
$\phi$ is a twisting function. 

By construction, $\phi_1([g]) = t(g)  = \tau\circ t_U(g)$, and so $\phi$ 
satisfies (i). Upon setting $n=2$ in (\ref{formulaphin}), we may identify 
$\phi_2([g|h])$ with $\chi(g,h)$, and so (ii) holds by Lemma \ref{l:chi}.

\smallskip

\noindent\textbf{(b) } Fix a twisting function $\phi$, and define 
$t\:G\Right2{}\Aut\typ^I(\call)$ by $t(g)\defeq\phi_1([g])$.  Then 
$t(1)=1$.  

For each $g,h\in{}G$, the formulas in \eqref{twisting.function} for faces 
of $\phi_2([g|h])$ take the form
	\[ d_1\phi_2([g|h])=t(g) \qquad\textup{and}\qquad 
	d_0\phi_2([g|h])=t(gh)\cdot t(h)^{-1} \,. \]
So there is $\chi(g,h)\in\caln_1\AUT(\call)=\Aut_\call(S)$ such that 
	\[ \phi_2([g|h]) = \Bigl( t(g) \Left5{\chi(g,h)} t(gh)t(h)^{-1} 
	\Bigr)\,. \]
By \eqref{e:Mor(Aut)}, $t(g)=c_{\chi(g,h)}\circ t(gh) t(h)^{-1}$, and hence 
	\begin{equation}\label{e:tp1}
	t(g)t(h)=c_{\chi(g,h)}\cdot t(gh)\end{equation}
Also, by the degeneracy relations 
in \eqref{twisting.function}, for each $g\in G$, $\chi(g,1)=1=\chi(1,g)$.  

Similarly, for each $g,h,k\in{}G$, $\phi_3([g|h|k])$ has faces 
	\beqq \begin{split} 
	d_2\phi_3([g|h|k]) &= \phi_2([g|h])\,, \\
	d_1\phi_3([g|h|k]) &= \phi_2([g|hk])\,, \textup{ and} \\
	d_0\phi_3([g|h|k]) &= \phi_2([gh|k])\cdot\phi_2([h|k])^{-1} 
	\end{split} \label{e:di-phi3} \eeqq
by \eqref{twisting.function}.  Hence $\phi_3([g|h|k])$ takes the form
	\[ \phi_3([g|h|k]) = \Bigl( t(g) \Left5{\chi(g,h)} t(gh)t(h)^{-1} 
	\Left5{\chi_1} t(ghk)t(hk)^{-1} \Bigr)  \]
for some $\chi_1\in\Aut_\call(S)$.  
The second formula in \eqref{e:di-phi3} implies that 
$\chi(g,h)\circ\chi_1=\chi(g,hk)$, and hence $\chi_1=\chi(g,h)^{-1}\chi(g,hk)$.  
The third formula in \eqref{e:di-phi3} now implies that 
	\begin{multline*} 
	\Bigl( t(gh)t(h)^{-1} \xfrom{\chi(g,h)^{-1}\chi(g,hk)} 
	t(ghk)t(hk)^{-1} \Bigr) \cdot
	\Bigl( t(h) \Left5{\chi(h,k)} t(hk)t(k)^{-1} \Bigr) \\
	= \Bigl( t(gh) \Left5{\chi(gh,k)} t(ghk)t(k)^{-1} \Bigr) .
	\end{multline*}
Hence by definition of the monoidal structure on $\AUT(\call)$ 
(see \eqref{composition}), 
	\[ \chi(gh,k) = 
	\chi(g,h)^{-1}\chi(g,hk) \cdot (t(ghk)t(hk)^{-1})(\chi(h,k)).\]
By \eqref{e:tp1}, $c_{\chi(g,hk)}\circ t(ghk)t(hk)^{-1} = t(g)$. Hence 
	\begin{align} 
	\chi(g,h)\chi(gh,k) 
	&= \chi(g,hk) \cdot (c_{\chi(g,hk)}^{-1}\circ 
	t(g))(\chi(h,k)) \tag*{} \\
	&= t(g)(\chi(h,k))\chi(g,hk)\,. \label{e:tp2}
	\end{align}

Now define $\widehat{\Gamma}=\Aut_{\call}(S)\times{}G$, with group 
multiplication 
	\[ (a,g)\cdot(b,h) = (a\cdot t(g)(b)\cdot\chi(g,h),gh) \]
for each $g,h\in{}G$ and each $a,b\in\Aut_{\call}(S)$.  Define 
$\tau\:\widehat{\Gamma}\Right2{}\Aut\typ^I(\call )$ by setting 
$\tau(a,g)=c_a\circ{}t(g)$. 
For each $g,h,k\in{}G$ and each $a,b,c\in\Aut_{\call}(S)$,
	\begin{align*} 
	\bigl((a,g)\cdot(b,h)\bigr)\cdot(c,k) 
	&= \bigl(a\cdot t(g)(b)\cdot\chi(g,h)\cdot t(gh)(c)\cdot\chi(gh,k),
			ghk \bigr) \\
	&= \bigl(a\cdot t(g)(b)\cdot t(g)t(h)(c)\cdot\chi(g,h)\cdot\chi(gh,k),
			ghk \bigr) \\
	&= \bigl(a\cdot t(g)(b)\cdot t(g)t(h)(c)\cdot
			t(g)(\chi(h,k))\cdot\chi(g,hk), ghk \bigr) \\
	&= (a,g)\cdot\bigl(b\cdot t(h)(c)\cdot\chi(h,k), hk\bigr) \\
	&= (a,g)\cdot\bigl((b,h)\cdot(c,k)\bigr) \,,
	\end{align*}
where the second and third equalities follow from (\ref{e:tp1}) and 
(\ref{e:tp2}) respectively.  Thus multiplication in $\widehat{\Gamma}$ is 
associative, and $\widehat{\Gamma}$ is a group, with the obvious identity 
and inverses.  Similarly,
	\begin{align*} 
	\tau\bigl((a,g)\cdot(b,h)\bigr) 
	&= \tau(a\cdot t(g)(b)\cdot\chi(g,h),gh)
	= c_a\circ c_{t(g)(b)}\circ c_{\chi(g,h)}\circ t(gh)\\
	&= c_a\circ t(g)c_b t(g)^{-1}\circ c_{\chi(g,h)}t(gh) 
	= c_a\circ t(g)c_b t(g)^{-1}\circ t(g)t(h) \\
	&=\tau(a,g)\tau(b,h) \,, 
	\end{align*}
and so $\tau$ is a homomorphism.  

The  triangle \pplus\ in Definition \ref{d:U(L,G)})(a) commutes by 
definition of $\tau$. Triangle \ppluss\ commutes since 
	\[ (a,g)(b,1)(a,g)^{-1}=(at(g)(b)a^{-1},1) = (\tau(a,g)(b),1) \in\Gamma 
	\]
for each $(a,g)\in\widehat{\Gamma}$ and each $b\in\Aut_\call(S)$.  This 
shows that $U\defeq(\widehat{\Gamma},\tau)$ is an extension pair for 
$\call$ and $G$. 

Finally,  for $g\in G$, set $t_U(g)=(1,g)\in \widehat{\Gamma}$. Then $t_U$ 
is a regular section of the obvious projection 
$\rho\colon\widehat{\Gamma}\to G$.  By the definitions, $\tau(t_U(g)) 
= t(g) = \phi_1([g])$, while for $g,h\in G$, 
	\[t_U(g)t_U(h) = (1,g)\cdot(1,h) = (\chi(g,h), gh) 
	= (\chi(g,h),1)\cdot(1,gh). \] 
Upon identifying $(\chi(g,h),1)\in\widehat{\Gamma}$ with $\phi_2([g|h])$, 
we obtain (i) and (ii).  
\end{proof}

A twisting function $\phi\colon \caln\calb(G)\Right1{}\caln\AUT(\call)$ 
determines a map $BG\to B|\caln\AUT(\call)|$ (see (\ref{assoc.map})), and 
hence determines a fibre bundle over $BG$ with fibre $|\call|$ and 
structure group $|\AUT(\call)|$. This pullback bundle is the realization of 
a simplicial set $E(\phi)=|\call|\times_\phi BG$:  a twisted cartesian 
product which is described as follows (cf.\ \cite[(6.4)]{Curtis}).

\begin{Defi} \label{d:E(phi)}
For any twisting function $\phi\colon\caln\calb(G)\to\caln\AUT(\call)$, let 
$E(\phi)$ be the simplicial set with $n$-simplices $E(\phi)_n = 
\caln_n\call \times \caln_n\calb(G)$, and with face and degeneracy maps:
	\begin{align*}
	d_i(\7\xi,\7g) &= \begin{cases} 
	\bigl(\phi_n(\7g)^{-1}\cdot d_0\7\xi, d_0\7g\bigr) 
	& \textup{if $i=0$} \\
	(d_i\7\xi,d_i\7g) & \textup{if $i>0$} \end{cases} \\
	s_i(\7\xi, \7g) &=
                   (s_i\7\xi, s_i\7g)  
	\qquad\qquad\qquad \textup{for all $i$} 
	\end{align*}
for all  $\7\xi\in \caln_n\call$ and all $\7g\in \caln_n\calb(G)$. Let 
$p_\phi\colon E(\phi)\Right2{}\caln\calb(G)$ be the natural projection 
which maps a pair $(\7\xi,\7g)\in E(\phi)$  to $\7g$. 
\end{Defi}

By \cite[\S\,6.4]{Curtis}, in the situation of Definition 
\ref{d:E(phi)}, $p_\phi$ induces a fibre bundle $|p_\phi|\: |E(\phi)|\to 
BG$ with fibre $|\call|$ and structure group $|\AUT(\call)|$. 
We want to identify the nerve of the extension $\call_U$ of Definition 
\ref{d:U(L,G)} with $E(\phi)$ for the associated twisting function $\phi$. 
Before we can do this, one more technical lemma is required. 

\begin{Lem} \label{X.=N(X)}
Let $X$ be a simplicial set such that for each $n\ge2$, the map 
	\[ D_n =  \: X_n \Right8{((d_2)^{n-1},d_0)} 
	\bigl\{ (x_1,x_2)\in X_1 \times X_{n-1} 
	\,\big|\, d_0x_1 = (d_1)^{n-1}(x_2)\in X_0 \bigr\} \]
is a bijection.  Then there is a category $\calx$ with 
$\Ob(\calx)=X_0$ and $\Mor(\calx)=X_1$, where $f\in{}X_1$ is a morphism 
from $d_0f$ to $d_1f$ and $\Id_x=s_0x$ for $x\in{}X_0$, in which 
composition is defined as follows. If $f_1, f_2\in X_1$ are composable 
morphisms (i.e., if $d_0f_1=d_1f_2$), then their composite is 
defined by setting $f_1\circ{}f_2=d_1x$, where $x\in{}X_2$ is the unique 
element such that $D_2(x)=(f_1,f_2)$. In other words,  $f_1\circ f_2 = 
d_1\circ{}D_2^{-1}(f_1,f_2)$.  Furthermore, $\caln\calx\cong{}X$ as 
simplicial sets.
\end{Lem}

\begin{proof} For each $0\le i<n$, let  $e^n_i\:X_n\Right2{}X_1$ be the 
``edge map'' induced by the morphism $[1]\Right2{}[n]$ in $\Delta$ with 
image $\{i,i+1\}$.  Thus, for example, $e^2_0=d_2$, $e^2_1=d_0$, 
$e^n_0=(d_2)^{n-1}$ for $n>2$, etc.  

For each $n\ge2$, set 
	\[ \caln_n\calx = 
	\bigl\{(f_0,\dots,f_{n-1})\in (X_1)^n \,\big|\, d_0f_i=d_1f_{i+1} 
	\textup{ for all $0\le i\le n-2$} \bigr\}\,: \]
the set of $n$-tuples of composable morphisms.  Set 
	\[ E_n = (e^n_0,e^n_1,\ldots,e^n_{n-1}) \: X_n \Right5{} \caln_n\calx 
	\subseteq (X_1)^n. \]
It is easy to check that $E_n=(\Id\times D_2)\circ(\Id\times D_3)\circ\cdots
\circ D_n$, and so $E_n$ is a 
bijection since the $D_i$ are.  In particular, for each 
$(f_1,f_2,f_3)\in\caln_3\calx$, there is a unique $y\in{}X_3$ such that 
$E_3(y)=(f_1,f_2,f_3) = (d_2^2y, d_2d_0y, d_0^2y) = (d_2d_3y, d_0d_3y, d_0d_1y)$, and 
	\[ (f_1\circ f_2)\circ f_3 = \2{=d_2d_1y}{d_1d_3y} \circ 
	\2{=d_0d_1y}{f_3} = d_1d_1y = d_1d_2y 
	= \2{=d_2d_2y}{f_1} \circ \2{=d_0d_2y}{d_1d_0y} 
	= f_1\circ(f_2\circ f_3). \]
The first equality follows since $D_2(d_3y) = (d_2d_3y, d_0d_3y) = 
(f_1,f_2)$, and the fifth since $D_2(d_0y) = (d^2_0y, d_0d_0y) = 
(f_2,f_3)$, the second and fourth equalities hold by definition of $D_2$, 
and the third is a simplicial identity.   This shows that composition in 
$\calx$ is associative, and hence that $\calx$ is a category.  

Set  $E_0=\Id_{X_0}$, $E_1=\Id_{X_1}$, and 
$E=\{E_n\}\:X\Right2{}\caln\calx$.  By construction, $E$ commutes with face 
maps on $X_n$ for $n\le2$.  We claim that $d_i\circ{}E_n=E_{n-1}\circ{}d_i$ 
for $n\ge3$ and $0\le i\le n$.  This is clear from the definition of $E_n$, 
when $i=0$ or $i=n$, since in these two cases $d_0$ and $d_n$ on 
$\caln_n\calx$ are restrictions of the obvious projections 
$p_0,p_n\:(X_1)^n\Right2{}(X_1)^{n-1}$.  So assume $0<i<n$, and consider 
the map
	\[ T\: X_n \xto{(e^n_0,\dots,e^n_{i-2},t,e^n_{i+1},\dots,e^n_{n-1})} 
	(X_1)^{i-1} \times X_2 \times (X_1)^{n-i-1}\,, \]
where $t$ is induced by the morphism $[2]\Right2{}[n]$ with image 
$\{i-1,i,i+1\}$.  Then $E_n=(\Id,(d_2,d_0),\Id)\circ T$, while 
$E_{n-1}\circ{}d_i=(\Id,d_1,\Id)\circ T$.  In other words, for each 
$x\in{}X_n$ with $E_n(x)=(f_1,\dots,f_n)$, we have $E_{n-1}(d_ix)= 
(f_1,\dots,f_{i-1},f_i\circ{}f_{i+1},\dots,f_n)=d_i(E_n(x))$.

This proves that $E$ commutes with face maps, and it is easily seen (by the 
choice of identity morphisms) to commute with degeneracies. Since all of 
the $E_i$ are bijections, $E$ is an isomorphism of simplicial sets. 
\end{proof}

For any extension pair $U=(\widehat{\Gamma},\tau)$ for $\call$ and $G$, we 
let $\pr\:\call_U\Right2{}\calb(G)$ be the functor defined by setting 
$\pr(\5{\varphi,\gamma})=\rho(\gamma)$.

\begin{Prop} \label{LU=Ephi}
Fix a linking triple $\SFL$, a finite group $G$, and an extension 
pair $U=(\widehat{\Gamma},\tau)$ for $\call$ and $G$ with a regular 
section $t_U\colon G\to \widehat{\Gamma}$.  Let $\phi=\phi_U$ be the 
associated twisting function, as in Proposition \ref{constwistingfunc}(a).  
Then $\caln\call_U\cong E(\phi)$, via a simplicial isomorphism which 
commutes with the projections to $\caln\calb(G)$.  In particular, 
$|\call_U|\Right3{|\pr_U|}BG$ is a fibre bundle over $BG$ with fibre 
$|\call|$ and structure group $|\AUT(\call)|$.
\end{Prop}

\begin{proof} We first claim, for each $n\ge2$, that the map 
	\[ D_n = (d_2^{n-1},d_0) \: E(\phi)_n \Right5{} 
	\bigl\{ (\7\eta_1,\7\eta_2)\in E(\phi)_1 \times E(\phi)_{n-1} 
	\,\big|\, d_0\7\eta_1 = d_1^{n-1}\7\eta_2\in E(\phi)_0 \bigr\} \]
is a bijection.  For $\7\xi=(P_0\xfrom{\varphi_1} P_1\from\cdots\from P_n)$ 
and $\7g=[g_1|\cdots|g_n]$, 
	\[ D_n(\7\xi,\7g) = \Bigl( \bigl(P_0\xfrom{\varphi_1} P_1,[g_1]\bigr) 
	\,,\, \bigl(\phi_n(\7g)^{-1}\cdot d_0\7\xi,[g_2|\cdots|g_n]\bigr) 
	\Bigr) \,, \]
and 
	\begin{multline*} 
	d_1^{n-1}\bigl(\phi_n(\7g)^{-1}\cdot d_0\7\xi,[g_2|\cdots|g_n]\bigr) = 
	\bigl(\phi_1(d_2^{n-1}\7g)^{-1}\cdot d_1^{n-1}d_0\7\xi, *) \\ 
	=\bigl(\phi_1([g_1])^{-1}P_1, *)= d_0\bigl(P_0\xfrom{\varphi_1} 
	P_1,[g_1]\bigr).
	\end{multline*}
For fixed $\7g$, $D_n$ restricts to a function from $\caln_n\call$ to the 
set of pairs $(\7\xi_1,\7\xi_2)\in \caln_1\call\times\caln_{n-1}\call$ such 
that 
	\begin{equation}\label{e:bij}
	\phi_1([g_1])^{-1}(d_0\7\xi_1)=d_1^{n-1}\7\xi_2.
	\end{equation}
This function is in fact a bijection, since by \eqref{e:bij}, the last term 
of $\7\xi_2$ is the first term of $\phi_1([g_1])^{-1}(d_0\7\xi_1)$, and 
since $\phi_1([g_1])$ is invertible. Since $\7g = [g_1|\cdots|g_n]$ is 
determined by $[g_1]$ and $d_0\7g = [g_2|\cdots|g_n]$, it follows that 
$D_n$ is a bijection.  

By Lemma \ref{X.=N(X)}, $E(\phi)\cong\caln\cale(\phi)$, where 
$\cale(\phi)$ is the category with 
	\[ \Ob(\cale(\phi))=E(\phi)_0=\Ob(\call)\times\{*\} 
	\quad\textup{and}\quad
	\Mor(\cale(\phi))=E(\phi)_1=\Mor(\call)\times{}G\,, \]
with source and target defined by $d_0$ and $d_1$, respectively, with 
$\Id_{(P,*)}=s_0(P,*)=(\Id_P,1)$, and with composition defined by 
$d_1\circ{}D_2^{-1}$. Let $t$  and $\chi$ be the functions associated to 
$U$ and $t_U$ via \eqref{e:tvchi} (we do not use $v$ here).  For each 
$\varphi\in\Mor_\call(P,Q)$ and $g\in{}G$,
	\begin{align*} 
	d_0\bigl(Q\Left2{\varphi}P,[g]\bigr) &= 
	\bigl(\phi_1([g])^{-1}(P),*\bigr)=\bigl(t(g)^{-1}(P),*\bigr) \,, \\
	d_1\bigl(Q\Left2{\varphi}P,[g]\bigr) &= (Q,*) \,, 
	\end{align*}
so $(\varphi,g)\in\Mor_{\cale(\phi)}\bigl((t(g)^{-1}(P),*),(Q,*)\bigr)$.  

We next describe composition in $\cale(\phi)$.  For each 
	\[ \7\eta = \bigl( R \Left2{\varphi} Q \Left2{\lambda} P \,,\,
	[g|h] \bigr) \in E(\phi)_2, \]
$d_1\7\eta=(\varphi\circ\lambda,[gh])$ and $d_2\7\eta=(\varphi,[g])$, while 
	\begin{align*} 
	d_0\7\eta &= \bigl( \phi_2([g|h])^{-1}(Q\Left2{\lambda}P),[g] 
	\bigr) \\
	&=\Bigl( \bigl(t(g)\Left4{\chi(g,h)}t(gh)t(h)^{-1}\bigr)^{-1} 
	\cdot \bigl( Q\Left2{\lambda}P \bigr) , [g] \Bigr) \\
	&= \bigl( t(g)^{-1}(Q) \Left2{\psi} t(h)t(gh)^{-1}(P), [g] \bigr)
	\end{align*}
where by formula \eqref{e:action} for the action of $\Mor(\AUT(\call))$ on 
$\Mor(\call)$, 
	\[ \lambda = \bigl( t(gh)t(h)^{-1}\Right4{\chi(g,h)}t(g) \bigr) 
	\cdot \psi = t(g)(\psi)\circ\chi(g,h)\,. \]
In other words, composition in $\cale(\phi)$ satisfies 
	\beqq (\varphi,g) \circ (\psi,h) = (\varphi\circ\lambda,gh) = 
	\bigl(\varphi\circ t(g)(\psi)\circ\chi(g,h),gh\bigr)\,. 
	\label{e:compos} \eeqq

It remains to show that $\cale(\phi)\cong\call_U$.  
By construction, each morphism in $\call_U$ has the form 
$\5{\varphi,t_U(g)}$ for some unique $\varphi$ and $g$.  Define 
$\omega\:\call_U\Right2{}\cale(\phi)$ by setting $\omega(P)=(P,*)$ and 
$\omega(\5{\varphi,t_U(g)})=(\varphi,g)$.  By \eqref{e:compos}, this 
preserves composition, and hence is an isomorphism of categories. 
\end{proof}

\bigskip


\section{The proof of Theorem A}

Before proving Theorem \ref{ThA}, we need one more result, which 
allows us to compare fibrations with fibre $|\call|\pcom$ and fibre 
bundles with fibre $|\call|$.  

\begin{Prop} \label{Phi-isom}
For each $p$-local compact group $\SFL $ and each finite group $G$, 
there is a bijection $\Phi$ from 
the set of equivalence classes of fibre 
bundles over $BG$ with fibre $|\call|$ and structure group $|\AUT(\call)|$
to the set of equivalence classes of 
fibrations over $BG$ with fibre homotopy equivalent to $|\call |\pcom$:  a 
bijection which sends 
the class of a fibre bundle to 
the equivalence class of its fibrewise $p$-completion.
\end{Prop}

\begin{proof} For any space $X$, let $\Aut(X)$ denote the space of its self homotopy equivalences.  Let 
	\[ \Omega\colon |\AUT(\call)| \Right2{} \Aut(|\call|)
	\qquad\textup{and}\qquad
	\widehat{\Omega}\colon |\AUT(\call)| \Right2{} \Aut(|\call|\pcom) 
	\]
be the homomorphisms induced by the evaluation functor 
$\AUT(\call)\times\call \Right2{}\call$ and by $p$-completion.

By \cite[Theorem 5.6]{BGM}, equivalence classes of fibrations over $BG$ 
with fibre $|\call|\pcom$ are in one-to-one correspondence with the set 
$\bigl[BG,B\Aut(|\call|\pcom)\bigr]_*$ of homotopy classes of pointed maps.  
Also, equivalence classes of $|\AUT(\call)|$-bundles over $BG$ with fibre 
$|\call|$ are in one-to-one correspondence with the set of homotopy classes 
of pointed maps $\bigl[BG,B|\AUT(\call )|\bigr]_*$.  

If a map $f\: BG\Right2{}B|\AUT(\call)|$ classifies an 
$|\AUT(\call)|$-bundle $\xi_f$ with fibre $|\call|$, then 
$B\widehat{\Omega}\circ f$ classifies the 
fibrewise $p$-completion $\widehat{\xi}_f$ of $\xi_f$.  So we must show 
that the map 
	\[ \Phi\: \bigl[BG,B|\AUT(\call )|\bigr]_* 
	\Right7{B\widehat{\Omega}\circ-} 
	\bigl[BG, B\Aut(|\call |\pcom)\bigr]_* \]
is a bijection. 

By \cite[Theorem 7.1]{BLO3}, $\widehat \Omega$ induces a homotopy 
equivalence after (componentwise) $p$-completion.  Hence 
$B\Aut(|\call|\pcom)$ is the fibrewise $p$-completion of $B\AUT(\call)$ 
over $B\Out(|\call|\pcom)\cong B\Out\typ(\call)$, and there is a map of 
fibration sequences 
	$$ \xymatrix{
	K(Z,2) \ar[r]\ar[d]^{\kappa_p} & 
	B|\AUT(\call)|\ar[r]\ar[d]^{B\widehat\Omega} & 
	B\Out\typ(\call ) \ar[d]^{\simeq} \\
	K(Z,2)\pcom \ar[r]  &  B\Aut(|\call |\pcom) \ar[r] & 
	B\Out(|\call |\pcom ) \rlap{\,.}
	} $$
Here, $Z=\Aut_{\AUT(\call)}(\Id)$, and by \eqref{e:Mor(Aut)} in 
Section 6, it can be identified as a subgroup of $Z(S)$.  (In fact, 
$Z=Z(\calf)$:  the center of the fusion system $\calf$.)  Hence it is an 
abelian discrete $p$-toral group of the form $Z\cong (\Z/p^\infty)^r \times 
A$ for some finite abelian $p$-group $A$.


By \cite[Proposition VI.5.1]{BK}, $K(\Z/p^\infty,2)\pcom\simeq K(\Z_p,3)$, 
$K(\Q_p,2)\pcom\simeq*$, and $K(\Z_p,2)$ is $p$-complete.  Hence by the 
fibre completion lemma \cite[Lemma II.4.8]{BK}, there is a map of fibration 
sequences
	$$ \xymatrix{
	K(\Z_p,2) \ar[r]\ar[d]^{\cong} & K(\Q_p,2) \ar[r]\ar[d] 
	& K(\Z/p^\infty,2) \ar[d]^{\kappa} \\
	K(\Z_p,2) \ar[r] & \* \ar[r] & K(\Z_p,3) \rlap{\,,}
	} $$
and so $\hofibre(\kappa)\simeq K(\Q_p,2)$.  Thus 
$K(Z,2)\pcom \simeq K((\Z_p)^r,3)\times K(A,2)$, and the homotopy fibre 
of $B\widehat{\Omega}$ is equivalent to $\textup{hofibre}(\kappa_p)\simeq 
K((\Q_p)^r,2)$.

Since $G$ is finite, $H^i(G;(\mathbb{Q}_p)^r)=0$ for all $i\geq1$.  
Thus $\Phi$ is a bijection by obstruction theory.
\end{proof}

\begin{Rmk}
In particular,  we have shown that each class of fibrations over $BG$ with 
fibre $|\call|\pcom$ has a representative which is actually an 
$|\AUT(\call)|$-bundle. In other words, the structure group of a fibration 
can always be reduced to $|\AUT(\call )|$. 
\end{Rmk}

We are now ready to prove Theorem \ref{ThA}, in the following slightly more 
precise form.

\begin{Thm} \label{ThmA}
Assume $f\:X\Right2{}Y$ is a finite regular covering space with covering 
group $G$, where $X\simeq|\dec\call|\pcom$ is the classifying space of a 
$p$-local compact group $(\dec{S},\dec\calf,\dec\call)$.  Then $Y\pcom$ is 
the classifying space of a $p$-local compact group 
$(S,\calf,\call)$. Furthermore, there is a transporter system 
$\calt$ associated to $\calf$ and $\call$ such that 
$\dec\call\nsg\calt$, and such 
that $\Aut_{\calt}(\dec{S})\big/\Aut_{\dec\call}(\dec{S})\cong{}G$.
\end{Thm}

\begin{proof}  By Proposition \ref{Phi-isom}, there is a fibre bundle 
$|\dec\call|\Right2{}E\Right2{}BG$ with structure group $|\AUT(\dec\call)|$ 
whose fibrewise $p$-completion is equivalent up to homotopy to the 
fibration $X\Right2{}X\times_GEG\Right2{}BG$.  This bundle is classified by 
a map $\varphi$ from $BG$ to 
$B|\AUT(\dec\call)|=|\widebar{W}\AUT(\dec\call)|$.  Since 
$\widebar{W}\AUT(\dec\call)$ is fibrant (cf. \cite[Corollary V.6.8]{GJ}), 
we can assume that $\varphi$ is the realization of a simplicial map, and 
hence is determined by a twisting function $\phi$.  

By Proposition \ref{constwistingfunc}(b), $\phi=\phi_U$ for some 
extension pair $U=(\widehat{\Gamma},\tau)$ with regular section 
$t_U\colon G\to\widehat{\Gamma}$.  By Proposition \ref{LU=Ephi}, 
$|\dec\call_U|\cong E$ as bundles over $BG$. In particular, there is a mod 
$p$ equivalence from $|\dec\call_U|$ to $X\times_GEG\simeq{}Y$. 

By Theorem \ref{T:new-ext}, there is a saturated fusion system $\calf$ over 
a discrete $p$-toral group $S$, and a transporter system $\calt$ associated 
to $\calf$, such that $\dec\call\nsg\calt$, 
$\Aut_{\calt}(\dec{S})/\Aut_{\dec\call}(\dec{S})\cong{}G$, and such 
that $|\calt|$ contains $|\dec\call_U|$ as deformation retract.  Hence 
$|\calt|\pcom \simeq |\dec\call_U|\pcom \simeq Y\pcom$.  Let 
$\calt_1\subseteq\calt$ be the full subcategory whose objects are the 
objects of $\calt$ which are $\calf$-centric.  Then $\Ob(\calt_1)$ 
is invariant under $\calf$-conjugacy, closed under overgroups, and contains 
all subgroups which are $\calf$-centric and $\calf$-radical. Hence 
$|\calt_1|\pcom\simeq|\calt|\pcom$ by Proposition \ref{dropmore}(a).  By 
Proposition \ref{centric-E(P)}, there is a linking system $\call_1$ 
associated to $\calt_1$ and to $\calf$, with $\Ob(\call_1)=\Ob(\calt_1)$, 
and such that $|\call_1|\pcom\simeq|\calt_1|\pcom$.  By Proposition 
\ref{existsLc}, there is a centric linking system $\call$ associated to 
$\calf$ which contains $\call_1$ as a full subcategory, and 
$|\call|\pcom\simeq|\call_1|\pcom$ by Corollary \ref{L0<L}.  So 
$|\call|\pcom\simeq{}Y\pcom$, and this finishes the 
proof of the theorem.  
\end{proof}


\appendix

\section{Transporter systems over discrete $p$-toral groups}

\newcommand{\oc}[2]{#1\,{\downarrow}\,#2} 
\newcommand{\iii}{i}
\newcommand{\6}[1]{#1{}^{\vee}}
\renewcommand{\8}{{}^\bullet}

Transporter systems are a generalization of linking systems, which were 
first defined (over finite $p$-groups) in \cite{OV1}.  We need them here in 
order to state our main theorem on extensions of fusion and linking systems 
(Theorem \ref{T:new-ext}) in sufficient generality.  In the first half of 
this section, we define and prove the basic properties of transporter 
systems over discrete $p$-toral groups; especially those properties needed 
to prove and apply Theorem \ref{T:new-ext}.  

Afterwards, we give some conditions under which an inclusion of transporter 
or linking systems (one a full subcategory of the other) induces an 
equivalence or mod $p$ equivalence of geometric realizations.  We also 
prove that every linking system all of whose objects are centric can be 
embedded in a centric linking system.  When it is a question of adding only 
finitely many conjugacy classes, these results can be proven using 
arguments similar to those already used in \cite{BCGLO1} and 
\cite{OV1} for linking and transporter systems over finite $p$-groups.  
What is new here (and makes the proofs harder) is the necessity of handling 
infinitely many classes at a time.

Let $G$ be a group, and let $\calh$ be a family of subgroups of $G$ which 
is invariant under $G$-conjugacy and overgroups. The transporter system of 
$G$ with respect to $\calh$ is the category $\calt_\calh(G)$ with object 
set $\calh$, and which has morphism sets 
	\[ \Mor_{\calt_\calh(G)}(P,Q) = T_G(P,Q) \defeq 
	\{x\in G\;|\; xPx^{-1} \le Q\} \]
(the \emph{transporter set}) for each pair of subgroups $P,Q\in\calh$.

\begin{Defi} \label{tport}
Let $\calf$ be a fusion system over a discrete $p$-toral group $S$. A \emph{transporter 
system} associated to a fusion system $\calf$ is a nonempty category 
$\calt$ such that $\Ob(\calt)\subseteq\Ob(\calf)$  is closed under 
$\calf$-conjugacy and overgroups, together with a pair of functors
    \[ \calt_{\Ob(\calt)}(S) \Right5{\gee} \calt
    \Right5{\rho} \calf\,, \]
satisfying the following conditions:
\begin{itemize}
\item[(A1) ]
The functor $\gee$ is the identity on objects and the functor $\rho$
is the inclusion on objects.

\item[(A2) ]
For each $P,Q\in\Ob(\calt)$, the kernel
    \[ E(P)\defeq\Ker\bigl[\rho_{P}\:
    \Aut_{\calt}(P)\longrightarrow\autf(P)\bigr] \]
acts freely on $\Mor_{\calt}(P,Q)$ by right composition, and
$\rho_{P,Q}\colon\Mor_\calt(P,Q)\to\Hom_\calf(P,Q)$ is the orbit map
for this action. Also, $E(Q)$ acts freely on $\Mor_{\calt}(P,Q)$ by
left composition.

\item[(B) ]
For each $P,Q\in\Ob(\calt)$,
$\gee_{P,Q}\:N_S(P,Q)\Right2{}\Mor_\calt(P,Q)$ is injective, and the
composite $\rho_{P,Q}\circ\gee_{P,Q}$ sends $g\in N_S(P,Q)$ to
$c_g\in\Hom_{\calf}(P,Q)$.

\item[(C) ]
For all $\varphi\in\Mor_{\calt}(P,Q)$ and all $g\in P$, the diagram
\[
\xymatrix@C+2ex{
\I02P \ar[r]^{\varphi} \ar[d]_{\gee_{P}(g)} & \I02Q
\ar[d]^{\gee_{Q}(\rho(\varphi)(g))} \\
\I02P \ar[r]^{\varphi} & \I02Q
}
\]
commutes in $\calt$.

\item[(I) ]
Each $\calf$-conjugacy class of subgroups in $\Ob(\calt)$ contains a 
subgroup $P$ such that $\gee_P(N_S(P))\in\sylp{\Aut_{\calt}(P)}$; i.e., 
such that $[\Aut_\calt(P):\gee_P(N_S(P))]$ is finite and prime to $p$.

\item[(II) ] Let $\varphi\in\Iso_{\calt}(P,Q)$, $P\normal\widebar{P}\le S$, 
and $Q\normal\widebar{Q}\le S$ be such that 
$\varphi\circ\gee_{P}(\widebar{P})\circ\varphi^{-1}\le 
\gee_{Q}(\widebar{Q})$. Then there is some 
$\widebar\varphi\in\Mor_{\calt}(\widebar{P},\widebar{Q})$ such that 
$\widebar\varphi\circ\gee_{P,\widebar P}(1) = 
\gee_{Q,\widebar{Q}}(1)\circ\varphi$.

\item[(III) ] Assume $P_1\le P_2\le P_3\le\cdots$ in $\Ob(\calt)$ and 
$\psi_i\in\Mor_\calt(P_i,S)$ are such that for all $i\ge1$, 
$\psi_i=\psi_{i+1}\circ\gee_{P_i,P_{i+1}}(1)$.  Set $P=\bigcup_{i=1}^\infty 
P_i$.  Then there is $\psi\in\Mor_\calt(P,S)$ such that 
$\psi\circ\gee_{P_i,P}(1)=\psi_i$ for each $i$.

\end{itemize}
\end{Defi}

When $P\le{}Q$ are both in $\Ob(\calt)$, we write 
$\iota_{P,Q}=\gee_{P,Q}(1)$, considered to be the ``inclusion'' of $P$ into 
$Q$. By axiom (B), $\rho$ sends $\iota_{P,Q}\in\Mor_\calt(P,Q)$ to 
$\incl_P^Q\in\homf(P,Q)$ (the inclusion in the usual sense).

\begin{Prop} \label{T-prop}
The following hold for any transporter system $\calt$ associated to a 
fusion system $\calf$ over a discrete $p$-toral group $S$.
\begin{enuma}

\item Fix morphisms $\varphi\in\homf(P,Q)$ and 
$\psi\in\homf(Q,R)$, where $P,Q,R\in\Ob(\calt)$.  Then for any pair of 
liftings $\til\psi\in\rho^{-1}_{Q,R}(\psi)$ and $\til{\psi\varphi}\in 
\rho^{-1}_{P,R}(\psi\varphi)$, there is a unique lifting $\til\varphi\in 
\rho^{-1}_{P,Q}(\varphi)$ such that $\til\psi\circ\til\varphi= 
\til{\psi\varphi}$.

\item For every morphism $\psi\in\Mor_\calt(P,Q)$, and every 
$P_*,Q_*\in\Ob(\calt)$ such that $P_*\le{}P$, $Q_*\le{}Q$, and 
$\rho(\psi)(P_*)\le{}Q_*$, there is a unique morphism 
$\psi_*\in\Mor_\calt(P_*,Q_*)$ such that 
$\psi\circ\iota_{P,P_*}=\iota_{Q,Q_*}\circ\psi_*$.

\item For each $\psi\in\Mor(\calt)$, $\psi$ is an isomorphism in $\calt$ if 
and only if $\rho(\psi)$ is an isomorphism in $\calf$.

\item All morphisms in $\calt$ are monomorphisms and epimorphisms in the
categorical sense.

\end{enuma}
\end{Prop}

\begin{proof} Point (a) follows from \cite[Lemma A.7(a)]{OV1}, and (b) is a 
special case of (a).  All morphisms in $\calt$ are monomorphisms by 
\cite[Lemma A.7(b)]{OV1} and since morphisms in $\calf$ are 
monomorphisms.

If $\psi\in\Mor_\calt(P,Q)$ is such that $\rho(\psi)\in\isof(P,Q)$, then by 
(a), there are $\psi'\in\Mor_\calt(Q,P)$ and $\psi''\in\Mor_\calt(P,Q)$ 
such that $\psi\circ\psi'=\Id_Q$ and $\psi'\circ\psi''=\Id_P$.  Then 
$\psi=\psi\psi'\psi''=\psi''$ is an isomorphism in $\calt$ with $\psi'$ as 
inverse. This proves (c).

It remains to prove that all morphisms in $\calt$ are epimorphisms.  Fix 
$\psi\in\Mor_\calt(P,Q)$ and $\varphi_1,\varphi_2\in\Mor_\calt(Q,R)$ such 
that $\varphi_1\circ\psi=\varphi_2\circ\psi$; we must show that 
$\varphi_1=\varphi_2$.  Since $\psi$ is the composite of an isomorphism 
followed by an inclusion by (b) and (c), it suffices to prove this when 
$P\le{}Q$ and $\psi=\iota_{P,Q}$ is the inclusion.  

Assume we can show that $\rho(\varphi_1) = \rho(\varphi_2)$. By axiom 
(A2),  $\varphi_2=\varphi_1\circ\alpha$ for some $\alpha\in{}E(Q)$.  Hence
	\[ \varphi_1\circ\alpha\circ\iota_{P,Q}=
	\varphi_2\circ\iota_{P,Q} = 
	\varphi_1\circ\iota_{P,Q}, \]
so $\alpha\circ\iota_{P,Q}=\iota_{P,Q}$ by (a), and $\alpha=\Id_Q$ since 
$E(Q)$ acts freely on $\Mor_\calt(P,Q)$ (axiom (A2)).  Thus 
$\varphi_1=\varphi_2$.

To complete the proof, we need to show that $\rho(\varphi_1) = 
\rho(\varphi_2)$. Assume otherwise, and let $\widebar{P}<Q$ be the subgroup 
of all $g\in{}Q$ such that $\rho(\varphi_1)(g)=\rho(\varphi_2)(g)$.  Then 
$P<N_Q(\widebar{P})$ by Lemma \ref{N_Q(P)>P}, and 
$\rho(\varphi_1)|_{N_Q(\widebar{P})}\ne\rho(\varphi_2)|_{N_Q(\widebar{P})}$ 
by definition of $\widebar{P}$.  So upon replacing $P$ by $\widebar{P}$ and 
$Q$ by $N_Q(\widebar{P})$, we can arrange that $P\nsg Q$.  

Set $P'=\rho(\varphi_1)(P)=\rho(\varphi_2)(P)$.  By (b), 
$\varphi_1\circ\iota_{P,Q}=\varphi_2\circ\iota_{P,Q}$ has a unique 
restriction $\beta=\varphi_1|_{P,P'}=\varphi_2|_{P,P'} 
\in\Iso_\calt(P,P')$.  Fix $x\in{}Q$, set $y_i= \rho(\varphi_i)(x)$, and 
consider the following two squares for each $i=1,2$:
	\[ \vcenter{
	\xymatrix@C+3ex@!{
	Q \ar[r]^{\varphi_i} \ar[d]_{\gee_{Q}(x)} & 
	R \ar[d]^{\gee_{R}(y_i)} &&
	P \ar[r]^{\beta} \ar[d]_{\gee_{P}(x)} 
	& P' \ar[d]^{\gee_{P'}(y_i)} \\
	Q \ar[r]^{\varphi_i} & R &&
	P \ar[r]^{\beta} & P' \\
	}} \]
The first square commutes by axiom (C), and the second square is defined to 
be a restriction of the first.  Note that $\gee_P(x)$ sends $P$ to itself 
because $P\nsg Q'$. Note also that $\gee_{P}(x)$ is the restriction of 
$\gee_{Q'}(x)$ since $\gee$ is a functor (and since 
$\iota_{U,V}=\gee_{U,V}(1)$ for all $U\le{}V$ in $\Ob(\calt)$).   Hence the 
second square commutes by the uniqueness of restriction morphisms (point 
(b)).  Thus $\gee_{P'}(y_1)= \gee_{P'}(y_2)= 
\beta\circ\gee_{P}(x)\circ\beta^{-1}$.  Since $\gee_{P'}$ is injective 
(axiom (B)), this shows that $y_1=y_2$.  Since this holds for all 
$x\in{}Q$, $\rho(\varphi_1)=\rho(\varphi_2)$, contradicting our assumption 
that they are distinct.
\end{proof}

As usual, in the situation of Proposition \ref{T-prop}(b), we write 
$\psi_*=\psi|_{P_*,Q_*}\in\Mor_\calt(P_*,Q_*)$, regarded as the restriction 
of $\psi$.  

\begin{Prop} \label{T=>Fsat}
For any fusion system $\calf$ over a discrete $p$-toral group $S$, and any 
transporter system $\calt$ associated to $\calf$,  $\calf$ is 
$\Ob(\calt)$-saturated.
\end{Prop}

\begin{proof} Set $\calh=\Ob(\calt)$ for short.  
By axiom (I) for $\calt$, each subgroup in $\calh$ is $\calf$-conjugate to 
some $P$ such that $\gee_P(N_S(P))\in\sylp{\Aut_\calt(P)}$.  Hence  
	\beqq \Aut_S(P)\in\sylp{\autf(P)} \qquad\textup{and}\qquad 
	\gee_P(C_S(P))\in\sylp{E(P)} \label{e:A.2a} \eeqq
where as usual, $E(P)=\Ker\bigl[\Aut_\calt(P)\Right2{\rho_P}\autf(P)\bigr]$.  
In particular, $P$ is fully automized in $\calf$.

We claim that $P$ is also receptive.  Fix any $Q\in{}P^\calf$ and 
$\varphi\in\isof(Q,P)$.  By axiom (A2), there is $\psi\in\Iso_\calt(Q,P)$ 
such that $\rho_{Q,P}(\psi)=\varphi$.  Let $N_\varphi\le N_S(Q)$ be the 
subgroup of all $g\in{}N_S(Q)$ such that 
$\varphi{}c_g\varphi^{-1}\in\Aut_S(P)$.  Then 
$\psi\gee_Q(N_\varphi)\psi^{-1}\le\gee_P(N_S(P))\cdot{}E(P)$.  Since 
$\gee_P(C_S(P))\in\sylp{E(P)}$, 
$\gee_P(N_S(P))\in\sylp{\gee_P(N_S(P))\cdot{}E(P)}$.  So there is 
$\chi\in{}E(P)$ such that 
	\[ (\chi\psi)\gee_Q(N_\varphi)(\chi\psi)^{-1} \le \gee_P(N_S(P))~. 
	\]
Axiom (II) now implies that there is 
$\widebar{\psi}\in\Mor_\calt(N_\varphi,S)$ such that 
$\widebar{\psi}\circ\iota_{Q,N_\varphi}=\iota_{P,S}\circ\chi\psi$, and so 
$\rho(\widebar{\psi})\in\homf(N_\varphi,S)$ is an extension of 
$\varphi=\rho(\psi)=\rho(\chi\psi)$.  

It remains to prove that axiom (III) (for a fusion system) holds for all 
subgroups in $\calh$.  So assume $P_1\le P_2\le P_3\le\cdots$ are subgroups 
in $\calh$, $P=\bigcup_{i=1}^\infty{}P_i$, and $\varphi\in\Hom(P,S)$ is 
such that $\varphi_i\defeq\varphi|_{P_i}\in\homf(P_i,S)$ for each $i$.  For each $i$, 
$E(P_i)$ contains a discrete $p$-toral group with index prime to $p$, and 
thus is an extension of a discrete $p$-torus by a finite group.  Also, 
restriction defines a homomorphism from $E(P_{i+1})$ to $E(P_i)$, and this 
is injective by Proposition \ref{T-prop}(d).  We can thus regard the 
$E(P_i)$ as a decreasing sequence of discrete-$p$-toral-by-finite groups, 
and any such sequence becomes constant, since discrete $p$-toral groups are artinian (see \cite[Proposition 1.2]{BLO3}).  In other words, there is $N$ such 
that for all $j>i\ge N$, the restriction of $E(P_j)$ to $E(P_i)$ is an 
isomorphism.

For each $i\le{}N$, choose $\psi'_i\in\Hom_\calt(P_i,S)$ such that 
$\varphi_i=\rho_{P_i,S}(\psi'_i)$.  Set $\psi_N=\psi'_N$.  By (A2), 
for each $i>N$, there is $\chi_i\in{}E(P_N)$ such that 
$\psi_N=\psi'_i|_{P_N,S}\circ\chi_i$.  We just saw that there is 
$\widebar{\chi}_i\in{}E(P_i)$ such that $\widebar{\chi}_i|_{P_N}=\chi_i$.  
So if we set $\psi_i=\psi'_i\circ\widebar{\chi}_i$, then 
$\psi_i|_{P_N,S}=\psi'_i|_{P_N,S}\circ\chi_i=\psi_N$.  Since morphisms in 
$\calt$ are epimorphisms by Proposition \ref{T-prop}(d), 
$\psi_{i+1}|_{P_i,S}=\psi_i$ for each $i>N$, and by axiom (III) for $\calt$, 
there is $\psi\in\Mor_\calt(P,S)$ such that $\psi|_{P_i,S}=\psi_i$ for each 
$i\ge{}N$.  Then $\rho(\psi)|_{P_i,S}=\varphi_i$ for each $i$, so 
$\varphi=\rho(\psi)\in\homf(P,S)$.  
\end{proof}

We next show that linking systems over saturated fusion systems are 
transporter systems and characterize linking systems among transporter 
systems. This is a generalization to linking systems over discrete 
$p$-toral groups of \cite[Proposition 4]{link-ext}. We will first list the 
properties of linking systems that we need as a separate proposition that 
might be useful for future reference and then state the result as a 
corollary.

\begin{Prop} \label{L-prop}
The following hold for any linking system $\call$ associated to a saturated 
fusion system $\calf$ over a discrete $p$-toral group $S$.
\begin{enuma}  
\item\label{x:ext(A)}
 For each $P,Q\in\Ob(\call)$, the subgroup 
$E(P)\defeq\Ker\bigl[\Aut_\call(P)\Right2{}\autf(P)\bigr]$ acts freely on 
$\Mor_{\call}(P,Q)$ via right composition, and $\pi_{P,Q}$ induces a 
bijection 
	\[ \Mor_\call(P,Q)/E(P) \Right5{\cong} \homf(P,Q) ~. \]


\item \label{x:deltainj} The functor $\delta$ is injective on all morphism sets.

\item\label{x:leftaction}   The  action of $E(Q)$ 
on $\Hom_{\calf}(P,Q)$ via left composition is free.

\item \label{x:uniqrest}
For every morphism $\psi\in\Mor_\call(P,Q)$, and every 
$P_*,Q_*\in\Ob(\call)$ such that $P_*\le{}P$, $Q_*\le{}Q$, and 
$\pi(\psi)(P_*)\le{}Q_*$, there is a unique morphism 
$\psi|_{P_*,Q_*}\in\Mor_\call(P_*,Q_*)$ (the ``restriction'' of $\psi$) 
such that $\psi\circ\iota_{P_*,P}=\iota_{Q_*,Q}\circ\psi|_{P_*,Q_*}$.  

\item \label{x:fnorm}
If $P\in\Ob(\call)$ is fully normalized in $\calf$, then 
$\delta_P(N_S(P))\in\sylp{\Aut_\call(P)}$.

\item\label{x:uniqext} Let $P\nsg \widebar{P}\le S$ and 
$Q\le\widebar{Q}\le S$ be objects in $\call$.  Let 
$\psi\in\Mor_\call(P,Q)$ be such that for each $g\in\widebar{P}$, 
there is $h\in\widebar{Q}$ satisfying 
$\iota_{Q,\widebar{Q}}\circ\psi\circ\delta_P(g)
=\delta_{Q,\widebar{Q}}(h)\circ\psi$.   Then there is a unique morphism 
$\widebar{\psi}\in\Mor_\call(\widebar{P},\widebar{Q})$ such that 
$\widebar{\psi}|_{P,Q}=\psi$. 

\item\label{x:monoepi} All morphisms in $\call$ are monomorphisms and 
epimorphisms in the categorical sense.

\item\label{x:axiomIII} Assume $P_1\le P_2\le P_3\le\cdots$ in $\Ob(\call)$ and 
$\psi_i\in\Mor_\call(P_i,S)$ are such that for all $i\ge1$, 
$\psi_i=\psi_{i+1}|_{P_i,S}$.  Set $P=\bigcup_{i=1}^\infty P_i$.  Then 
there is $\psi\in\Mor_\call(P,S)$ such that $\psi|_{P_i,S}=\psi_i$ for each 
$i$.
\end{enuma}
\end{Prop}

\begin{proof} Points (a) and (b) are exactly the same as points (a) and 
(c), respectively, in \cite[Proposition 4]{link-ext}, and the proofs go 
through unchanged. The proof of (c) is contained in that of  
\cite[Proposition 4(f)]{link-ext}, again with no modification necessary. We 
prove the remaining points.

\smallskip

\noindent\textbf{(d) }  This is a special case of \cite[Lemma 
4.3(a)]{BLO3} (which is Proposition \ref{T-prop}(a) for linking 
systems).

\noindent\textbf{(e) } 
For each $P\in\Ob(\call)$ which is fully centralized in $\calf$, 
	\[ \autf(P)\cong\Aut_\call(P)/\delta_P(C_S(P)) 
	\qquad\textup{and}\qquad
	\Aut_S(P)\cong N_S(P)/C_S(P)\,. \]
Hence $[\autf(P):\Aut_S(P)]=[\Aut_\call(P):\delta_P(N_S(P))]$, since 
$\delta_P$ is injective. If $P\in\Ob(\call)$ is fully 
normalized in $\calf$, then $\Aut_S(P)\in\sylp{\autf(P)}$, and so 
$\delta_P(N_S(P))\in\sylp{\Aut_\call(P)}$.

\smallskip

\noindent\textbf{(f) } The proof of existence of an extension 
$\widebar{\psi}$ is identical in our case to the proof of the corresponding 
statement \cite[Propostion 4(e)]{link-ext}. It remains to prove uniqueness. 

Assume 
$\widebar{\psi}_1,\widebar{\psi}_2\in\Iso_\call(\widebar{P},\widebar{Q})$ are 
two extensions of $\psi$.  We must show that 
$\widebar{\psi}_1=\widebar{\psi}_2$.  It suffices to do this when $\psi$ is 
an isomorphism, and also (after composing by an isomorphism, if necessary) 
when $\widebar{P}$ is fully centralized in $\calf$.  

Fix $x\in\widebar{P}$, set $y_i=\pi(\widebar{\psi}_i)(x)$, and consider the 
following two squares of morphisms in $\call$:
	\[ \vcenter{
	\xymatrix@C+3ex@!{
	\widebar{P} \ar[r]^{\widebar{\psi}_i} 
	\ar[d]_{\delta_{\widebar{P}}(x)} & 
	\widebar{Q} \ar[d]^{\delta_{\widebar{Q}}(y_i)} &&
	P \ar[r]^{\psi}_{\cong} \ar[d]_{\delta_{P}(x)} 
	& Q \ar[d]^{\delta_{Q}(y_i)} \\
	\widebar{P} \ar[r]^{\widebar{\psi}_i} & \widebar{Q} &&
	P \ar[r]^{\psi}_{\cong} & Q \\
	}} \]
The first square commutes by axiom (C), and the second square is defined to 
be a restriction of the first.  Note that $\delta_{P}(x)$ is the 
restriction of $\delta_{\widebar{P}}(x)$ since $\delta$ is a 
functor (and since $\iota_{U,V}=\delta_{U,V}(1)$ for all $U\le{}V$ in 
$\Ob(\call)$).  
Hence the second square commutes by the uniqueness of restriction morphisms 
(point (\ref{x:uniqrest})).  Thus $\delta_{Q}(y_1)= \psi\circ\delta_{P}(x)\circ\psi^{-1} 
=\delta_{Q}(y_2)$.  Since $\delta_{Q}$ is 
injective by (\ref{x:deltainj}), this shows that $y_1=y_2$.  Since this holds for all 
$x\in\widebar{P}$, $\pi(\widebar{\psi}_1)=\pi(\widebar{\psi}_2)$. 

By axiom (A) (and since $\widebar{P}$ is fully centralized), 
$\widebar{\psi}_2=\widebar{\psi}_1\circ\delta_{\widebar{P}}(a)$ for some 
$a\in{}C_S(\widebar{P})$.  Hence
	\[ \iota_{Q,\widebar{Q}}\circ\psi 
	= \widebar{\psi}_2|_{P,\widebar{Q}} 
	= \bigl(\widebar{\psi}_1\circ\delta_{\widebar{P}}(a)\bigr)
	\big|_{P,\widebar{Q}}
	= \widebar{\psi}_1\circ\delta_{P,\widebar{P}}(a)
	= \widebar{\psi}_1|_{P,\widebar{Q}}\circ\delta_P(a)
	= \iota_{Q,\widebar{Q}}\circ\psi\circ\delta_P(a),  \]
so $\delta_P(a)=1$ since $E(P)$ acts freely on $\Mor_\call(P,\widebar{Q})$ 
by (\ref{x:ext(A)}) (note that $a\in{}C_S(P)$).  Then $a=1$ since $\delta_P$ is injective 
by (\ref{x:deltainj}), and so $\widebar{\psi}_1=\widebar{\psi}_2$.

\smallskip

\noindent\textbf{(g) } By \cite[Lemma A.6]{OV1}, for each 
$\varphi\in\Mor(\call)$, $\varphi$ is an isomorphism in $\call$ if and only 
if $\pi(\varphi)$ is an isomorphism in $\calf$.  Hence by (d), each 
morphism in $\call$ is the composite of an isomorphism followed by an 
inclusion.  So it suffices to prove that inclusions are monomorphisms and 
epimorphisms.  That they are monomorphisms follows from the uniqueness of 
restrictions in (d).  

By the uniqueness of extensions in (f), any inclusion $\iota_{P,Q}$ for 
$P\nsg{}Q$ is an epimorphism.  Since a composite of epimorphisms is an 
epimorphism, this proves that $\iota_{P,Q}$ is an epimorphism if $P<Q$ with 
finite index.  

Assume this is not true in general.  Then there are 
$\psi_1,\psi_2\in\Iso_\call(P,Q)$ and $P_*<P$ such that $\psi_1\ne\psi_2$ 
and $\psi_1|_{P_*}=\psi_2|_{P_*}$.  Choose $P_*<P_1\le{}P$ which is minimal 
among subgroups of $P$ containing $P_*$ for which 
$\psi_1|_{P_1}\ne\psi_2|_{P_1}$ ($S$ is artinian by \cite[Proposition 
1.2]{BLO3}).  There is no proper subgroup $P_2<P_1$ of finite index which 
contains $P_*$ (otherwise $\psi_1|_{P_2}=\psi_2|_{P_2}$ and the result 
follows from (f)), so $P_1$ is the union of its proper subgroups which 
contain $P_*$.  In particular, $\pi(\psi_1)|_{P_1}=\pi(\psi_2)|_{P_1}$.  
The rest of the argument to show that $\psi_1|_{P_1}=\psi_2|_{P_1}$ (thus 
giving a contradition) goes through exactly as in the last paragraph of the 
proof of (\ref{x:uniqext}). (Note that we can easily arrange for $P$ to be 
fully centralized.)

\smallskip

\noindent\textbf{(h) } 
Assume $P_1\le P_2\le P_3\le\cdots$ in $\Ob(\call)$ 
and $\psi_i\in\Mor_\call(P_i,S)$ are such that for each $i\ge1$, 
$\psi_i=\psi_{i+1}|_{P_i,S}$.  Set $P=\bigcup_{i=1}^\infty P_i$, and set 
$\varphi_i=\pi(\psi_i)\in\homf(P_i,S)$ for each $i$.  Then 
$\varphi_{i+1}|_{P_i}=\varphi_i$ for each $i$, so the union of the 
$\varphi_i$ is a homomorphism $\varphi\in\Hom(P,S)$.  By axiom (III) in 
Definition \ref{sat.Frob.}, $\varphi\in\homf(P,S)$.

By Lemma \ref{dot-props}(d), there is $N$ such that $P_i\8=P\8$ for each 
$i\ge{}N$.  Let $Q_N\in(P_N)^\calf$ be fully centralized in $\calf$, fix 
$\nu\in\isof(P_N,Q_N)$, and set $Q=\nu\8(P)$.  Thus 
$E(Q_N)=\delta_{Q_N}(C_S(Q_N))$, each element of $C_S(Q_N)$ also centralizes 
$Q_N\8$ by the uniqueness of extensions to $Q_N\8$, so each element of 
$E(Q_N)$ extends to an element of $E(Q)$, which implies that each element 
of $E(P_N)$ extends to an element of $E(P)$.  

Choose $\psi'\in\Mor_\call(P,S)$ such that $\pi(\psi')=\varphi$.  By axiom 
(A), there is $\chi_N\in{}E(P_N)$ such that 
$\psi_N=\psi'|_{P_N,S}\circ\chi_N$.  Let $\chi\in{}E(P)$ be such that 
$\chi|_{P_N}=\chi_N$, and set $\psi=\psi'\circ\chi$.  Then 
$\psi|_{P_N,S}=\psi'|_{P_N,S}\circ\chi_N=\psi_N$, and 
$\psi|_{P_i,S}=\psi_i$ for each $i$:  via composition with inclusions when 
$i\le N$, and by \eqref{x:monoepi} (morphisms in $\call$ are 
epimorphisms) when $i>N$. 
\end{proof}

We are now ready to show that all linking systems are transporter 
systems.

\begin{Cor}\label{L-Cor}
Fix a saturated fusion system $\calf$ over a discrete $p$-toral group 
$S$. Then each linking system $\call$ associated to $\calf$ is also a 
transporter system. Conversely, a transporter system $\calt$ associated to 
$\calf$ which contains all $\calf$-centric $\calf$-radical subgroups as 
objects, and such that the kernel subgroups $E(P)$ are all discrete 
$p$-toral groups, is a linking system. 
\end{Cor}

\begin{proof} Assume $\call$ is a linking system associated to 
$\calf$.  Axiom (A2) in Definition \ref{tport} follows from 
Proposition \ref{L-prop}(\ref{x:ext(A)},\ref{x:leftaction}), (B) and (C) 
follow from the corresponding axioms in Definition \ref{L-cat}, and (I), 
(II), and (III) follow from points (\ref{x:fnorm}), (\ref{x:uniqext}), and 
(\ref{x:axiomIII}), respectively, in Proposition \ref{L-prop}.  Thus 
$\call$ is a transporter system. 

Assume now that $\calt$ is a transporter system associated to $\calf$ such 
that $\Ob(\calt)$ contains all $\calf$-centric $\calf$-radical subgroups of 
$S$. By axiom (A2), $\rho_{P,Q}$ sends $\Mor_\calt(P,Q)$ surjectively onto 
$\homf(P,Q)$ for each $P,Q\in\Ob(\calt)$, so every object of $\calt$ is 
isomorphic in $\calt$ to an object which is fully centralized.  If $P$ is 
fully centralized, then $\gee_P$ sends  $C_S(P)$ injectively to a Sylow 
$p$-subgroup of $E(P)$, and hence $E(P)=\gee_P(C_S(P))$ in that case since 
we are assuming that $E(P)$ is discrete $p$-toral. It follows that $\calt$ 
is a linking system associated to $\calf$.  
\end{proof}

\begin{Prop} \label{centric-E(P)}
Let $\calf$ be a saturated fusion system over a discrete $p$-toral group 
$S$.  Let $\calt$ be a transporter system associated to $\calf$, all of 
whose objects are $\calf$-centric, which contains all $\calf$-centric 
$\calf$-radical subgroups of $S$.  Then for every $P\in\Ob(\calt)$, 
$E(P)\cong E_0(P)\times Z(P)$ for some $E_0(P)\nsg E(P)$ which is finite of 
order prime to $p$ and normal in $\Aut_\calt(P)$.  There is a linking 
system $\call$ associated to $\calf$, defined by setting 
$\Ob(\call)=\Ob(\calt)$ and 
	\[ \Mor_\call(P,Q) = \Mor_\calt(P,Q)/E_0(P)\,, \] 
and the natural functor $\calt\Onto2{}\call$ induces a mod $p$ equivalence 
between the geometric realizations and hence a homotopy equivalence 
$|\calt|\pcom\simeq|\call|\pcom$.
\end{Prop}

\begin{proof}  By axiom (C), for all $P\in\Ob(\calt)$, $E(P)$ commutes with 
$\gee_P(P)$ in the group $\Aut_\calt(P)$.  Hence 
$\gee_P(Z(P))=\gee_P(P)\cap{}E(P)$ is central in $E(P)$, and it has finite 
index prime to $p$ by axiom (I) and since $P$ is $\calf$-centric.  So 
$E(P)\cong{}Z(P)\times{}E_0(P)$ for some unique $E_0(P)\nsg E(P)$ which is 
finite of order prime to $p$.  

It is now straightforward to check that $\call$, when defined as above, is 
a quotient category of $\calt$ (i.e., composition is well defined).  
Axioms (A1), (A2), (B), and (C) for a transporter system imply that $\call$ 
satisfies the corresponding axioms for a linking system, and thus is a 
linking system associated to $\calf$.  The induced map of spaces 
$|\calt|\Right2{}|\call|$ is a mod $p$ equivalence (and hence 
$|\calt|\pcom\simeq|\call|\pcom$) by \cite[Lemma 1.3]{BLO1}.
\end{proof}

We next define normal transporter subsystems.  Recall that for any fusion 
system $\calf$ over a discrete $p$-toral group $S$, a normal subgroup 
$\dec{S}\nsg S$ is \emph{weakly closed} in $\calf$ if for each $g\in\dec{S}$, 
$g^\calf\subseteq\dec{S}$.  In other words, no element of $\dec{S}$ is 
$\calf$-conjugate to any element of $S{\sminus}\dec{S}$.

\begin{Defi} \label{d:T0<|T}
Fix a pair of saturated fusion systems $\dec\calf\subseteq\calf$ over 
discrete $p$-toral groups $\dec{S}\nsg{}S$, and let $\dec\calt\subseteq\calt$ be 
associated transporter systems.  Then $\dec\calt$ is \emph{normal} in $\calt$ 
($\dec\calt\nsg\calt$) if
\begin{enumerate}[\rm(i) ]
\item $\dec{S}$ is strongly closed in $\calf$ and 
$\Ob(\dec\calt)=\{P\cap \dec{S} \,|\, P\in\Ob(\calt) \}$;

\item for all $P,Q\in\Ob(\dec\calt)$ and $\psi\in\Mor_{\calt}(P,Q)$,
there are morphisms $\gamma\in\Aut_\calt(\dec{S})$ and
$\psi_*\in\Mor_{\dec\calt}(\gamma(P),Q)$ such that
$\psi=\psi_*\circ\gamma|_{P,\gamma(P)}$; and

\item for all $P,Q\in\Ob(\dec\calt)$, $\psi\in\Mor_{\dec\calt}(P,Q)$, and 
$\gamma\in\Aut_\calt(\dec{S})$, $\gamma|_{Q,\gamma(Q)}\circ\psi\circ 
\gamma|_{P,\gamma(P)}^{-1}$ is in $\Mor_{\dec\calt}(\gamma(P),\gamma(Q))$.

\end{enumerate}
Here, in (ii) and (iii), we write $\gamma(P)=\rho(\gamma)(P)$ and
$\gamma(Q)=\rho(\gamma)(Q)$ for short.  In this situation, we define
    \[ \calt/\dec\calt = \Aut_\calt(\dec{S})/\Aut_{\dec\calt}(\dec{S})\,. \]
\end{Defi}

Let $\calt$ be a transporter system.  A subgroup $P\in\Ob(\calt)$ is 
defined to be \emph{$\calt$-radical} if $O_p(\Aut_\calt(P))=\gee_P(P)$.  
When $\calt$ is a linking system associated to $\calf$, then $P$ is 
$\calt$-radical if and only if it is $\calf$-centric and $\calf$-radical.  
If $\calt=\calt_S(G)$ for a finite group $G$ and $S\in\sylp{G}$, then 
$P\le{}S$ is $\calt$-radical if and only if it is a radical $p$-subgroup of 
$G$ in the usual sense.

\newcommand{\FpFun}{\underline{\F}_p}

\begin{Lem} \label{drop1}
Let $\calf$ be a saturated fusion system over a discrete $p$-toral group 
$S$, and let $\calt$ be a transporter system associated to $\calf$.  Let 
$P\in\Ob(\calt)$ be a minimal object, let $\calp$ be the $\calf$-conjugacy 
class of $P$, and let $\calt_0\subseteq\calt$ the full subcategory with 
object set $\Ob(\calt){\sminus}\calp$.  Assume that 
\begin{itemize} 
\item $P$ is fully normalized, 
\item $N_S(P)/P$ is finite, and 
\item $P$ is either not $\calf$-centric or not $\calt$-radical.  
\end{itemize}
Then the inclusion of $|\calt_0|$ into $|\calt|$ induces an isomorphism in 
mod $p$ cohomology. 
\end{Lem}

\begin{proof} Let $\Phi\:\calt\op\Right2{}\zploc\mod$ be the functor which 
sends objects in $\calp$ to $\F_p$ and other objects to $0$, regarded as a 
subfunctor of the constant functor $\FpFun$.  By the minimality of $P$, and 
since the quotient functor $\FpFun/\Phi$ vanishes on $\calp$, it follows that 
	\[H^*(\calt, \FpFun/\Phi) \cong H^*(\calt_0,\FpFun) 
	\cong H^*(|\calt_0|, \F_p),\] 
where the second isomorphism holds  because $\FpFun$ is constant on 
$\calt_0$. Upon applying $H^*(\calt, -)$ to the short exact sequence of 
functors associated to the inclusion of $\Phi$ in $\FpFun$, and using 
the above isomorphism, we conclude that 
$H^*(\calt, \Phi) \cong H^*(|\calt|,|\calt_0|,\F_p)$.

Let $\calo(\calt)$ denote the orbit category associated to $\calt$: 
the category with the same objects as $\calt$, and with morphism sets 
$\Mor_{\calo(\calt)}(P,Q) = \Mor_\calt(P,Q)/Q$, where $Q\cong\gee_Q(Q)$ acts 
freely by left composition. The projection 
$\calt\to\calo(\calt)$ is target regular in the sense of \cite[Definition 
A.5(b)]{OV1}, so by \cite[Proposition A.11]{OV1}, there is a spectral sequence
	\[ E_2^{ij} = H^i\bigl(\calo(\calt);H^j(-;\Phi(-))\bigr)
	\quad\implies\quad H^{i+j}(\calt; \Phi)\,. \]

Since $N_S(P)/P$ is finite, and $[\Aut_\calt(P):\gee_P(N_S(P))]<\infty$ by axiom
(I) (recall $P$ is fully normalized), $\Aut_\calt(P)/\gee_P(P)$ is also 
finite.  By \cite[Proposition 5.4]{BLO3}, and since $\Phi(Q)=0$ for 
$Q\notin\calp$,
	\begin{align*} 
	H^i\bigl(\calo(\calt);  H^j(-;\Phi(-)) \bigr) &\cong 
	\Lambda^i\bigl(\Aut_\calt(P)/\gee_P(P); H^j(P;\Phi(P))\bigl) \\
	&\cong \Lambda^i\bigl(\Aut_\calt(P)/\gee_P(P);H^j(P;\FpFun)\bigr) 
	\end{align*}
for each $i$ and $j$.  If $P$ is not $\calf$-centric (and since it is fully 
normalized), there is $g\in{}C_S(P){\sminus}P$, and the class of 
$\gee_P(g)\ne1$ in $\Aut_\calt(P)/\gee_P(P)$ acts trivially on 
$H^j(P;\FpFun)$.  If $P$ is not $\calt$-radical, then 
$O_p(\Aut_\calt(P)/\gee_P(P))\ne1$.  In either case, by \cite[Proposition 
6.1(ii)]{JMO}, 
$\Lambda^i\bigl(\Aut_\calt(P)/\gee_P(P);H^j(P;\FpFun)\bigr)=0$ for each $i$ 
and $j$, so 
	\[H^*(|\calt|,|\calt_0|;\F_p)\cong H^*(\calt;\Phi) =0.\] 
Thus the inclusion $|\calt_0|\to|\calt|$ induces an isomorphism on 
mod $p$ cohomology as claimed.
\end{proof}

Lemma \ref{drop1} gives us a tool for proving that under certain 
conditions, adding or subtracting one conjugacy class in a transporter 
system does not change the (mod $p$) homotopy type of its realization.  
However, more is needed when dealing with infinitely many conjugacy 
classes. One problem when doing this is that the ``bullet functor'' 
$P\mapsto{}P\8$ need not lift to a functor on a transporter system.  This 
is why we need to assume the existence of a normal linking subsystem in the 
following proposition.

\begin{Prop} \label{dropmore}
Let $\calf$ be a saturated fusion system over a discrete $p$-toral group 
$S$, and let $\calt$ be a transporter system associated to $\calf$. 
Let $\calt_0\subseteq\calt$ be a 
full subcategory which is a transporter subsystem associated to a full 
subcategory $\calf_0\subseteq\calf$.  Assume there is a normal linking 
subsystem $\dec\call\nsg\calt$ over a subgroup $\dec{S}\le S$ of finite 
index, and let $\dec\calf\subseteq\calf$ be the corresponding fusion 
subsystem. Assume also that either 
\begin{enuma} 
\item all objects in $\dec\call$ are $\dec\calf$-centric and 
$\Ob(\calt_0)\supseteq\Ob(\calt)\cap\Ob(\calf^c)$; or 

\item each object in $\Ob(\calt){\sminus}\Ob(\calt_0)$ is $\calf$-centric 
and not $\calt$-radical.

\end{enuma}
Then the inclusion of $|\calt_0|$ in $|\calt|$ is a mod $p$ homology 
equivalence.
\end{Prop}

\begin{proof} For each $P\le{}S$, we write $\dec{P}=P\cap{}\dec{S}$, and 
set $\6P=P\dec{P}\8$ (where $(-)\8$ means the bullet construction for 
$\calf$).  By Lemma \ref{dot-props}(b), $N_S(\dec{P})\le N_S(\dec{P}\8)$,  
and since $\dec{P}\nsg P$, it follows that  $P$ normalizes $\dec{P}\8$ (so 
$\6P$ is a subgroup).  Also, since $\dec{P}\8\le\dec{S}$, 
$\dec{\6P}=\dec{P}\cdot\dec{P}\8=\dec{P}\8$, and so $\6{(\6P)} = \6P$.

For each $P,Q\in\Ob(\calt)$ and each $\psi\in\Mor_\calt(P,Q)$, let 
$\dec{\psi}=\psi|_{\dec{P},\dec{Q}}$.  Then $\dec\psi$ extends to 
$\psi\in\Mor_\calt(P,Q)$ and to 
$\dec\psi\8\in\Mor_\calt(\dec{P}\8,\dec{Q}\8)$ (Lemma 
\ref{reduce2dots}(\ref{4.6.a})), and hence by axiom (II) to 
$\6\psi\in\Mor_\calt(\6P,\6Q)$. By the uniqueness of these extensions, 
$\6\psi|_P=\psi$, and thus $\6{(-)}$ defines an idempotent functor from 
$\calt$ to itself. 

Let $\6\calt\subseteq\calt$ be the full subcategory whose objects are those 
$P$ such that $\6P=P$ (and use the same notation for the corresponding 
subcategory of any sub-transporter system of $\calt$). For each $P\le 
S$, $\6P/\dec{P}\8$ is a subgroup of order at most $|S/\dec{S}|$ in 
$N_S(\dec{P}\8)/\dec{P}\8$, and by \cite[Lemma 1.4(a)]{BLO3}, there are 
only finitely many conjugacy classes of such subgroups.  Since $\calf\8$ 
has finitely many $S$-conjugacy classes of objects, $\6\calt$ also has 
finitely many $S$-conjugacy classes of objects.

Let $\mathfrak{H}$ be the set of all transporter subsystems 
$\calt'\subseteq\calt$ such that 
\begin{itemize} 
\item $\calt'\supseteq\calt_0$, 
\item $\calt'$ is a full subcategory of $\calt$, and 
\item the inclusion $|\calt_0|\subseteq|\calt'|$ is a mod $p$ homology 
equivalence in case (a), or a homotopy equivalence in case (b).  
\end{itemize}
We must show that $\calt\in\mathfrak{H}$.

Assume otherwise, and choose $\calt_1\in\mathfrak{H}$ for which 
$\Ob(\calt_1)\cap\Ob(\6\calt)$ contains the largest possible number of 
$\calf$-conjugacy classes.  Let $\calt_2\subseteq\calt$ be the full 
subcategory with $\Ob(\calt_2)=\{P\in\Ob(\calt)\,|\,\6P\in\Ob(\calt_1)\}$.  
By the above discussion, $|\6{\calt_1}|=|\6{\calt_2}|$ is a strong deformation 
retract of $|\calt_1|$ and of $|\calt_2|$, so $|\calt_1|\simeq|\calt_2|$, 
and $\calt_2\in\mathfrak{H}$.  

Since $\calt_2\subsetneqq\calt$ by assumption, 
$\Ob(\6\calt)\not\subseteq\Ob(\calt_2)$. Let $P$ be maximal among objects 
in $\6\calt$ not in $\calt_2$.  By definition of $\Ob(\calt_2)$, $P$ is 
maximal among all objects of $\calt$ not in $\calt_2$.  Let 
$\calt_3\subseteq\calt$ be the full subcategory with 
$\Ob(\calt_3)=\Ob(\calt_2)\cup{}P^\calf$.  

Now, $N_S(P)/PC_S(P)\cong\Out_S(P)$ is finite by \cite[Lemma 2.5]{BLO3}.  
In case (a), $PC_S(P)/P$ is finite since $\dec{P}$ is $\dec\calf$-centric 
and $[S:\dec{S}]<\infty$, while $C_S(P)\le P$ in case (b).  So 
$|N_S(P)/P|<\infty$ in either case.

By assumption, $P$ is either not $\calf$-centric (case (a)) or not 
$\calt$-radical (case (b)).  So by Lemma \ref{drop1}, the inclusion of 
$|\calt_2|$ into $|\calt_3|$ is a mod $p$ equivalence.  Hence 
$\calt_3\in\mathfrak{H}$, contradicting our maximality assumption on 
$\calt_1$.  
\end{proof}

\begin{Cor} \label{L0<L}
Let $\SFL$ be a linking triple, and let $\call_0\subseteq\call$ be a full 
subcategory which is a linking subsystem associated to some 
$\calf_0\subseteq\calf$.  Then the inclusion of $|\call_0|$ into $|\call|$ 
is a mod $p$ homology equivalence.
\end{Cor}

\begin{proof} Since $\call$ is a linking system, a subgroup $P\le S$ is 
$\call$-radical exactly when it is $\calf$-centric and $\calf$-radical. So 
by definition, $\call_0$ contains all $\call$-radical subgroups.  Hence the 
lemma is a special case of Proposition \ref{dropmore}(b).  
\end{proof}

It remains to look at the problem of constructing a centric linking system 
which contains as full subcategory a linking system over a smaller set of 
objects.  As usual, we first check what happens when we add one conjugacy 
class of objects.

\begin{Lem} \label{existsL1}
Let $\calf$ be a saturated fusion system over a discrete $p$-toral group 
$S$, and let $\calf^c$ be the full subcategory of $\calf$-centric objects. 
Let $\calf_0\subseteq\calf_1\subseteq\calf^c$ be full subcategories such 
that $\Ob(\calf_0)$ and $\Ob(\calf_1)$ are invariant under 
$\calf$-conjugacy and closed under overgroups, $\calf_0$ contains all 
$\calf$-centric $\calf$-radical subgroups of $S$, and 
$\Ob(\calf_1)=\Ob(\calf_0)\cup\calp$ for some $\calf$-conjugacy class 
$\calp$. If $\call_0$ is a linking system associated to $\calf_0$, then 
there is a linking system $\call_1$ associated to $\calf_1$ such that 
$\call_0$ is isomorphic as a linking system to the full subcategory of 
$\call_1$ with same set of objects as $\call_0$.
\end{Lem}

\begin{proof} This follows by the same proof as in \cite[Proposition 
III.4.8]{AKO} (Steps 1--3), together with Lemma \ref{l:centric}.  

Alternatively, let $\Phi$ be the functor: $\Phi(P)=Z(P)$ when $P\in\calp$, 
and $\Phi(P)=0$ when $P\in\Ob(\calf_0)$.  For $P\in\calp$, 
$\Lambda^*(\outf(P);\Phi(P))=0$ by \cite[Proposition 6.1(ii)]{JMO} and since 
$O_p(\outf(P))\ne1$.  Hence by \cite[Proposition 5.4]{BLO3}, 
$H^3(\calo(\calf_1);\Phi)=0$. By an argument similar to that used to prove 
\cite[Proposition 3.1]{BLO2}, the obstruction to extending $\call_0$ to 
$\call_1$ lies in $H^3(\calo(\calf_1);\Phi)$, and hence $\call_1$ does 
exist.
\end{proof}

This is now generalized as follows.

\begin{Prop} \label{existsLc}
Let $\calf$ be a saturated fusion system over a discrete $p$-toral group 
$S$, and let $\calf^c$ be the full subcategory of $\calf$-centric objects. 
Let $\calf_0\subseteq \calf^c$ be a full subcategory such that 
$\Ob(\calf_0)$ is invariant under $\calf$-conjugacy and closed under 
overgroups, and contains all $\calf$-centric $\calf$-radical 
subgroups of $S$. If $\call_0$ is a linking system associated to $\calf_0$, 
then there is a centric linking system $\call$ associated to $\calf$ 
(associated to $\calf^c$) which contains $\call_0$ as a linking subsystem.
\end{Prop}

\begin{proof} Let $\mathfrak{H}$ be the set of all families of subgroups 
$\calh\subseteq\Ob(\calf^c)$ such that 
\begin{itemize} 
\item $\calh\supseteq\Ob(\calf_0)$, 
\item $\calh$ is invariant under $\calf$-conjugacy and closed under 
overgroups, and 
\item there is a linking system with object set $\calh$ 
which contains $\call_0$ as a linking subsystem.  
\end{itemize}
We must show that $\Ob(\calf^c)\in\mathfrak{H}$.

Choose $\calh_1\in\mathfrak{H}$ for which 
$\calh_1\cap\calf\8$ contains the largest possible number of 
$\calf$-conjugacy classes.  Set $\calh_2=\{P\in\Ob(\calf^c)\,|\,P\8\in\calh_1\}$, 
and let $\calf_1\subseteq\calf_2\subseteq\calf$ be the full subcategories 
with $\Ob(\calf_i)=\calh_i$.  Let $\call_2$ be the pullback of $\call_1$ 
and $\calf_2$ via the functors
	\[ \call_1 \Right5{\pi_1} \calf_1 \Left5{(-)\8} \calf_2\,. \]
Then $\call_2$ is a linking system associated to $\calf_2$ (recall that 
$Z(P)=Z(P^\bullet)$), 
and hence $\calh_2\in\mathfrak{H}$.  

Assume $\Ob((\calf^c)\8){\sminus}\calh_1\ne\emptyset$, and let 
$P$ be maximal among subgroups in this set.  Then $P$ is 
maximal among subgroups in $\Ob(\calf^c){\sminus}\calh_2$, so $\calh_2\cup 
P^\calf\in\mathfrak{H}$ by Lemma \ref{existsL1}, and this contradicts the 
choice of $\calh_1$. Hence $\Ob((\calf^c)\8)\subseteq\calh_1$, 
and so $\calh_2=\Ob(\calf^c)\in\mathfrak{H}$. 
\end{proof}


\bigskip\bigskip

\begin{thebibliography}{BLO2}

\bibitem[AOV]{AOV1} K. Andersen, B. Oliver, \& J. Ventura, Reduced, tame,
and exotic fusion systems, Proc. London Math. Soc. 105 (2012), 87--152

\bibitem[AKO]{AKO} M. Aschbacher, R. Kessar, \& B. Oliver, Fusion systems 
in algebra and topology, Cambridge Univ. Press (2011)

\bibitem[BGM]{BGM} M. Barratt, V. Gugenheim, \& J. Moore, On semisimplicial 
fibre-bundles, Amer. J. Math. 81 (1959), 639--657

\bibitem[BK]{BK} A. K. Bousfield \& D. M. Kan,  Homotopy limits, 
completions and localizations. Lecture Notes in Mathematics, Vol. 304. 
Springer-Verlag, Berlin-New York, 1972. v+348 pp.

\bibitem[BCGLO1]{BCGLO1} C. Broto, N. Castellana, J. Grodal, R. Levi, \&
B. Oliver, Subgroup families controlling $p$-local finite groups, Proc.
London Math. Soc. 91 (2005), 325--354

\bibitem[BCGLO2]{BCGLO2} C. Broto, N. Castellana, J. Grodal, R. Levi, \&
B. Oliver, Extensions of $p$-local finite groups, Trans. Amer. Math. Soc.
359 (2007), 3791--3858

\bibitem[BLO1]{BLO1} C. Broto, R. Levi, \& B. Oliver, Homotopy equivalences
of $p$-completed classifying spaces of finite groups, Invent. math. 151
(2003), 611--664

\bibitem[BLO2]{BLO2} C. Broto, R. Levi, \& B. Oliver, The homotopy theory
of fusion systems, J. Amer. Math. Soc. 16 (2003), 779--856

\bibitem[BLO3]{BLO3} C. Broto, R. Levi, \& B. Oliver, Discrete models for 
the $p$-local homotopy theory of compact Lie groups and $p$-compact groups, 
Geometry \& Topology 11 (2007) 315--427

\bibitem[BLO4]{BLO4} C. Broto, R. Levi, \& B. Oliver, A geometric
construction of saturated fusion systems, An alpine anthology of homotopy
theory (proceedings Arolla 2004), Contemp. math. 399 (2006), 11-39

\bibitem[Cu]{Curtis}  E. Curtis, Simplicial homotopy theory, Adv. in 
Math. 6 (1971), 107--209

\bibitem[DW]{DW} W. G. Dwyer, \& C. W. Wilkerson,  Homotopy fixed-point 
methods for Lie groups and finite loop spaces. Ann. of Math. (2) 139 
(1994), no. 2, 395--442.

\bibitem[GZ]{GZ} P. Gabriel \& M. Zisman, Calculus of fractions and 
homotopy theory, Springer-Verlag (1967)

\bibitem[GJ]{GJ} P. Goerss \& R. Jardine, Simplicial Homotopy Theory, 
Birkh\"auser Verlag (1999)

\bibitem[G]{Gorenstein} D. Gorenstein, Finite groups, Harper \& Row (1968)

\bibitem[JMO]{JMO} S. Jackowski, J. McClure, \& B. Oliver, Homotopy 
classification of self-maps of $BG$ via $G$-actions, \emph{Annals of 
Math.}, {\bf135} (1992), 183--270.


\bibitem[O1]{link-ext} B. Oliver, Extensions of linking systems and fusion
systems, Trans. Amer. Math. Soc. 362 (2010), 5483--5500

\bibitem[OV1]{OV1} B. Oliver, \& J. Ventura, Extensions of linking systems
with $p$-group kernel, Math. Annalen 338 (2007), 983--1043

\bibitem[Q]{Qu2} D. Quillen, Algebraic $K$-theory I, Lecture notes in 
mathematics 341 (1973), 77--139

\bibitem[RS]{RS} K.~Roberts, \& S.~Shpectorov, On the definition of saturated 
fusion systems, J.~Group Theory 12 (2009), 679--687.

\end{thebibliography}
\end{document}